\documentclass[a4paper,10pt]{article}
\usepackage[a4paper]{geometry}
\usepackage{amsmath}
\usepackage[latin1]{inputenc}
\usepackage[T1]{fontenc}
\allowdisplaybreaks[2]
\usepackage{graphicx}
\usepackage{amssymb}
\usepackage{amsthm}
\def\var{\mathrm{Var}}

\def \Lam{\Lambda(\lambda,J_0^c)}
\def\argmin{\mathrm{argmin}}

\newcommand{\indic}{\mathbf{1}}
\newtheorem{Th}{Theorem}

\newtheorem{Prop}{Proposition}
\newtheorem{Lemma}{Lemma}

\newtheorem{Cor}{Corollary}
\renewenvironment{proof}{\noindent{\bf Proof.}}{\hfill
  $\blacksquare$\par\noindent}
 \newcommand{\com}[1]{}
\newcommand{\norm}[1]{\ensuremath{\vert\!\vert #1 \vert\!\vert}}
\newcommand{\normu}[1]{\ensuremath{\vert\!\vert #1 \vert\!\vert}_{\ell_1}}
\newcommand{\normd}[1]{\ensuremath{\vert\!\vert #1 \vert\!\vert}_{\ell_2}}
\newcommand{\norminf}[1]{\ensuremath{\vert\!\vert #1 \vert\!\vert}_{\ell_\infty}}
\newcommand{\normD}[1]{\ensuremath{\vert\!\vert #1 \vert\!\vert}_{2}}
\newcommand{\E}{\ensuremath{\mathbb{E}}}

\renewcommand{\P}{\ensuremath{\mathbb{P}}}

\newcommand{\R}{\ensuremath{\mathbb{R}}}

\renewcommand{\ln}{{\log}}
\renewcommand{\L}{\ensuremath{\mathbb{L}}}

\newcommand{\ga}{\ensuremath{\gamma}}
\newcommand{\al}{\ensuremath{\alpha}}
\newcommand{\la}{\ensuremath{\lambda}}
\newcommand{\si}{\ensuremath{\sigma}}
\newcommand{\e}{\ensuremath{\varepsilon}}
\newcommand{\p}{\ensuremath{\varphi}}

\newcommand{\be}{\ensuremath{\beta}}

\newcommand{\hb}{\ensuremath{\hat{\beta}}}

\newcommand{\Ne}{\ensuremath{\mathbb{N}}}
\newcommand{\fo}{\ensuremath{f}_0}
\newcommand{\Jo}{\ensuremath{J_0}}
\newcommand{\JoUn}{\ensuremath{J_{01}}}
\newcommand{\JUn}{\ensuremath{J_{1}}}
\newcommand{\Comp}[1]{\ensuremath{#1^C}}
\newcommand{\JoC}{\ensuremath{\Comp{\Jo}}}

\newcommand{\etam}{\eta^{(-)}}
\newcommand{\etaM}{\eta^{(+)}}

\newcommand{\parinc}[2]{\parbox[c]{#1}{\includegraphics[width=#1]{#2}}}

\graphicspath{{./fig/}}

\title{Adaptive Dantzig density estimation}
\makeatletter
\def\timenow{\@tempcnta\time
  \@tempcntb\@tempcnta
  \divide\@tempcntb60
  \ifnum10>\@tempcntb0\fi\number\@tempcntb
  \multiply\@tempcntb60
  \advance\@tempcnta-\@tempcntb
  \ifnum10>\@tempcnta0\fi\number\@tempcnta}
\makeatother
\author{K. \textsc{Bertin}\thanks{Supported by Project PBCT 13
    laboratorio ANESTOC and Project FONDECYT 1090285. Departamento de
    Estad\'{i}stica, CIMFAV, Universidad de Valpara\'{i}so, Avenida
    Gran Breta\~na 1091, Valpara\'{i}so, Chile. Tel
    0056-(0)32-2508324. Email: karine.bertin@uv.cl}, E. \textsc{Le
    Pennec}\thanks{Laboratoire de Probabilit\'e et Mod\`eles
    Al\'eatoires, Universit\'e Paris 7, 175 rue de Chevaleret, F-75013
    Paris, France. Email: lepennec@math.jussieu.fr},
  V. \textsc{Rivoirard}\thanks{
Laboratoire de Math\'ematique,
 U.M.R. C.N.R.S. 8628,
Universit\'e Paris Sud,
91405 Orsay Cedex,
France and
D\'epartement de Math\'ematiques et Applications,
U.M.R. C.N.R.S. 8553,
ENS-Paris,
45 Rue d'Ulm,
75230 Paris Cedex 05,
France, Email: vincent.rivoirard@math.u-psud.fr}}

\date{}
\begin{document}

\maketitle

\begin{abstract}
This paper deals with the problem of density estimation. We aim at building an estimate of an unknown density as a linear combination of
functions of a dictionary. Inspired by Cand\`es
and Tao's approach, we propose an $\ell_1$-minimization under an
adaptive Dantzig constraint coming from sharp concentration
inequalities. This allows to consider a wide class of 
dictionaries. Under local or global coherence assumptions, oracle inequalities are derived. These theoretical results are also proved to be valid for the natural Lasso estimate associated with our Dantzig procedure. Then, the issue of calibrating these procedures is studied from both theoretical and practical points of view. Finally, a numerical study shows the significant improvement obtained by our procedures when compared with other classical procedures.
\end{abstract}



\section{Introduction}
Various estimation procedures based on $l_1$ penalization (exemplified by the Dantzig procedure
in \cite{candes} and the LASSO procedure in \cite{tib}) have extensively been studied recently. These procedures are computationally efficient as shown in
\cite{efron,osborne1,osborne2},
and thus are adapted
to high-dimensional  data. They have been widely used in regression
models, but only the Lasso estimator has
been studied in the density model (see \cite{bunea3,bunea7,geer}).
Although we will mostly consider the Dantzig estimator in the density model
for which no result exists so far, we recall some of the classical results
obtained in different settings by procedures based on $l_1$ penalization.

The Dantzig selector has been introduced by Cand\`es and Tao
\cite{candes} in the linear regression model.
More precisely,
given
\begin{equation*}\label{model}
Y=A \lambda_0 + \e, \end{equation*} where $Y\in\R^n$, $A$ is a $n$
by $M$ matrix, $\e\in\R^n$ is the noise vector and
$\lambda_0\in\R^M$ is  the unknown  regression parameter  to
estimate,  the Dantzig  estimator is defined by
\begin{equation*}
\hat{\lambda}^D=\arg\min_{\lambda\in\R^M}\normu{\lambda} \text{ subject
to } \norminf{A^T(A\lambda-Y)} \le \eta,
\end{equation*}
where $\norminf{\cdot}$ is the sup-norm in $\R^M$, $\normu{\cdot}$ is the $\ell_1$ norm in
$\R^M$, and $\eta$ is a regularization parameter.
A natural companion of this estimator is the Lasso procedure or
more precisely its relaxed form
\begin{equation*}\label{estlasso}
\hat{\lambda}^L=\arg\min_{\lambda\in\R^M}\left\{\frac{1}{2}\normd{A\lambda-Y}^2\ + \eta \normu{\lambda}\right\},
\end{equation*}
where $\eta$ plays exactly the exact same role as for the Dantzig estimator.
This $\ell_1$ penalized method
is also called \emph{basis pursuit} in signal processing (see \cite{chen,don}).

Cand\`es and Tao \cite{candes}
have obtained a bound for the $\ell_2$ risk of the estimator
$\hat{\lambda}^D$, with large probability, under a global condition on
the matrix $A$ (the Restricted Isometry Property) and a sparsity
assumption on $\lambda_0$, even for $M\geq n$. Bickel et al. \cite{paralleling}
have obtained oracle inequalities and bounds of the $\ell_p$ loss for both estimators under weaker
assumptions. Actually, Bickel et al. \cite{paralleling} deal
with the non parametric regression framework
 in which one observes
\[
Y_i = f(x_i) + e_i,\quad i=1,\ldots,n
\]
where $f$ is an unknown function while $(x_i)_{i=1,\ldots,n}$ are known design points
and $(e_i)_{i=1,\ldots,n}$ is a noise vector. There is no intrinsic matrix
$A$ in this problem but for any dictionary of functions $\Upsilon=(
\p_m)_{m=1,\ldots,M}$
 one can search $f$ as a weighted sum $f_\lambda$ of elements
of $\Upsilon$
\[
f_\lambda = \sum_{m=1}^M \lambda_m \p_m
\]
and introduce the matrix  $A=(\p_m(x_i))_{i,m}$,
 which summarizes the information on the dictionary
and on the design. Notice that if there exists $\lambda_0$ such that
$f=f_{\lambda_0}$ then the model can be rewritten exactly as the
classical linear model. However, if it is not the case and if a model bias exists,
 the Dantzig and Lasso procedures can be after all
applied under similar assumptions on $A$. Oracle
inequalities are obtained for which approximation theory
plays an important role in \cite{paralleling,bunea2,bunea,geer}.

Let us also mention that in various settings, under various assumptions on the matrix $A$ (or more precisely on the associated
Gram matrix $G=A^T A$), properties of these
estimators have been established for subset selection (see \cite{bun,karim,meinbul,mein,binyu,zhanghuang}) and for prediction (see \cite{paralleling,knight,karim,mein,zou}).

\subsection{Our goals and results}
We consider in this paper the density estimation framework already studied for
the Lasso estimate by Bunea et al \cite{bunea3,bunea7} and van de Geer \cite{geer}.
Namely, our goal is to estimate $\fo$, an unknown density function, by using
the observations of an $n$-sample of  variables
$X_1,\dots,X_n$ of density $f_0$. As in the non parametric
regression setting, we introduce
a dictionary of functions $\Upsilon=(\p_m)_{m=1,\ldots,M}$, and search again estimates
of $\fo$ as linear combinations $f_\lambda$  of the dictionary
functions. We rely on the Gram matrix associated with $\Upsilon$ and
on the empirical scalar products of $\fo$ with $\p_m$
\[
 \hat{\beta}_m = \frac{1}{n}
\sum_{i=1}^n \p_m(X_i).
\]
The Dantzig  estimate $\hat{f}^D$ is then obtained by minimizing $\normu{\lambda}$
over the set of parameters $\lambda$ satisfying the adaptive Dantzig
constraint:
\[\forall\, m\in\{1,\ldots.M\},\quad\
|(G\lambda)_m-\hb_m|\leq\eta_{\ga,m}
\]
where for $m\in\{1,\ldots,M\}$, $(G\lambda)_m$ is the scalar
product of $f_\lambda$ with $\p_m$,
\[
\eta_{\ga,m}=\sqrt{\frac{2\tilde\si^2_m\ga\log
M}{n}}+\frac{2\norm{\p_m}_\infty\ga\log M}{3n},
\]
$\tilde\si^2_m$ is a sharp estimate of the variance of $\hat\be_m$ and $\gamma$ is a
constant to be chosen. Section \ref{def} gives precise definitions and heuristics for using this constraint. We just mention here that $\eta_{\ga,m}$ comes from sharp concentration inequalities to give tight constraints. Our idea is that if $f_0$ can be decomposed on $\Upsilon$ as \[\fo=\sum_{m=1}^M\la_{0,m}\p_m,\] then we force the set of feasible parameters $\la$ to contain $\la_0$ with large probability and to be as small as possible. Significant improvements in practice are expected.

Our goals in this paper are mainly twofold. First, we aim at
establishing sharp
oracle inequalities under very mild assumptions on the dictionary. Our starting point is that most of the papers in the literature assume that the functions of the
dictionary are bounded by a constant independent of $M$ and $n$, which
constitutes a strong limitation, in particular for dictionaries based
on histograms or wavelets (see for instance \cite{bunea4},
\cite{bunea3}, \cite{bunea2},  \cite{bunea},  \cite{bun} or
\cite{geer}). Such assumptions on the functions of $\Upsilon$ will not be considered in our
paper. Likewise, our methodology does not rely on the knowledge of
$\norm{f_0}_\infty$ that can even be infinite (as noticed by Birg\'e
\cite{bir} for the study of the integrated $\L_2$-risk, most of the
papers in the literature typically assume that the sup-norm of the
unknown density is finite with a known or estimated bound for this
quantity). Finally, let us mention that, in contrast with what Bunea et al \cite{bunea7} did, we obtain oracle inequalities with leading constant 1, and furthermore these are established under much weaker assumptions on the dictionary than in \cite{bunea7}.

The second goal of this paper deals with the problem of calibrating the so-called \emph{Dantzig constant} $\ga$: how should this constant be chosen to obtain good results in both theory and practice? Most of the time, for Lasso-type estimators,
the regularization parameter is of the form $a\sqrt{\frac{\log
M}{n}}$ with $a$ a positive constant (see \cite{paralleling}, \cite{bunea3}, \cite{bunea4}, \cite{bunea}, \cite{candes_plan}, \cite{karim} or \cite{mein} for instance). These results are obtained
with large probability that depends on the tuning coefficient $a$.
In practice, it is not  simple to calibrate the constant $a$. Unfortunately, most of the time, the theoretical choice of the
regularization parameter is not suitable for practical issues. This fact is true for Lasso-type estimates but also for many algorithms for which the
regularization parameter provided by the theory  is often too conservative for practical purposes (see \cite{jll} who clearly explains and illustrates this point for their thresholding procedure).
So, one of the main goals of this paper is to fill the gap between the optimal parameter choice provided by theoretical results on the one hand and by a simulation study on the other hand. Only a few papers are devoted to this problem. In the model selection setting, the issue of calibration has been addressed by Birg\'e and Massart \cite{birgemassart} who considered $\ell_0$-penalized estimators in a Gaussian homoscedastic regression framework and showed that there exists a minimal penalty in the sense that taking smaller penalties leads to inconsistent estimation procedures. Arlot and Massart \cite{arlotmassart} generalized these results for non-Gaussian or heteroscedastic data and Reynaud-Bouret and Rivoirard \cite{reyriv1} addressed this question for thresholding rules in the Poisson intensity framework.

Now, let us describe our results. By using the previous data-driven Dantzig constraint, oracle inequalities are derived under local conditions on the dictionary that are valid under classical assumptions on the structure of the dictionary. We extensively discuss these assumptions and we show their own interest in the context of the paper. Each term of these oracle inequalities is easily interpretable. Classical results are recovered when we further assume:
\[\norm{\p_m}_\infty^2\leq c_1\left(\frac{n}{\log M}\right)\norm{f_0}_\infty~,\]
where $c_1$ is a constant.
This assumption is very mild and, unlike in classical works, allows to consider dictionaries based on wavelets. Then, relying on our Dantzig estimate, we build an adaptive Lasso procedure whose oracle performances are similar. This illustrates the closeness between Lasso and Dantzig-type estimates.

Our results are proved for $\ga>1$. For the theoretical calibration issue, we study the performance of our procedure when $\gamma<1$. We show that in a simple framework, estimation of the straightforward signal $f_0=\indic_{[0,1]}$ cannot be performed  at a convenient
rate  of convergence when $\gamma<1$. This result proves that the assumption $\gamma>1$  is thus not too conservative.

Finally, a simulation study illustrates how dictionary-based methods
outperform classical ones. More precisely, we show that our
Dantzig and Lasso procedures with $\gamma>1$, but close to 1, outperform classical ones, such as
simple histogram procedures, wavelet thresholding or Dantzig procedures based on the knowledge of $\norm{f_0}_\infty$ and less tight Dantzig constraints.

\subsection{Outlines}
Section \ref{def} introduces the density estimator of $f_0$ whose
theoretical performances are studied in Section \ref{results}. Section
\ref{sec:conn-betw-dantz} studies the Lasso estimate proposed in this
paper. The calibration issue is studied in Section \ref{calib} and
numerical experiments are performed in Section
\ref{numerical}. Finally, Section \ref{proofs} is devoted to the
proofs of our results.

\section{The Dantzig estimator of the density $f_0$}\label{def}
As said in Introduction, our goal is to build an estimate of $\fo$  as a linear combination of
functions of $\Upsilon=(\p_m)_{m=1,\ldots,M}$, where we assume without any loss of
generality that, for any $m$, $\|\p_m\|_2=~1$:
\[
f_\lambda = \sum_{m=1}^M \lambda_m \varphi_m.
\] For this purpose, we naturally rely on natural estimates of the $\L_2$-scalar products between $f_0$ and the $\varphi_m$'s. So, for $m\in\{1,\ldots,M\}$, we set
\begin{equation}\label{betazero}
\be_{0,m}=\int  \p_m(x)\fo(x)dx,
\end{equation}
and we consider its empirical counterpart
\begin{equation}\label{empirestim}
\hb_m=\frac{1}{n}\sum_{i=1}^n\p_m(X_i)
\end{equation}
that is an unbiased estimate of $\be_{0,m}$.
The variance of this estimate is
$\var(\hb_m)=\frac{\si_{0,m}^2}{n}$ where
\begin{equation}\label{sigmazero}
\si_{0,m}^2=\int\p_m^2(x)\fo(x)dx-\be_{0,m}^2.
\end{equation}
Note also that for any $\lambda$ and any $m$, the $\L_2$-scalar product between $f_\lambda$ and $\varphi_m$ can be easily computed:
\[
\int  \p_m(x)f_\lambda(x)dx =  \sum_{m'=1}^M \lambda_{m'} \int
\p_{m'}(x) \p_m(x) dx = (G\lambda)_m
\]
where $G$ is the Gram matrix associated to the
dictionary $\Upsilon$ defined for any $1\leq m,m'\leq M$ by
\[
G_{m,m'}=\int \p_{m}(x)\p_{m'}(x)dx.
\]
Any  reasonable choice of  $\lambda$ should ensure that the
coefficients
$(G\lambda)_m$ are close to $\hat\beta_m$ for all $m$.
Therefore,
using Cand\`es and Tao's approach, we define the Dantzig constraint:
\begin{equation}\label{constraint}\forall\, m\in\{1,\ldots.M\},\quad\
|(G\lambda)_m-\hb_m|\leq\eta_{\ga,m}\end{equation}
and the Dantzig
estimate $\hat
f^D$  by \(\hat f^D=f_{\hat\lambda^{D,\gamma}}\) with
\[
\hat\lambda^{D,\gamma}=\argmin_{\lambda\in\R^M} \normu{\lambda}
\quad\mbox{such that $\lambda$ satisfies
the Dantzig
constraint (\ref{constraint})}
,
\]
where for
$\ga>0$ and $m\in\{1,\ldots,M\}$,
\begin{equation}\label{etam}
\eta_{\ga,m}=\sqrt{\frac{2\tilde\si^2_m\ga\log
M}{n}}+\frac{2\norm{\p_m}_\infty\ga\log M}{3n},
\end{equation}
with
\begin{equation}\label{sigmatilde}
\tilde\si^2_m=\hat\si^2_m+2\norm{\p_m}_\infty\sqrt{\frac{2\hat\si^2_m\ga\log
M}{n}}+\frac{8\norm{\p_m}_\infty^2\ga\log M}{n}
\end{equation}
and
\begin{equation}\label{sigmachapeau}
\hat\si^2_m=\frac{1}{n(n-1)}\sum_{i=2}^n\sum_{j=1}^{i-1}(\p_m(X_i)-\p_m(X_j))^2.
\end{equation}
Note that $\eta_{\gamma,m}$ depends on the data, so the constraint (\ref{constraint}) will be referred as the \textit{adaptive Dantzig constraint} in the sequel. We now justify the introduction of the density estimate $\hat f^D$.

The definition of $\eta_{\la,\ga}$ is based on the following heuristics. Given $m$, when  there exists  a constant
$c_0>0$ such that $f_0(x)\geq
c_0$ for $x$ in the support of $\p_m$ satisfying
 $\|\p_m\|_{\infty}^2=o_n(n(\log   M)^{-1})$, then,  with large probability, the   deterministic  term   of
(\ref{etam}) is negligible with respect to the random one. In this case, the random term is the main one and we asymptotically derive
\begin{equation}\label{approxseuil}
\eta_{\gamma,m}\approx\sqrt{2\gamma\,\ln M \,\frac{{\tilde{\si}}^2_{m}}{n}}.
\end{equation}
Having in mind that ${\tilde{\si}}^2_{m}/n$ is a convenient estimate for $\var(\hb_m)$ (see the proof of Theorem \ref{concentrationTh}), the shape of the right hand term of the formula (\ref{approxseuil}) looks like the bound proposed by Cand\`es and Tao \cite{candes} to define the Dantzig constraint in the linear model.
Actually, the deterministic term of
(\ref{etam}) allows to get sharp concentration inequalities. As often done in the literature, instead of estimating $\var(\hb_m)$, we could use the inequality
$$\var(\hb_m)=\frac{{\si}^2_{0,m}}{n}\leq \frac{\norm{f_0}_\infty}{n}$$
and we could replace ${\tilde{\si}}^2_{m}$ with
$\norm{f_0}_\infty$ in the definition of the $\eta_{\gamma,m}$.
But this requires a strong assumption: $f_0$ is bounded and
$\norm{f_0}_\infty$ is known. In our paper, $\var(\hb_m)$ is
estimated, which allows not to impose these conditions. More
precisely, we slightly overestimate $\si^2_{0,m}$ to control large
deviation terms and this is the reason why we introduce
${\tilde{\si}}^2_{m}$ instead of using $\hat\si^2_m$, an unbiased
estimate of $\si_{0,m}^2$. Finally, $\gamma$ is a constant that
has to to be suitably calibrated and plays a capital role in
practice.

The following result justifies previous heuristics by showing that, if $\gamma>1$, with high probability, the
quantity $|\hat\beta_m-\beta_{0,m}|$ is smaller than $\eta_{\gamma,m}$
for all $m$. The parameter $\eta_{\gamma,m}$ with $\gamma$ close to $1$ can be viewed as the ``smallest'' quantity that ensures this property.
\begin{Th}\label{concentrationTh}
Let us assume that $M$ satisfies
\begin{equation}\label{Mn}
n\leq M\leq \exp(n^\delta)
\end{equation}
for $\delta<1$.   Let $\ga>1$. Then, for  any $\e>0$, there  exists a constant
$C_1(\e,\delta,\gamma)$ depending on $\e,$ $\delta$ and $\gamma$ such
that
\[
\P\left(\forall m\in\{1,\ldots,M\},\quad
|\be_{0,m}-\hb_m|\geq\eta_{\ga,m}\right)\leq
C_1(\e,\delta,\gamma)M^{1-\frac{\ga}{1+\e}}.
\]
In addition, there exists a constant $C_2(\delta,\gamma)$ depending on $\delta$ and $\gamma$ such
that
\[
\P\left(\forall m\in\{1,\ldots,M\},\quad \etam_{\ga,m} \leq
\eta_{\ga,m}\leq \etaM_{\ga,m} \right)\leq
C_2(\delta,\gamma)M^{1-\ga}
\]
where, for $m\in\{1,\ldots,M\}$,
\[
\etam_{\gamma,m}=\si_{0,m} \sqrt{\frac{8\ga\log
M}{7n}}+\frac{2\norm{\p_m}_\infty\ga\log M}{3n}
\]
and
\[
\etaM_{\gamma,m}=\si_{0,m}\sqrt{\frac{16 \ga\log
M}{n}}+\frac{10\norm{\p_m}_\infty\ga\log M}{n}.
\]
\end{Th}
This result is proved in Section \ref{proofconcentrationTh}. The first
part is a sharp concentration inequality proved by using Bernstein
type controls. The second part of the theorem proves that,
up to constants depending on $\gamma$,
$\eta_{\gamma,m}$
 is of order $\sigma_{0,m}
\sqrt{\frac{\log M}{n}}+
\norm{\p_m}_\infty\frac{\log M}{n}$ with high probability.
Note that the assumption $\ga>1$ is essential to obtain probabilities going to 0.

Finally, let
$\la_0=(\la_{0,m})_{m=1,\dots,M}\in\R^M$ such that
\[P_{\Upsilon}\fo=\sum_{m=1}^M\la_{0,m}\p_m\] where $P_{\Upsilon}$ is the
projection on the space spanned by $\Upsilon$. We have
\[(G\la_0)_m=\int (P_{\Upsilon}\fo)\p_m=\int \fo\p_m=\be_{0,m}.\]
So, Theorem \ref{concentrationTh} proves that $\la_0$ satisfies the adaptive Dantzig
constraint (\ref{constraint}) with  probability larger  than
$1-C_1(\e,\delta,\gamma)M^{1-\frac{\ga}{1+\e}}$ for any $\e>0$.
Actually, we force the set of parameters $\la$ satisfying the adaptive Dantzig constraint to contain $\la_0$ with large probability and to be as small as possible.
Therefore, $\hat f^D=f_{\hat\lambda^{D,\gamma}}$ is a good candidate among sparse estimates linearly decomposed on $\Upsilon$ for estimating $\fo$.

We mention that Assumption (\ref{Mn}) can be relaxed and we can take $M<n$ provided the definition of $\eta_{\gamma,m}$ is modified.
\section{Results for the Dantzig estimators}\label{results}
In the sequel, we will denote
$\hat\lambda^{D}=\hat\lambda^{D,\gamma}$ to simplify the notations,
but the Dantzig estimator $\hat f^D$ still depends on $\gamma$.
Moreover, we assume that (\ref{Mn}) is true and we denote the vector
$\eta_{\ga}=(\eta_{\ga,m})_{m=1,\ldots,M}$ considered with the Dantzig constant $\gamma>1$.
\subsection{The main result under local assumptions}
Let us state the main result of this paper.  For any $J\subset \{1,\dots,M\}$, we
set $\Comp{J}=\{1,\dots,M\}\smallsetminus J$ and define $\la_J$ the vector which has the same coordinates as $\la$ on $J$ and zero
coordinates on $\Comp{J}$. We introduce a local assumption indexed by
a subset $J_0$.
\begin{itemize}
\item \textbf{Local Assumption} Given $J_0\subset \{1,\dots,M\}$, for some constants $\kappa_{\Jo}>0$ and
$\mu_{\Jo}>0$ depending on $J_0$, we have for any $\lambda$,
\begin{equation}\label{condlocal}
\normD{f_\lambda} \geq \kappa_{\Jo}\normd{\lambda_{\Jo}} - \mu_{\Jo}
 \left( \normu{\lambda_{\Comp{\Jo}}}-\normu{\lambda_{\Jo}}\right)_+.
\tag{$LA(J_0,\kappa_{\Jo},\mu_{\Jo})$}
\end{equation}
\end{itemize}
We obtain the following oracle type inequality
without any assumption on $\fo$.
\begin{Th}\label{oracleadmiloc}
Let $\Jo\subset\{1,\ldots,M\}$ be fixed. We suppose that
\eqref{condlocal} holds. Then, with probability at least
$1-C_1(\e,\delta,\gamma)M^{1-\frac{\gamma}{1+\e}}$, we have for
any $\beta>0$,
\begin{equation}\label{eq2loc}
 \normD{\hat{f}^D-\fo}^2 \leq \inf_{\lambda \in \R^M} \left\{
\normD{f_\lambda-\fo}^2+\beta \frac{
\Lam
^2}{|\Jo|}\left(1+\frac{2\mu_{\Jo}\sqrt{|\Jo|}}{\kappa_{\Jo}}\right)^2+16|\Jo|
\left(\frac{1}{\beta}+\frac{1}{\kappa_{\Jo}^2}\right)
\norminf{\eta_{\gamma}}^2 \right\},
\end{equation}
with
\[
 \Lam= \normu{\lambda_{\Jo^C}}+\frac{\left(\normu{\hat
  \lambda^D}-\normu{\lambda}\right)_+}{2}.
\]
\end{Th}
Let us comment each term of the right hand side of (\ref{eq2loc}).
The first term is an approximation term which measures the closeness
between $\fo$ and $f_\lambda$. This term can vanish if $f_0$ can be
decomposed on the dictionary. The second term is a price to pay when
either $\lambda$ is not supported by the subset $\Jo$ considered or it
does not satisfy the condition $\normu{\hat
  \lambda^D}\leq\normu{\lambda}$ which holds as soon as $\lambda$ satisfy the
adaptive Dantzig constraint.
Finally, the last term, which does not depend on $\lambda$,
 can be viewed as a variance term corresponding to the estimation on the subset $\Jo$. Indeed, remember that $\eta_{\gamma,m}$ relies on an estimate of the variance of $\hat\be_m$. Furthermore, we have with high probability:
\begin{eqnarray*}
\norminf{\eta_{\gamma}}^2&\leq&2\left(\frac{16\si_{0,m}^2 \ga\log
M}{n}+\left(\frac{10\norm{\p_m}_\infty\ga\log M}{n}\right)^2\right).
\end{eqnarray*}
So,  if $f_0$ is bounded then,  $\si_{0,m}^2 \leq \norm{f_0}_\infty$ and if there exists a constant $c_1$ such that for any $m$,
\begin{equation}\label{conddicob}
\norm{\p_m}_\infty^2\leq c_1\left(\frac{n}{\log M}\right)\norm{f_0}_\infty,
 \end{equation}
(which is true for instance for a bounded dictionary), then
\[\norminf{\eta_{\gamma}}^2\leq C\norm{f_0}_\infty\frac{\log
M}{n},\]
(where $C$ is a constant depending on $\gamma$ and $c_1$) and tends to 0 when $n$
goes to $\infty$. We obtain thus the following result.
\begin{Cor}
Let $\Jo\subset\{1,\ldots,M\}$ be fixed. We suppose that \eqref{condlocal}
holds. If (\ref{conddicob}) is satisfied
then, with probability at least
$1-C_1(\e,\delta,\gamma)M^{1-\frac{\gamma}{1+\e}}$, we have for  any
$\beta>0$,
for any  $\lambda$ that satisfies the adaptive Dantzig constraint
\begin{equation}\label{eq2terbis}
\normD{\hat{f}^D-\fo}^2 \leq \normD{f_\lambda-\fo}^2+\beta c_2
(1+\kappa_{\Jo}^{-2}\mu_{\Jo}^2|\Jo|)\frac{\normu{\lambda_{\Jo^C}}
  ^2}{|\Jo|}+ c_3 (\beta^{-1}+\kappa_{\Jo}^{-2})|\Jo|\norm{f_0}_\infty\frac{\log
M}{n},
\end{equation}
where $c_2$ is an absolute constant and $c_3$ depends on $c_1$ and
$\gamma$.
\end{Cor}
The parameter $\beta$ calibrates the weights given
for the bias and variance terms.
Remark that if $f_0=f_{\lambda_0}$ and if \eqref{condlocal}
holds with $\Jo=J_{\lambda_0}$, under (\ref{conddicob}), the proof of
Theorem~\ref{oracleadmiloc}
yields the more classical
inequality
\[
 \normD{\hat{f}^D-\fo}^2 \leq
C'
|\Jo|
\norm{f_0}_\infty\frac{\log
M}{n},
\]
where $C'=c_3\kappa_{J_0}^{-2}$,
with at least the same probability
$1-C_1(\e,\delta,\gamma)M^{1-\frac{\gamma}{1+\e}}$.

Assumption~\eqref{condlocal} is local, in the sense that
the constants $\kappa_{\Jo}$ and
$\mu_{\Jo}$ (or their mere existence) may highly depend on the subset
$\Jo$. For a given $\lambda$, the best choice for $J_0$ in
Inequalities~\eqref{eq2loc} and~\eqref{eq2terbis} depends
thus on the interaction between these constants and the value of
$\lambda$ itself.
Note that the assumptions of Theorem~\ref{oracleadmiloc}
are reasonable as
the next section gives conditions for which
Assumption~\eqref{condlocal} holds simultaneously
with the same constant $\kappa$ and $\mu$ for all subsets $J_0$ of
the same size.

\subsection{Results under global assumptions}
As usual,  when $M>n$, properties of the Dantzig estimate can be derived from assumptions on the structure of the dictionary $\Upsilon$. For $l\in\Ne$, we denote
\begin{align*}
\phi_{\min}(l)&
=\min_{|J|\leq l}\min_{\substack{\lambda\in\R^M\\\lambda_J\neq0}}
\frac{\normD{f_{\lambda_J}}^2}{\normd{\lambda_J}^2}&\text{and}&&
\phi_{\max}(l)&
=\max_{|J|\leq l}\max_{\substack{\lambda\in\R^M\\\lambda_J\neq0}}
\frac{\normD{f_{\lambda_J}}^2}{\normd{\lambda_J}^2}.
\end{align*}
These quantities correspond to the ``restricted'' eigenvalues of the Gram matrix $G$. Assuming that $\phi_{\min}(l)$ and $\phi_{\max}(l)$ are close to 1 means that every set of columns of $G$ with cardinality less than $l$  behaves like an orthonormal system. We also consider the restricted correlations
\[
\theta_{l,l'}=
\max_{\substack{\ |J|\le l\\\ |J'|\le l'\\J\cap
    J'=\emptyset}}
\max_{\substack{\lambda,\lambda'\in \R^M\\
\lambda_J\neq0,\lambda'_{J'}\neq0
}
}
 \frac{\langle f_{\lambda_{J}}, f_{\lambda'_{J'}} \rangle}
{\normd{\lambda_{J}}\normd{\lambda'_{J'}}}.
\]
Small values of $\theta_{l,l'}$ mean that two disjoint sets of columns of $G$ with cardinality less than $l$ and $l'$ span nearly orthogonal spaces. We will use one of the following assumptions considered in \cite{paralleling}.
\begin{itemize}\item\textbf{Assumption 1}
For some integer $1\le s\le M/2$, we have
\begin{equation}\label{Ass1}
\phi_{\min}(2s)> \theta_{s,2s}.\tag{\mbox{A1}(s)}
\end{equation}
Oracle inequalities of the Dantzig selector were established  under this assumption in the parametric linear model by Cand\`es and Tao in \cite{candes}. It was also considered by Bunea, Ritov and Tsybakov \cite{paralleling} for non-parametric regression and for the Lasso estimate. The next assumption, proposed in \cite{paralleling}, constitutes an alternative to Assumption 1.
\item\textbf{Assumption 2}
For some integers $s$ and $l$  such that
\begin{equation}\label{sl}
1\le s\le \frac{M}{2},\quad l\ge s\quad \mbox{and}\quad s+l\le M,
\end{equation}
we have
\begin{equation}\label{Ass2}
l\phi_{\min}(s+l)> s\phi_{\max}(l).\tag{\mbox{A2}(s,l)}
\end{equation}
If Assumption 2 is
true for $s$ and $l$ such that $l\gg s$, then Assumption 2 means that $\phi_{\min}(l)$ cannot decrease at a rate faster than $l^{-1}$ and this condition is related to the ``incoherent designs'' condition stated in \cite{mein}.
\end{itemize}
In the sequel, we set, under Assumption 1,
\[
\kappa_1(s)=\sqrt{\phi_{\min}(2s)}\left(1-\frac{\theta_{s,2s}}{\phi_{\min}(2s)}\right)>0,\quad
\mu_1(s)=\frac{\theta_{s,2s}}{\sqrt{s\phi_{\min}(2s)}}
\]
and under Assumption 2,
\[
\kappa_2(s,l)=\sqrt{\phi_{\min}(s+l)}\left(1-\sqrt{\frac{s\phi_{\max}(l)}{l\phi_{\min}(s+l)}}\right)>0,\quad
\mu_2(s,l)=\sqrt{\frac{\phi_{\max}(l)}{l}}.
\]
Now, to apply Theorem \ref{oracleadmiloc}, we need to check (\ref{condlocal}) for some some subset $J_0$ of $\{1,\ldots,M\}$.
Either Assumption~1 or Assumption~2 implies this assumption. Indeed, we have the following result.
\begin{Prop}\label{lema0}
Let $s$ and $l$ two integers satisfying (\ref{sl}). We suppose that (\ref{Ass1}) or (\ref{Ass2}) is true.  Let $J_0\subset\{1,\ldots,M\}$ of size $|J_0|=s$ and $\lambda\in\R^M$,
then we have
\[\normD{f_\lambda} \geq\kappa
 \normd{\lambda_{\Jo}}-\mu \left( \normu{\lambda_{\Comp{\Jo}}}-\normu{\lambda_{\Jo}}\right)_+
\]
with $\kappa=\kappa_1(s)$ and $\mu=\mu_1(s) $  under (\ref{Ass1}) (respectively $\kappa=\kappa_2(s,l)$ and $\mu=\mu_2(s,l)$ under (\ref{Ass2}). If (\ref{Ass1}) and (\ref{Ass2})  are both satisfied,
$\kappa=\max(\kappa_1(s),\kappa_2(s,l))$ and $\mu=\min(\mu_1(s),\mu_2(s,l))$.
\end{Prop}
Proposition \ref{lema0} proves that Theorem \ref{oracleadmiloc} can be
applied under Assumptions 1 or 2. In addition, the constants
$\kappa_{\Jo}$ and $\mu_{\Jo}$ only depend on $|J_0|$.
 From Theorem \ref{oracleadmiloc}, we deduce the following result.
\begin{Th}\label{oracleadmi}
Let $s$ and $l$ two integers satisfying (\ref{sl}). We suppose
that (\ref{Ass1}) or (\ref{Ass2}) is true.  Then, with probability
at least $1-C_1(\e,\delta,\gamma)M^{1-\frac{\gamma}{1+\e}}$, we
have for any $\beta>0$,
\[\normD{\hat{f}^D-\fo}^2 \leq \inf_{\lambda\in\R^M}
\inf_{\substack{J_0 \subset \{1,\ldots,M\}\\
  |J_0|=s}} \left\{ \normD{f_\lambda-\fo}^2+\beta \frac{
\Lam^2}{s}\left(1+\frac{2\mu\sqrt
    s}{\kappa}\right)^2+16s
\left(\frac{1}{\beta}+\frac{1}{\kappa^2}\right)
\norminf{\eta_{\gamma}}^2 \right\}
\]
where
\[
\Lam = \normu{\lambda_{\Jo^C}}+\frac{\left(\normu{\hat
  \lambda^D}-\normu{\lambda}\right)_+}{2}.
\]
\end{Th}
Remark that the best subset $J_0$ of cardinal $s$ in Theorem~\ref{oracleadmi}
can be easily chosen for a given $\lambda$: it is given by the set of
the $s$ largest coordinates of $\lambda$. This was not necessarily the case in
Theorem~\ref{oracleadmiloc} for which a different subset may give a better local condition and then may provide
a smaller bound. If we further assume the mild assumption (\ref{conddicob}) on the sup norm of the
dictionary introduced in the previous section, we deduce the following result.
\begin{Cor}
Let $s$ and $l$ two integers satisfying (\ref{sl}). We suppose that (\ref{Ass1}) or (\ref{Ass2}) is true. If  (\ref{conddicob}) is satisfied,
with probability at least
$1-C_1(\e,\delta,\gamma)M^{1-\frac{\gamma}{1+\e}}$, we have for  any
$\beta>0$,
 any  $\lambda$ that satisfies the adaptive Dantzig constraint and
for the best subset $J_0$ of cardinal $s$ (that corresponds to the $s$
largest coordinates of $\lambda$ in absolute value),
\begin{equation}\label{eq2ter}
\normD{\hat{f}^D-\fo}^2 \leq \normD{f_\lambda-\fo}^2+
\beta c_2
(1+\kappa^{-2}\mu^2s)\frac{\normu{\lambda_{\Jo^C}}
  ^2}{s}+ c_3 (\beta^{-1}+\kappa^{-2})s\norm{f_0}_\infty\frac{\log
M}{n},
\end{equation}
where $c_2$ is an absolute constant and $c_3$ depends on $c_1$ and
$\gamma$.
\end{Cor}
Note that, when $\lambda$ is $s$-sparse so that $\lambda_{\Jo^C}=0$,
the oracle inequality (\ref{eq2ter}) corresponds to the classical
oracle inequality obtained in parametric frameworks (see
\cite{candes_plan} or  \cite{candes} for instance) or  in
non-parametric settings. See, for instance \cite{bunea4},
\cite{bunea3}, \cite{bunea2},  \cite{bunea},  \cite{bun} or
\cite{geer} but in these works, the functions of the dictionary are
assumed to be bounded by a constant independent of $M$ and
$n$. So, the adaptive Dantzig estimate requires weaker conditions
since under (\ref{conddicob}), $\norm{\p_m}_\infty$ can go to $\infty$
when $n$ grows.
 This point is capital for practical purposes, in particular when wavelet bases are considered.
\section{Connections between the Dantzig and Lasso estimates}\label{sec:conn-betw-dantz}
We show in this section the strong connections between Lasso and Dantzig estimates, which has already been illustrated in \cite{paralleling} for non-parametric regression models. By choosing convenient random weights depending on $\eta_{\gamma}$ for $\ell_1$-minimization, the Lasso estimate satisfies the adaptive Dantzig constraint. More precisely, we consider the Lasso estimator given by the solution of the
following minimization problem
\begin{equation}\label{lasso}
\hat \lambda^{L,\gamma} = \argmin_{\lambda\in\mathbb{R}^M}\left\{
R(\lambda)+2\sum_{m=1}^M \eta_{\gamma,m} |\lambda_m|\right\},
\end{equation}
where
\[
R(\lambda)=\normD{f_\lambda}^2-\frac{2}{n}\sum_{i=1}^nf_\lambda(X_i).
\]
Note that $R(\cdot)$ is the quantity minimized in unbiased
estimation of the risk. For simplifications, we write
$\hat{\lambda}^L=\hat \lambda^{L,\gamma}$. We denote
$\hat{f}^L=f_{\hat{\lambda}^L}$. As said in Introduction, classical Lasso estimates are defined as the minimizer of expressions of the form
\[\left\{
R(\lambda)+2\eta\sum_{m=1}^M |\lambda_m|\right\},\]
where $\eta$ is proportional to $\sqrt{\frac{\log M}{n}}$. So, $\hat{\lambda}^L$ appears as a data-driven version of classical Lasso estimates.

The first order condition for the minimization of the expression given in (\ref{lasso})
corresponds exactly to the adaptive Dantzig constraint and thus
Theorem~\ref{oracleadmi} always applies to $\hat{\lambda}^L$.
Working along the lines of the proof of Theorem~\ref{oracleadmi}
(Replace $f_\lambda$ by $\hat{f}^{D}$ and $\hat{f}^{D}$ by $\hat{f}^{L}$ in
(\ref{rela1}) and (\ref{maj2})), one can prove a slightly stronger
result.
\begin{Th}\label{equivalence}
Let us assume that assumptions of Theorem \ref{oracleadmi} are
true. Let $J_0\subset\{1,\ldots,M\}$ of size $|J_0|=s$. Then, with
probability at least
$1-C_1(\e,\delta,\gamma)M^{1-\frac{\gamma}{1+\e}}$, we have for
any $\beta>0$,
\[\left|\normD{\hat{f}^{D}-\fo}^2 -
\normD{\hat{f}^{L}-\fo}^2\right| \leq \beta \frac{
\normu{\hat{\lambda}^{L}_{\Jo^C}} ^2}{s}\left(1+\frac{2\mu\sqrt
    s}{\kappa}\right)^2+16s \left(\frac{1}{\beta}+\frac{1}{\kappa^2}\right)
\norminf{\eta_{\gamma}}^2.
\]
\end{Th}
To extend this theoretical result, numerical performances of the Dantzig and Lasso estimates will be compared in Section \ref{numerical}.
\section{Calibration and numerical experiments}
\subsection{The calibration issue}\label{calib}
In this section, we consider the problem of calibrating previous
estimates. In particular, we prove that the sufficient condition
$\ga>1$ is ``almost'' a necessary condition since we derive a
special and very simple framework in which Lasso and Dantzig
estimates cannot achieve the optimal rate if $\ga<1$ (``almost'' means that the case $\gamma=1$ remains an open question). Let us
describe this simple framework. The dictionary $\Upsilon$ considered in this section is the orthonormal Haar system:
\[\Upsilon=\left\{\phi_{jk}:\quad -1\leq  j\leq j_0,\ 0 \leq  k< 2^j\right\},\]
with $\phi_{-10}=\indic_{[0,1]}$, $2^{j_0+1}=n$, and for $ 0\leq j\leq
j_0,$ $0 \leq  k\leq 2^j-1$,
\[\phi_{jk}=2^{j/2}\left(1_{[k/2^j,(k+0.5)/2^j]}-1_{[(k+0.5)/2^j,(k+1)/2^j]}\right).\]
In this case, $M=n$. In this setting, since functions of $\Upsilon$ are orthonormal, the Gram matrix $G$ is the identity. Thus, the Lasso and Dantzig estimates both correspond to the soft thresholding rule:
$$\hat{f}^{D}=\hat{f}^{L}=\sum_{m=1}^M\mbox{sign}(\hat\beta_m)\left(|\hat\beta_m|-\eta_{\gamma,m}\right)1_{\{|\hat\beta_m|>\eta_{\gamma,m}\}}\p_m.$$
Now, our goal is to estimate
$\fo=\phi_{-10}=\indic_{[0,1]}$ by using $\hat{f}^{D}$ depending on $\gamma$ and to show the influence of this constant. Unlike previous results stated in probability, we consider the expectation of the $\L_2$-risk:
\begin{Th}\label{lower}
On the one hand, if $\ga>1$, there exists a constant $C$ such that
\begin{equation}\label{majorisk}
\E\norm{\hat f^D-\fo}_2^2\leq \frac{C\log n}{n},
\end{equation}
On the other hand, if $\ga<1$, there exists a constant $c$  and $\delta<1$  such that
\begin{equation}\label{minorisk}
\E\norm{\hat f^D-\fo}^2_2\geq \frac{c}{n^{\delta}}.
\end{equation}
\end{Th}
This result shows that choosing $\gamma<1$ is a bad choice in our
setting. Indeed, in this case, the Lasso and Dantzig estimates cannot
estimate a very simple signal ($f_0=\indic_{[0,1]}$)  at a convenient
rate  of convergence.

\begin{figure}
  \centering
  \parinc{8cm}{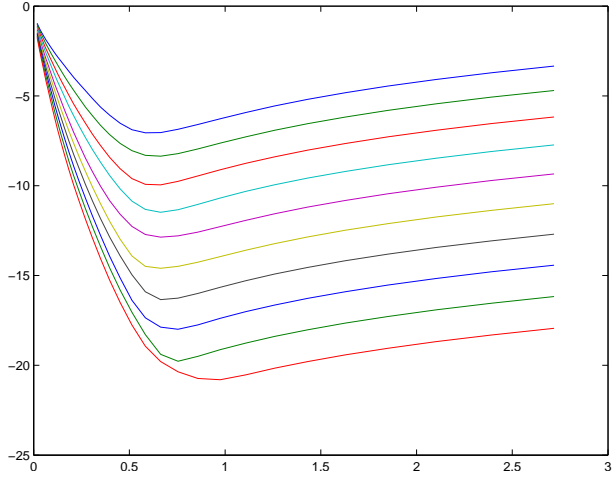}
  \caption{Graphs of $\gamma\mapsto\log_2(\overline{R_n}(\gamma))$ for $n=2^J$ with, from top to bottom,  $J=4,5,6,\ldots,13$}
\label{fig:penalmini}
\end{figure}

A small simulation study is carried out to strengthen this theoretical asymptotic result.
Performing our estimation procedure 100 times,
we compute the average risk $\overline{R_n}(\gamma)$ for
several values of the Dantzig constant $\gamma$ and several values of $n$. This computation is
summarized in Figure~\ref{fig:penalmini} which displays the logarithm of $\overline{R_n}(\gamma)$ for $n=2^J$ with, from top to bottom,  $J=4,5,6,\ldots,13$ on a grid of $\gamma$'s around $1$.
To discuss our results, we denote by $\gamma_{\min}(n)$ the best $\gamma$:
\(
\gamma_{\min}(n)=\argmin_{\gamma>0}\overline{R_n}(\gamma).
\)
 We note that
$1/2 \leq \gamma_{\min}(n)\leq 1$ for all values of $n$, with
$\gamma_{\min}(n)$ getting closer to $1$ as $n$ increases. Taking
$\gamma$ too small strongly deteriorates the performance while a value
close to $1$ ensures a risk withing a factor $2$ of the optimal risk.
The assumption $\gamma>1$ giving a theoretical control on the
quadratic error is thus not too conservative.
Following these
results, we set $\gamma=1.01$ in our numerical experiments in the next subsection.

\subsection{Numerical experiments}\label{numerical}

In this section, we present our numerical experiments with the Dantzig
density estimator and their results. We test our estimator with a
collection of 6 dictionaries, 4
densities described below and for 2 sample sizes.
We compare our procedure with the adaptive Lasso
introduced in Section~\ref{sec:conn-betw-dantz} and with a non
adaptive Dantzig
estimator. We also consider a two-step
estimation procedure, proposed by Cand\`es and Tao~\cite{candes}, which
improves the numerical results.

The numerical scheme for a given dictionary
$\Upsilon=(\p_m)_{m=1,\ldots,M}$
and a sample
$(X_i)_{i=1,\ldots,n}$ is the following.
\begin{enumerate}
\item Compute $\hat{\beta}_m$ for
  all $m$,
\item Compute $\hat{\sigma}^2_m$,
\item Compute $\eta_{\gamma,m}$ as defined in \eqref{etam} by
\begin{equation*}
\eta_{\ga,m}=\sqrt{\frac{2\tilde\si^2_m\ga\log
M}{n}}+\frac{2\norm{\p_m}_\infty\ga\log M}{3n},
\end{equation*}
with
\begin{equation*}
\tilde\si^2_m=\hat\si^2_m+2\norm{\p_m}_\infty\sqrt{\frac{2\hat\si^2_m\ga\log
M}{n}}+\frac{8\norm{\p_m}_\infty^2\ga\log M}{n}
\end{equation*}
and $\ga=1.01$.
\item Compute the coefficients $\hat\lambda^{D,\gamma}$ of the Dantzig
  estimate, $\hat\lambda^{D,\gamma}=\argmin_{\lambda\in\R^M}
  \normu{\lambda}$ such that $\lambda$ satisfies the Dantzig
constraint (\ref{constraint})
\[
\forall\, m\in\{1,\ldots.M\},\quad\
|(G\lambda)_m-\hb_m|\leq\eta_{\ga,m}
\]
 with the homotopy-path-following method proposed
by Asif and Romberg \cite{asif},
\item Compute the Dantzig estimate $\hat f^{D,\gamma}= \sum_{m=1}^M
  \hat\lambda^{D,\gamma}_m \phi_m$.
\end{enumerate}
Note that we have implicitly assumed that the Gram matrix $G$ used in the
definition of the Dantzig constraint has been precomputed.

For the Lasso estimator, the Dantzig minimization of step 4 is
replaced by the Lasso minimization~\eqref{lasso}
\[
\hat \lambda^{L,\gamma} = \argmin_{\lambda\in\mathbb{R}^M}\left\{
R(\lambda)+2\sum_{m=1}^M \eta_{\gamma,m} |\lambda_m|\right\},
\]
which is solved using the LARS algorithm. The non adaptive Dantzig estimate is obtained by replacing $\tilde{\sigma}_m^2$ in step $3$ by
$\|\fo\|_\infty.$ The two-step procedure of Cand\`es and Tao adds a least-square step
between step 4 and step 5. More precisely, let $\hat{J}^{D,\gamma}$ be the
support of the estimate $\hat{\lambda}^{D,\gamma}$. This defines a
subset of the dictionary on which the density is regressed
\[
\left(\hat{\lambda}^{D+LS,\gamma}\right)_{\hat{J}^{D,\gamma}}
= G_{\hat{J}^{D,\gamma}}^{-1} (\hat{\beta}_m)_{\hat{J}^{D,\gamma}}
\]
where $G_{\hat{J}^{D,\gamma}}$ is the submatrix of $G$ corresponding to
  the subset chosen. The values of  $\hat{\lambda}^{D+LS,\gamma}$
  outside $\hat{J}^{D,\gamma}$ are set to 0 and
  $\hat{f}^{D+LS,\gamma}$ is set accordingly.

We describe now the dictionaries we consider.
We focus
numerically on
densities defined on the interval $[0,1]$ so we use dictionaries
adapted to this setting.
The first four
are orthonormal systems, which are used as a benchmark, while the last
two are ``real'' dictionaries. More precisely, our dictionaries are
\begin{itemize}
\item the Fourier basis with $M=n+1$ elements (denoted ``Fou''),
\item the histogram collection  with the classical number $\sqrt{n}/2\leq
  M=2^{j_0} < \sqrt{n}$ of bins (denoted ``Hist''),
\item the Haar wavelet basis with maximal resolution $n/2 <
  M=2^{j_1} < n$ and thus $M=2^{j_1}$ elements (denoted ``Haar''),
\item the more regular Daubechies 6 wavelet basis with  maximal
  resolution $n/2 \leq  M=2^{j_1} < n$ and thus $M=2^{j_1}$ elements
  (denoted ``Wav''),
\item the dictionary made of the union of the Fourier basis and the
  histogram collection and thus comprising $M=n+1+2^{j_0}$ elements.
(denoted ``Mix''),
\item the dictionary which is the union of the Fourier basis, the
  histogram collection and the Haar wavelets of resolution greater than
  $2^{j_0}$ comprising $M=n+1+2^{j_1}$ elements (denoted ``Mix2'').
\end{itemize}
The orthonormal families we have chosen are often used by
practitioners. Our dictionaries combine very different orthonormal
families, sine and cosine with bins or Haar wavelets, which
ensures a sufficiently incoherent design.

\begin{figure}
\def\widthplot1{4.75cm}
  \centering
  \begin{tabular}{cccc}
& $n=500$ & $n=1000$ & $n=2000$\\
$f_1$ & \parinc{\widthplot1}{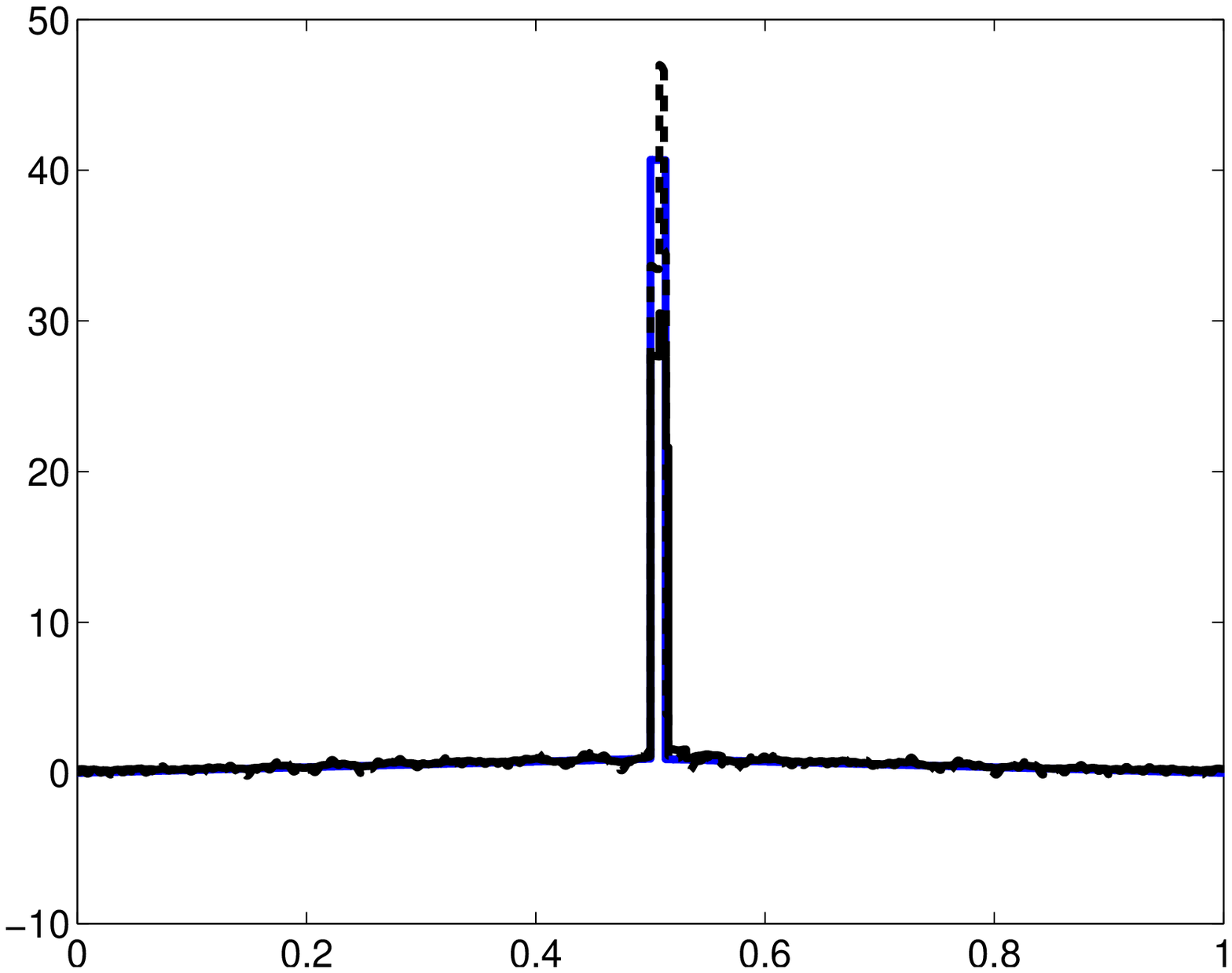} &
\parinc{\widthplot1}{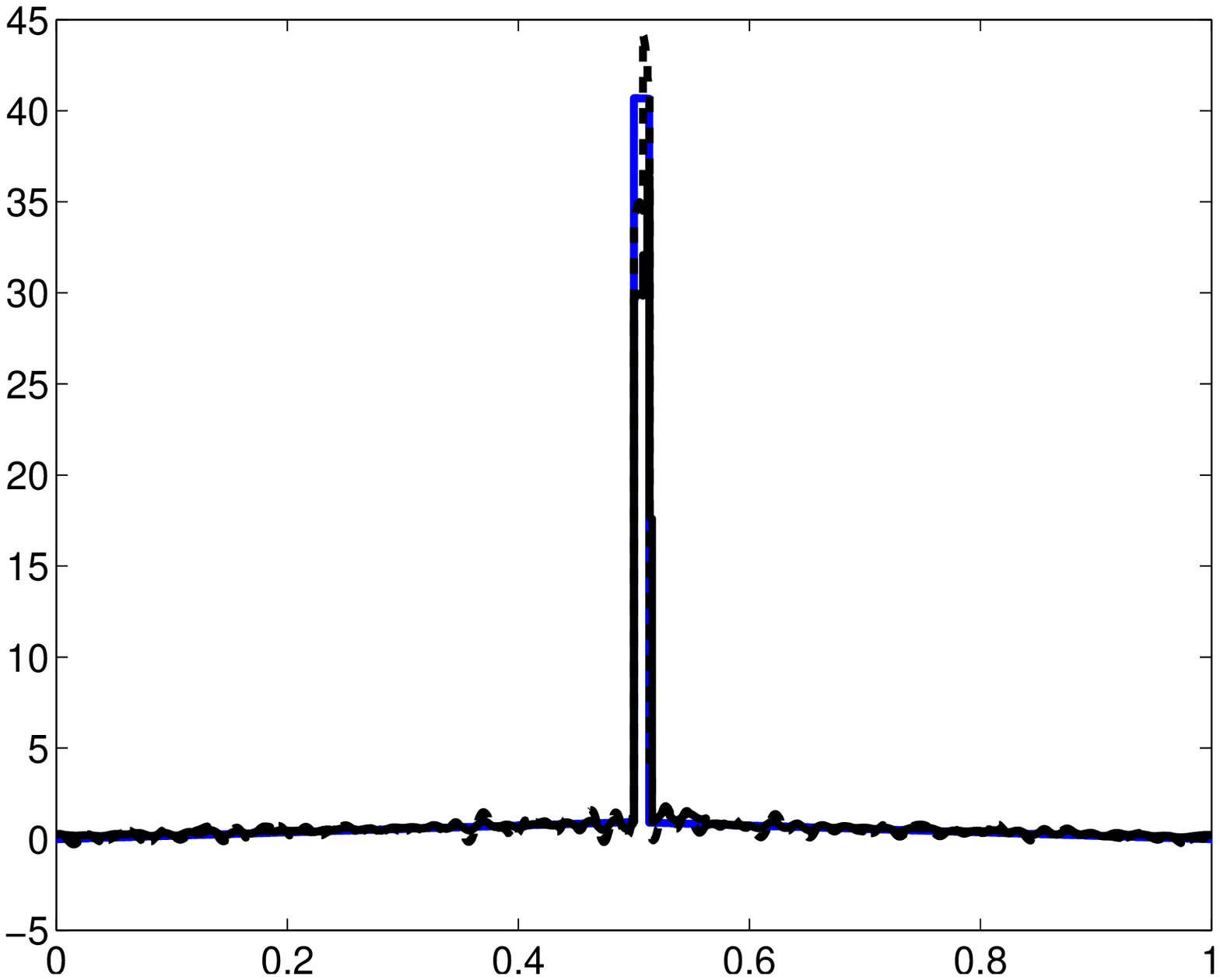} &  \parinc{\widthplot1}{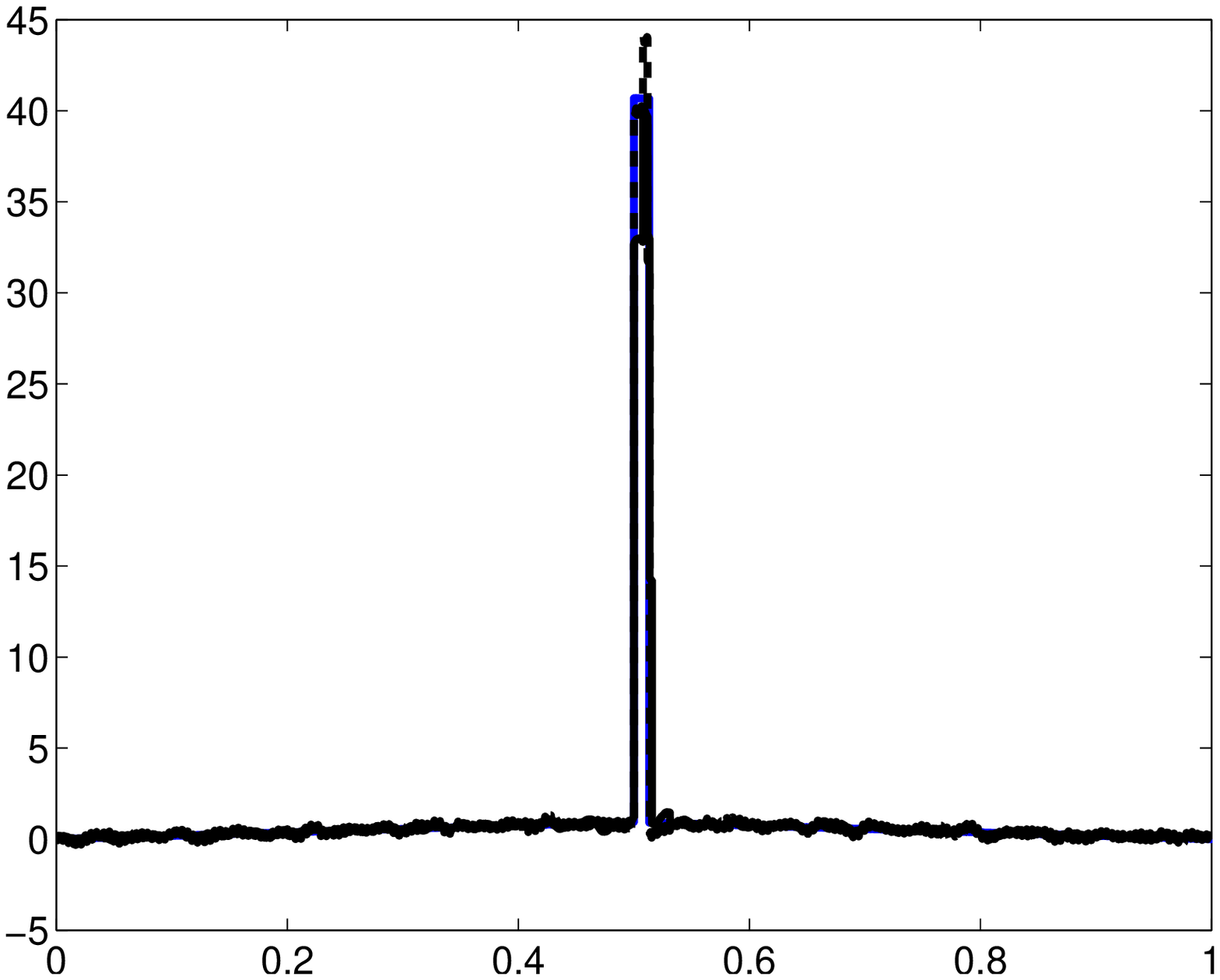}\\
$f_2$ & \parinc{\widthplot1}{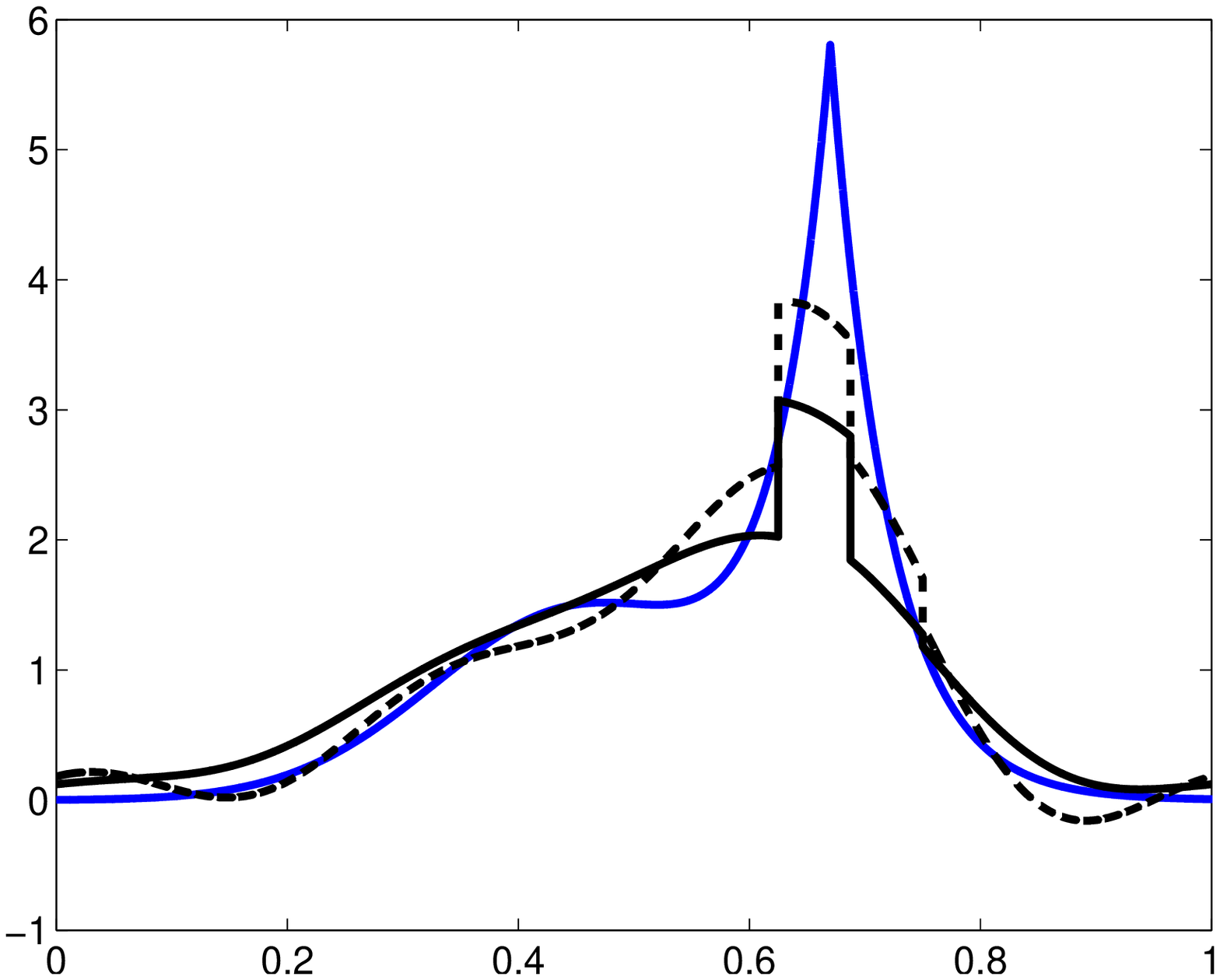} &
\parinc{\widthplot1}{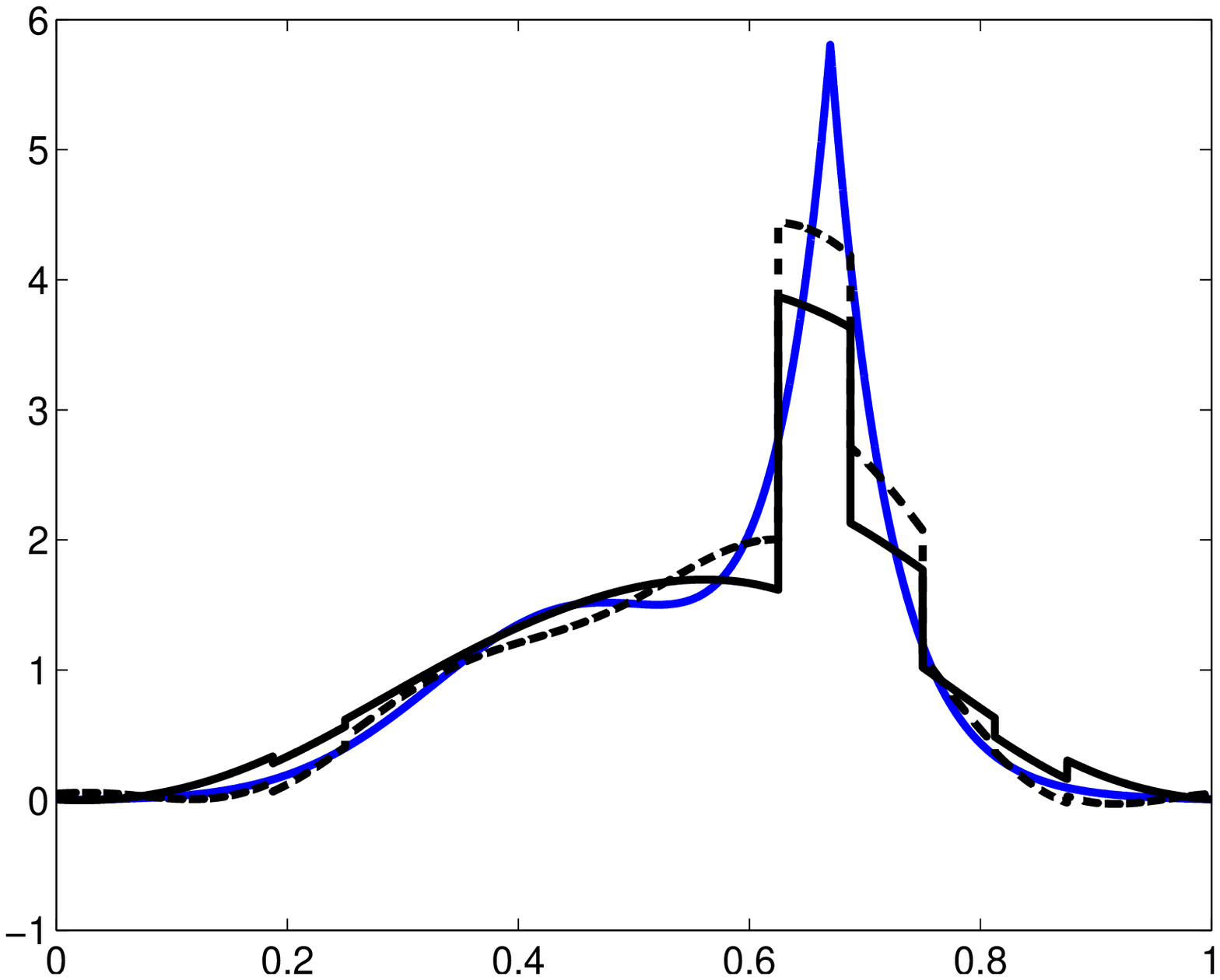} &  \parinc{\widthplot1}{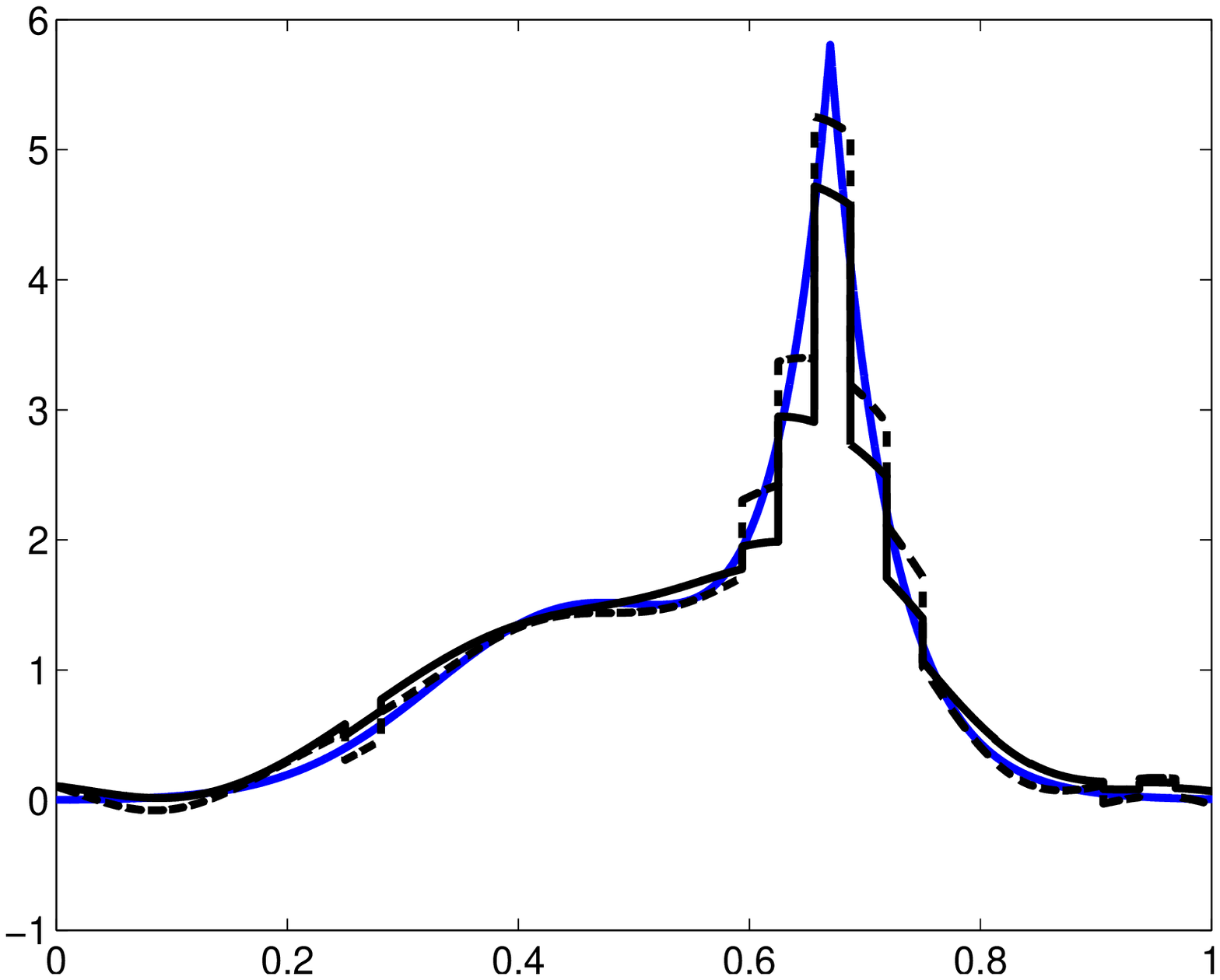}\\
$f_3$ & \parinc{\widthplot1}{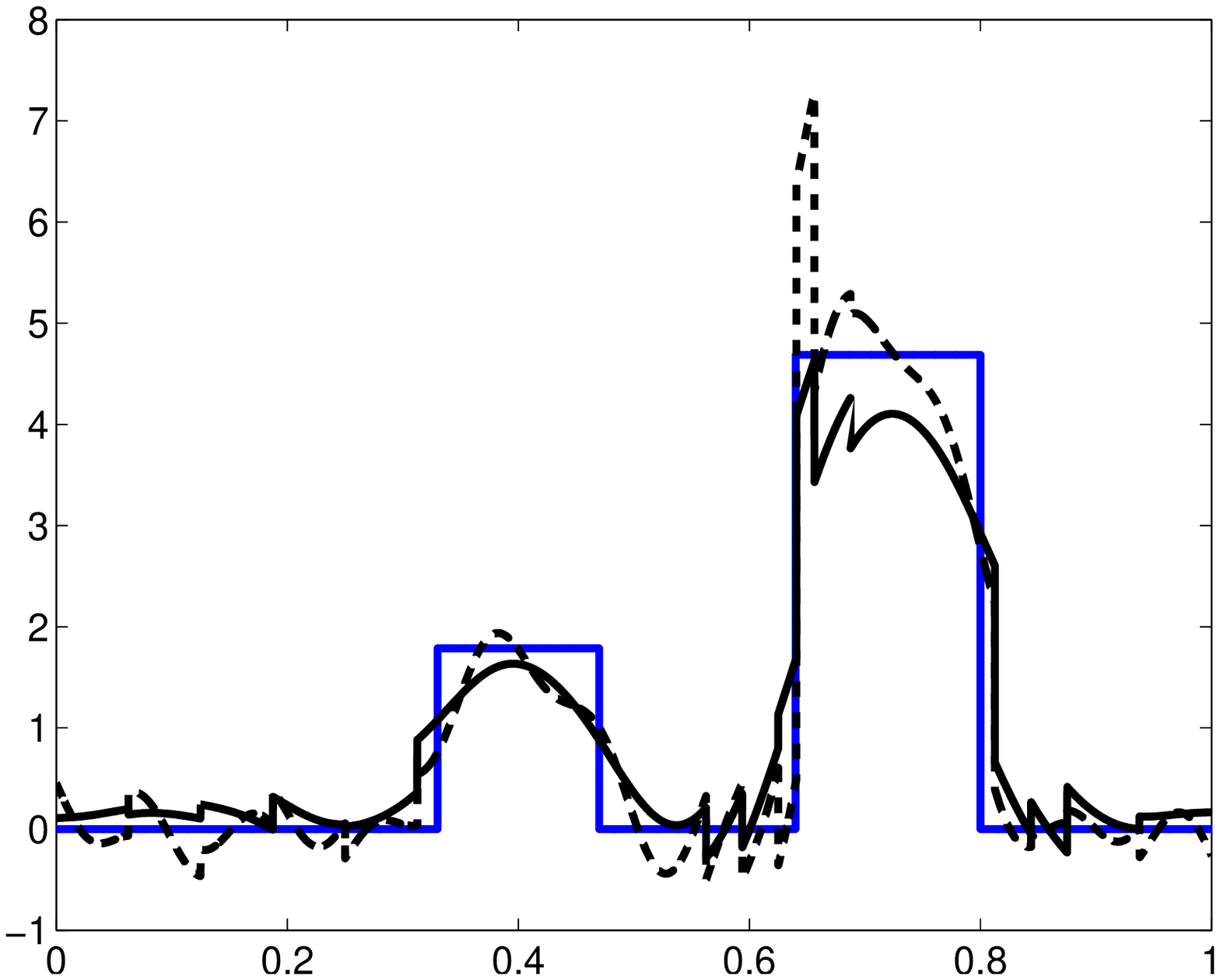} &
\parinc{\widthplot1}{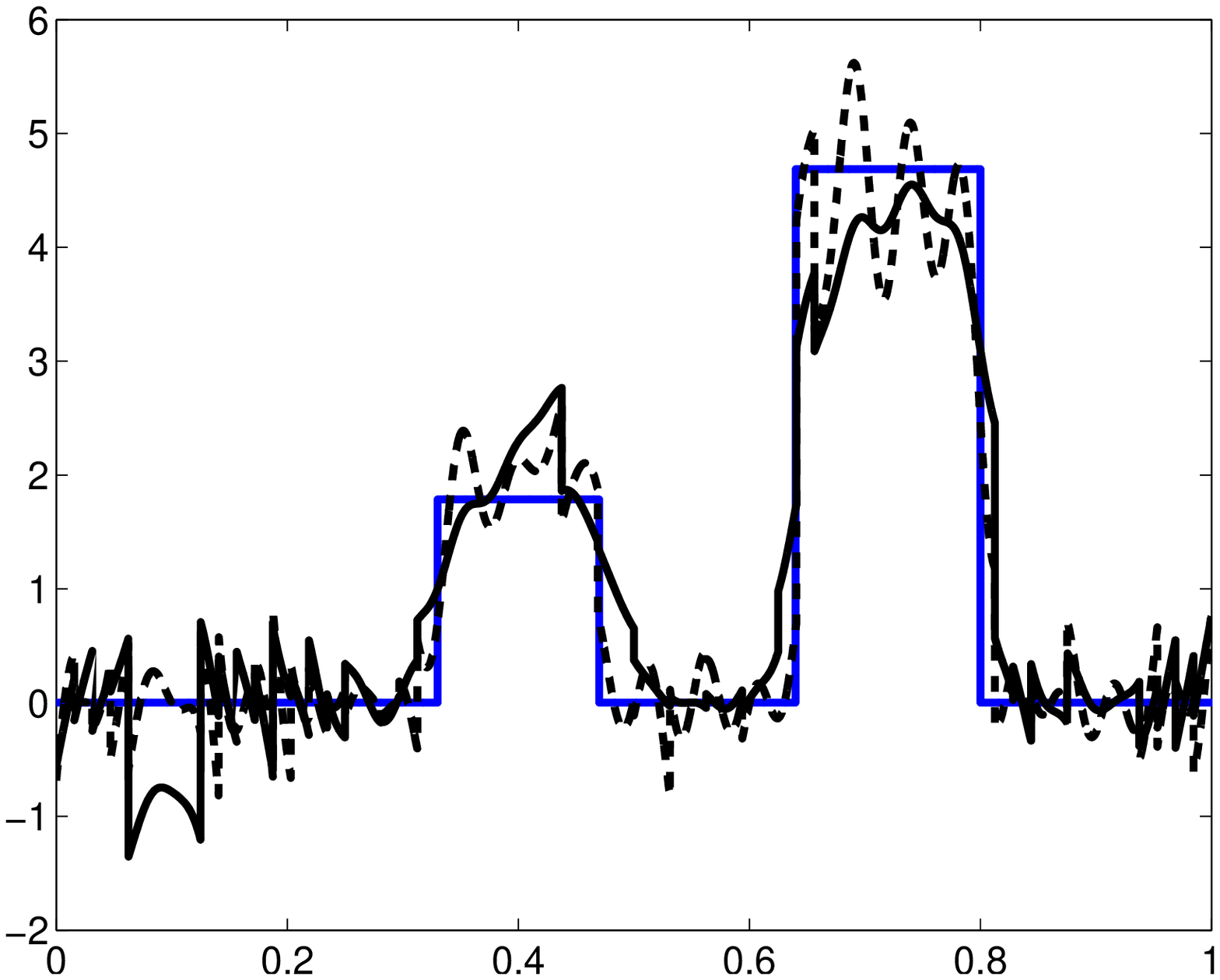} &  \parinc{\widthplot1}{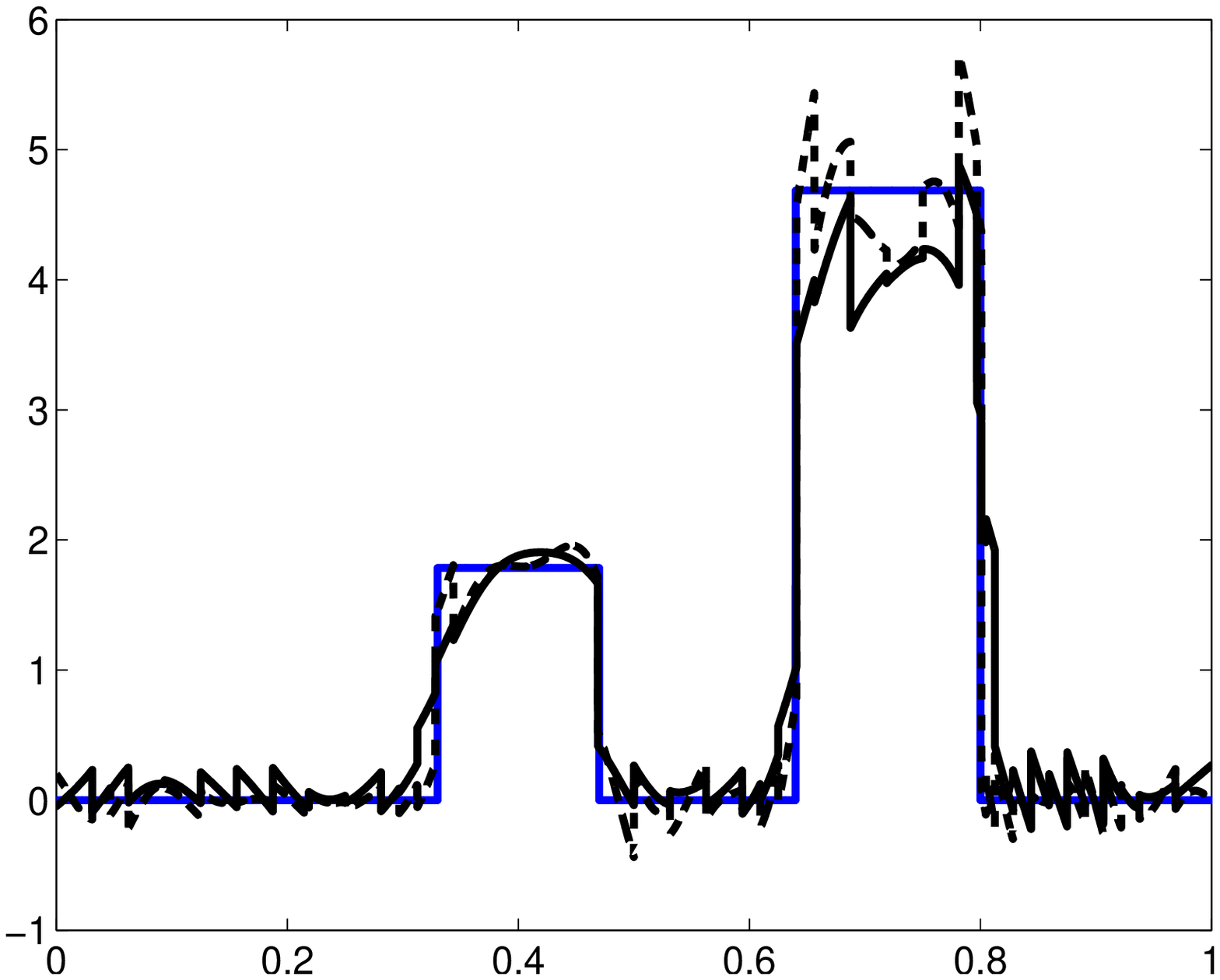}\\
$f_4$ & \parinc{\widthplot1}{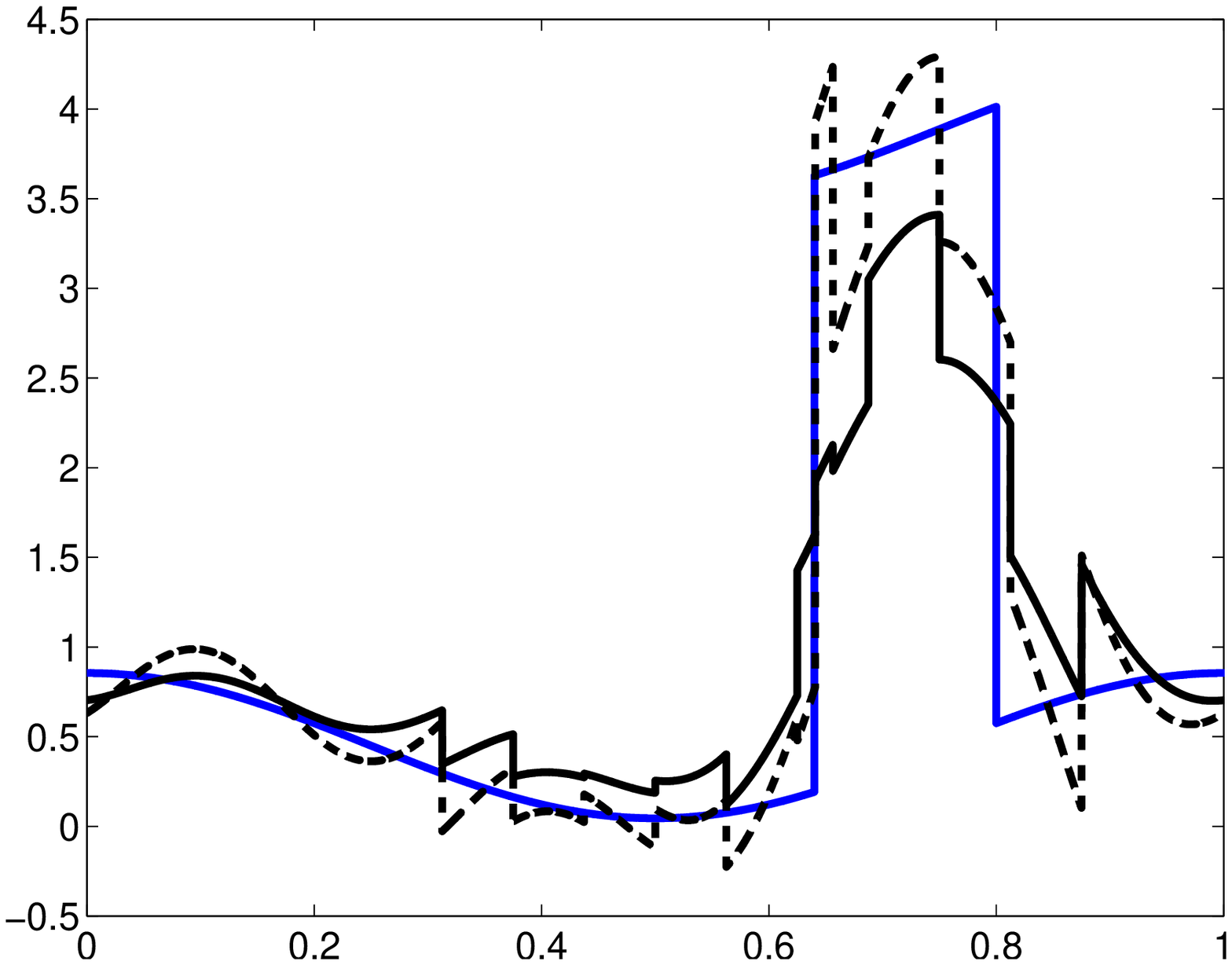} &
\parinc{\widthplot1}{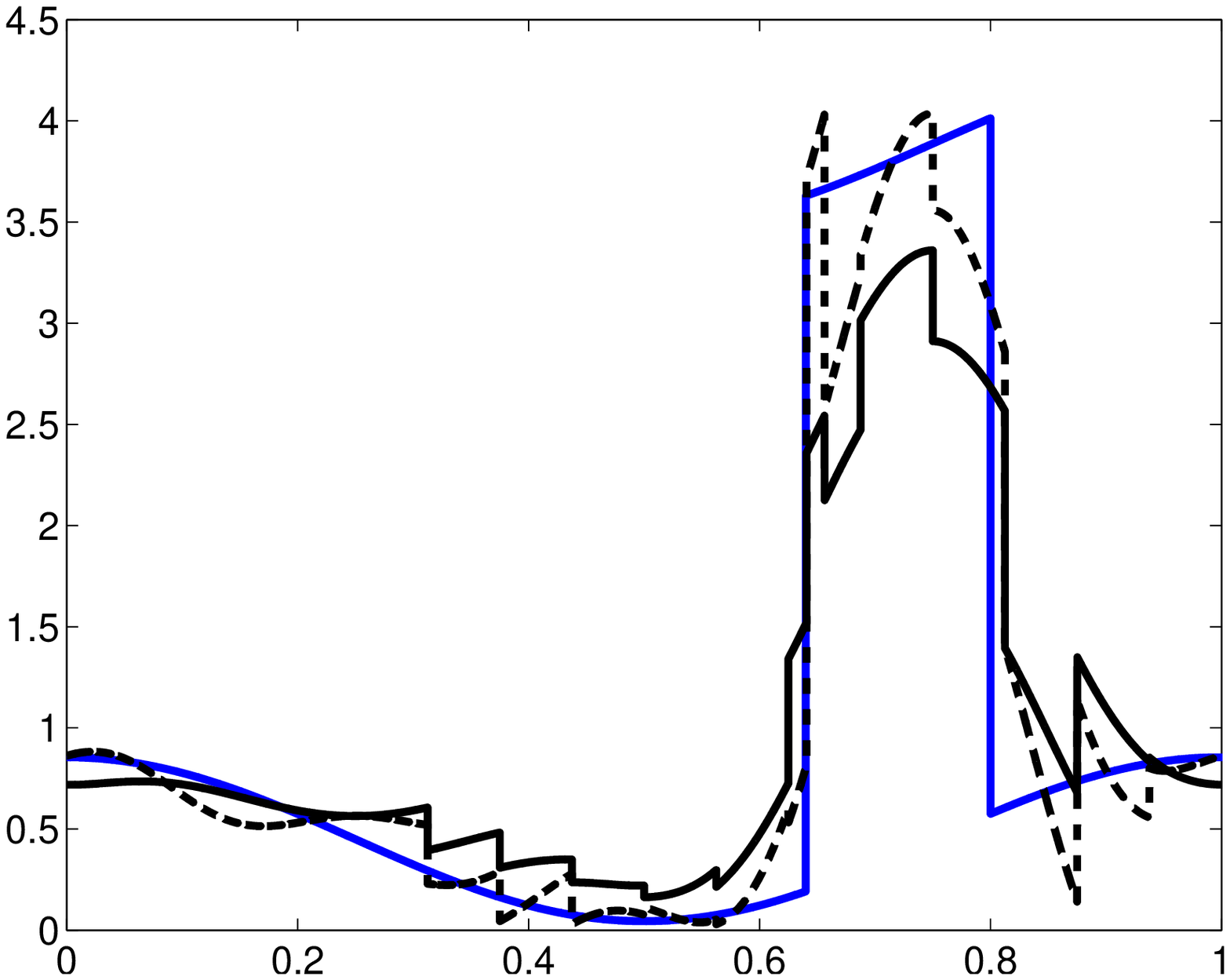} &  \parinc{\widthplot1}{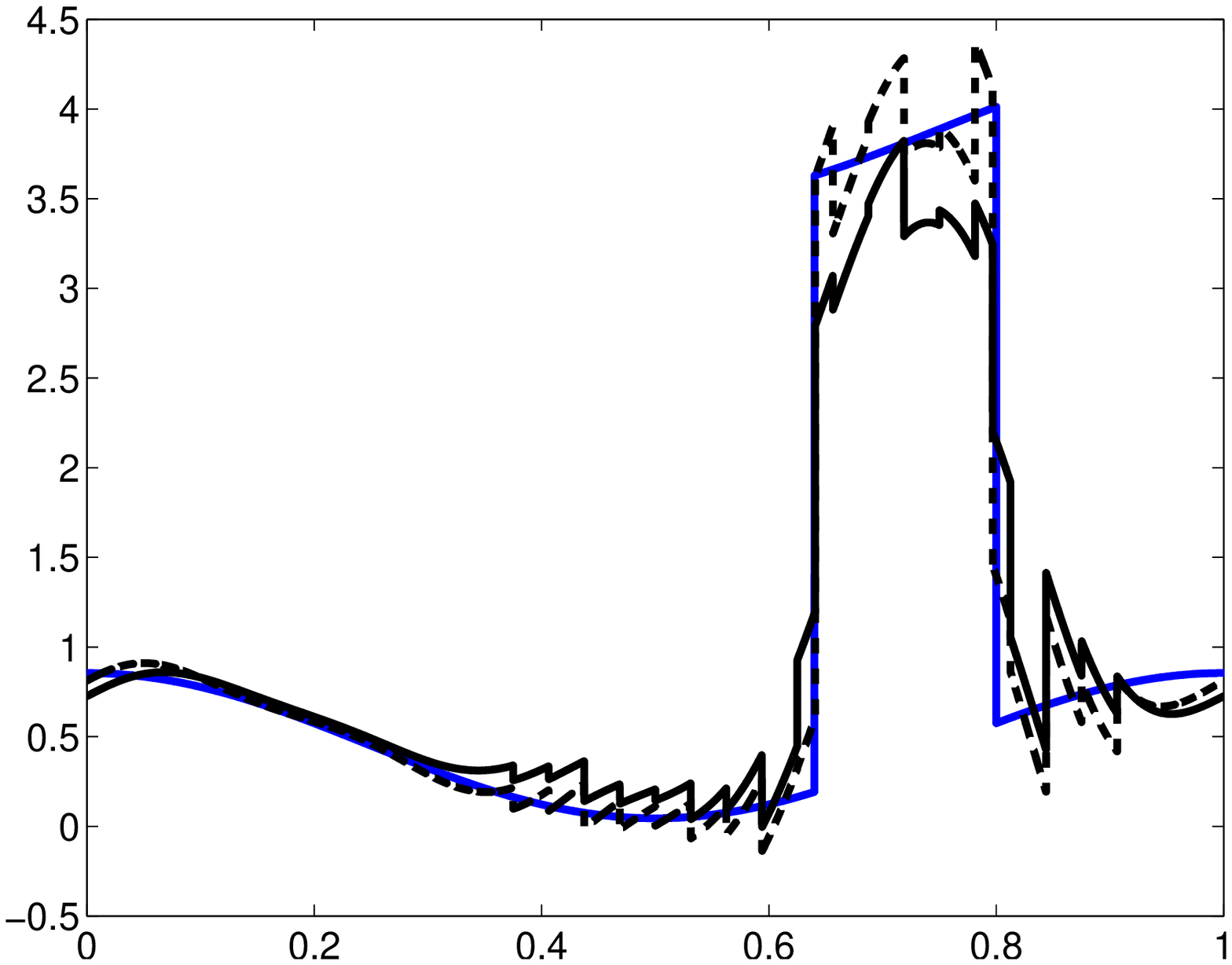}
  \end{tabular}
  \caption{The different densities and their ``Mix2''
    estimates. Densities are plotted in blue while their estimates
    are plotted in black. The full line corresponds to the
   adaptive Dantzig studied in this paper while the dotted line
    corresponds to its least square variant.}
  \label{fig:plots}
\end{figure}

We test the estimators of the following 4 functions shown in
Figure~\ref{fig:plots} (with their Dantzig and Dantzig+Least Square estimates with the ``Mix2'' dictionary):
\begin{itemize}
\item a very spiky density
\[
f_1(t) = .47 \times \left( 4t\times\indic_{t\leq.5} +
  4(1-t)\times\indic_{t>.5}\right) + .53 \times \left( 75 \times\indic_{.5\leq
  t \leq .5+\frac{1}{75}} \right),
\]
\item a mix of Gaussian and Laplacian type densities
\[
f_2(t)  = .45 \times \left( \frac{e^{-(t-.45)^2/(2(.125)^2)}}
{\int_0^1 e^{-(u-.45)^2/(2(.125)^2)} du}
 \right)
+ .55 \times \left( \frac{e^{20|t-.67|}}{\int_0^1e^{20|u-.67|} du}
\right),
\]
\item a mix of uniform densities on subintervals
\[
f_3(t)  = .25 \times \left( \frac{1}{.14} \indic_{.33 \leq t \leq .47}
\right)
+  .75 \times \left( \frac{1}{.16} \indic_{.64 \leq t \leq .80}
\right),
\]
\item a mix of a density easily described in the Fourier domain and
a uniform density on a subinterval
\[
f_4(t)  = .45 \times \left( 1+ .9 \cos (2\pi t) \right)
+  .55 \times \left( \frac{1}{.16} \indic_{.64 \leq t \leq .80}
\right).
\]
\end{itemize}

\begin{figure}
\def\widthboxplot1{4.75cm}
  \centering
  \begin{tabular}{cccc}
& Dantzig/Lasso & Dantzig/Non adapt. Dantzig & Dantzig/Dantzig+LS \\
$f_1$ & \parinc{\widthboxplot1}{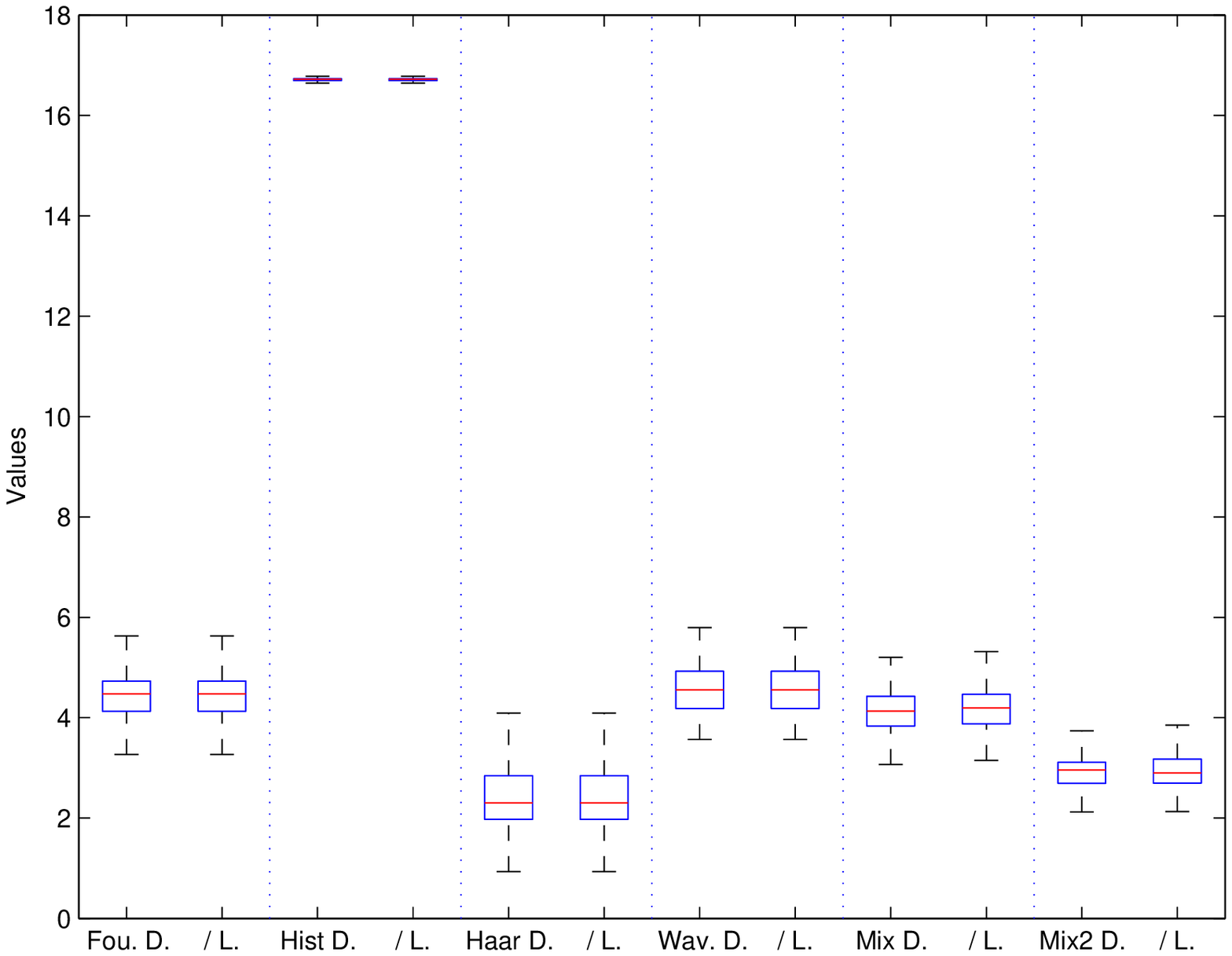} &
\parinc{\widthboxplot1}{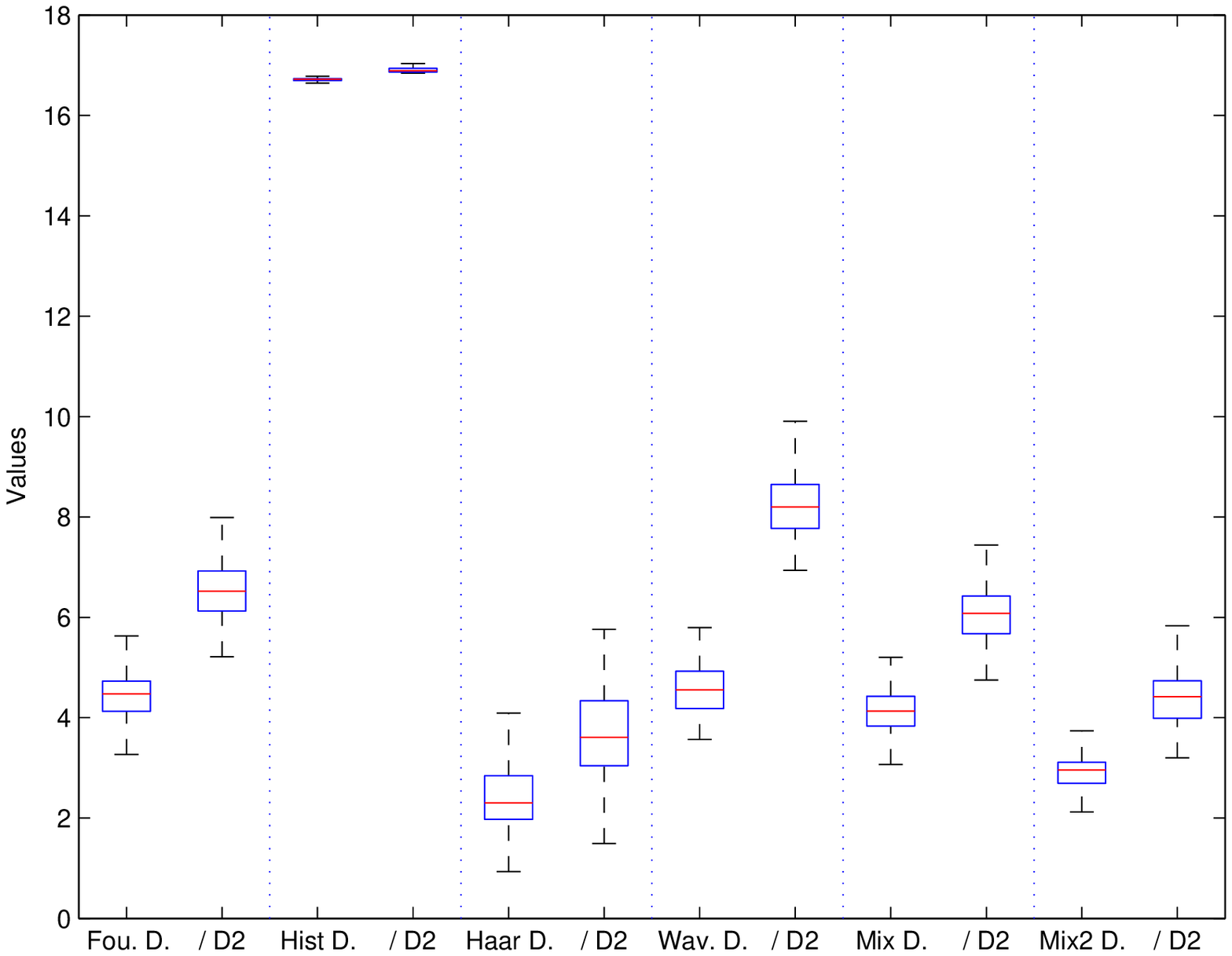} &
\parinc{\widthboxplot1}{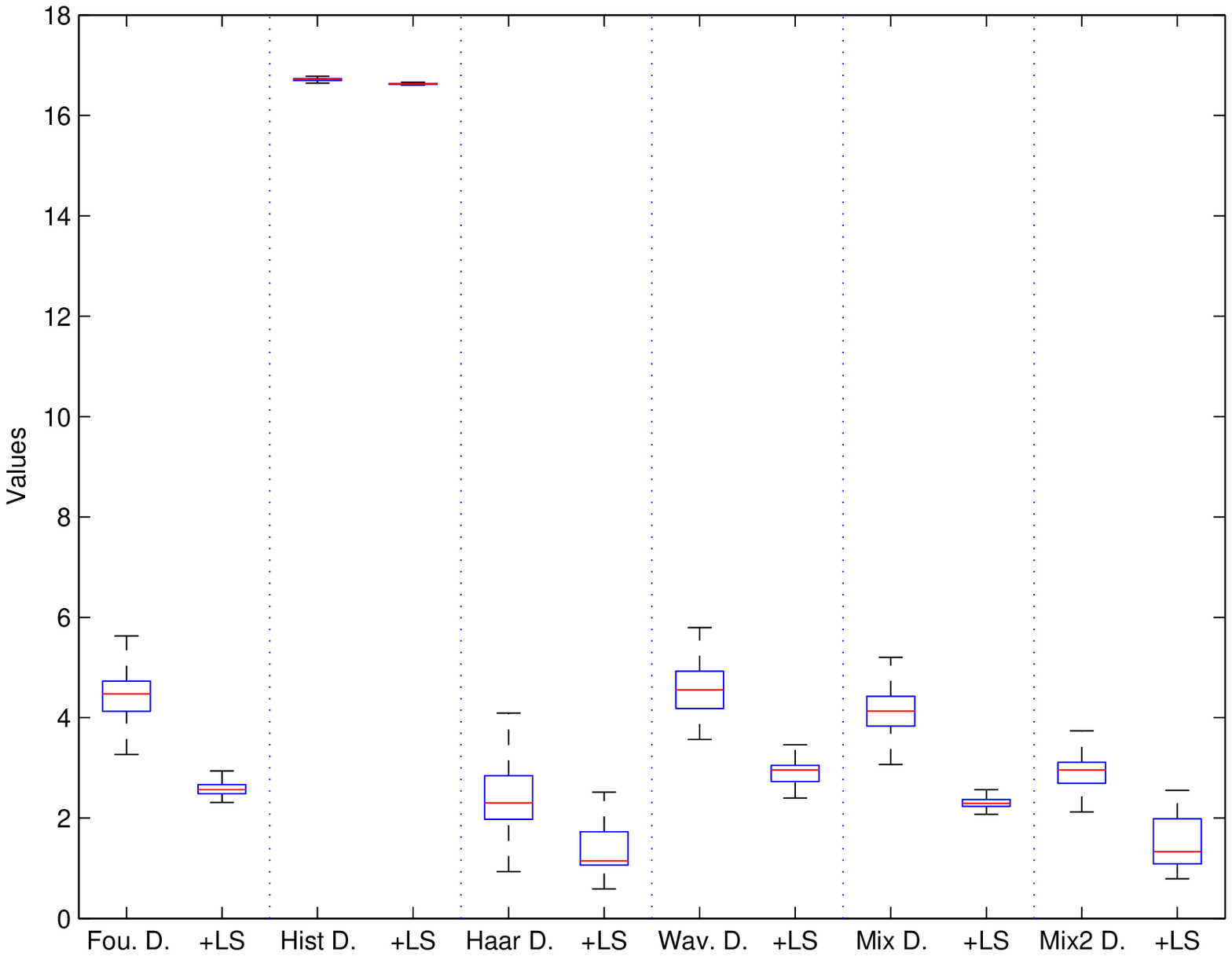}\\
$f_2$ & \parinc{\widthboxplot1}{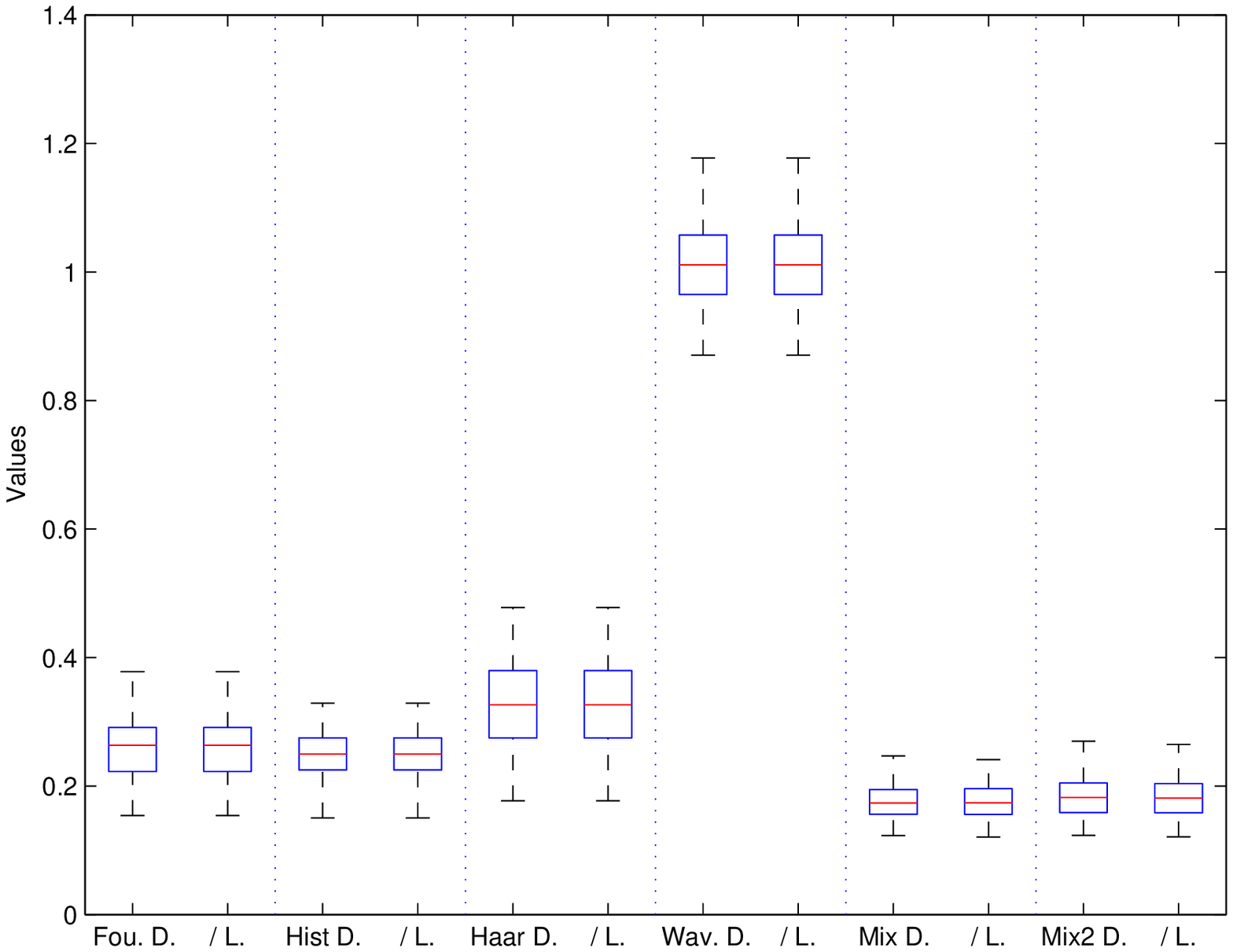} &
\parinc{\widthboxplot1}{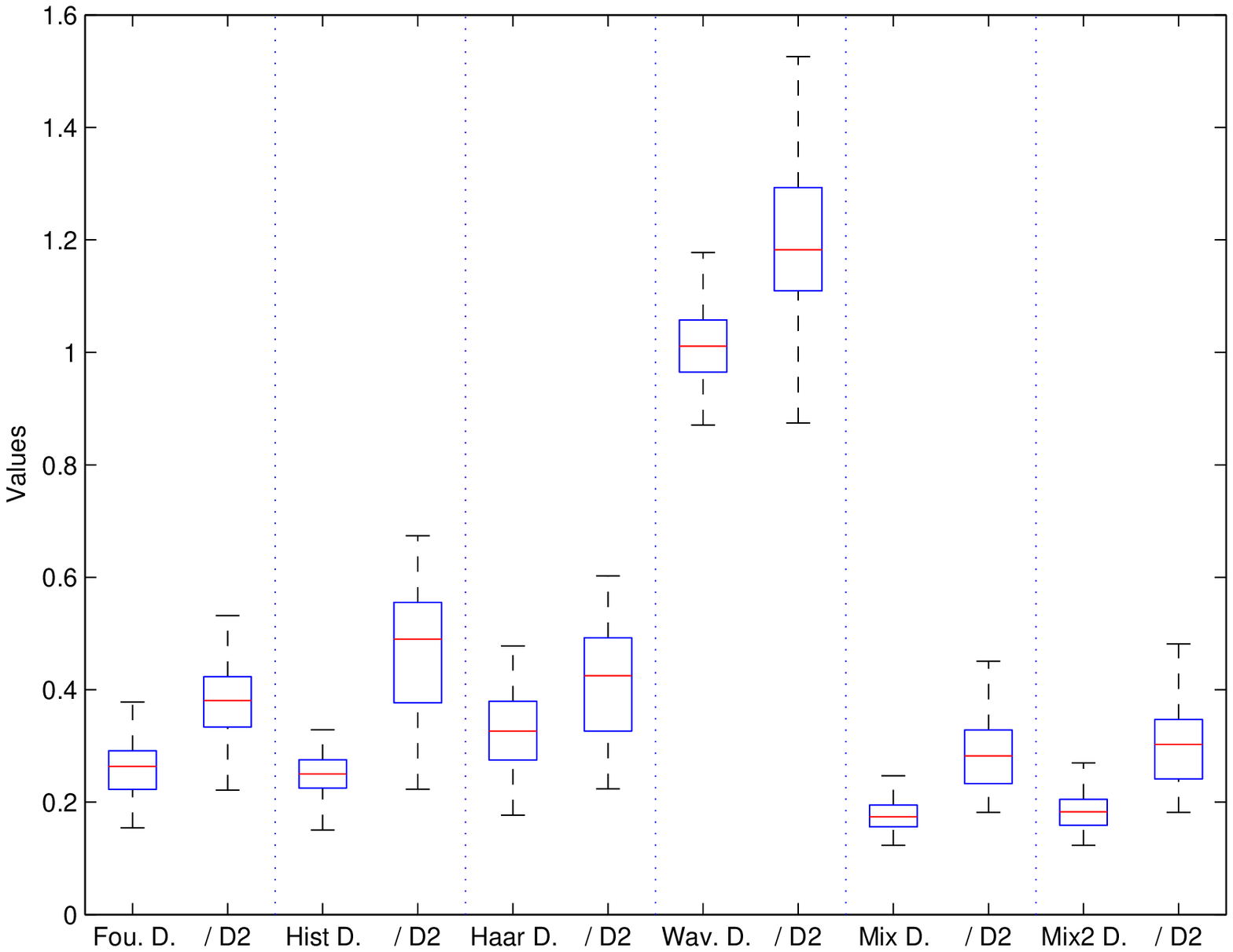} &
\parinc{\widthboxplot1}{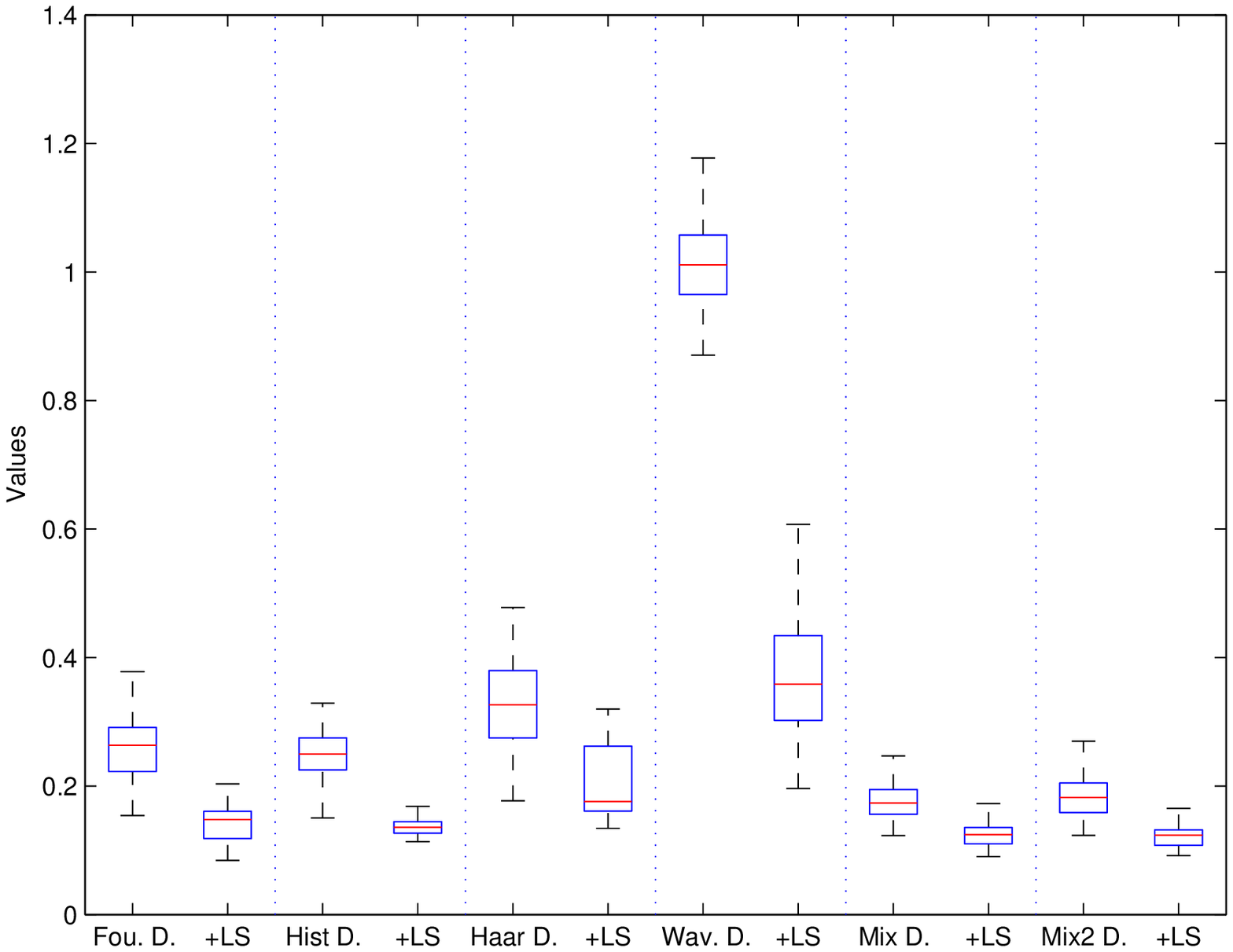}\\
$f_3$ & \parinc{\widthboxplot1}{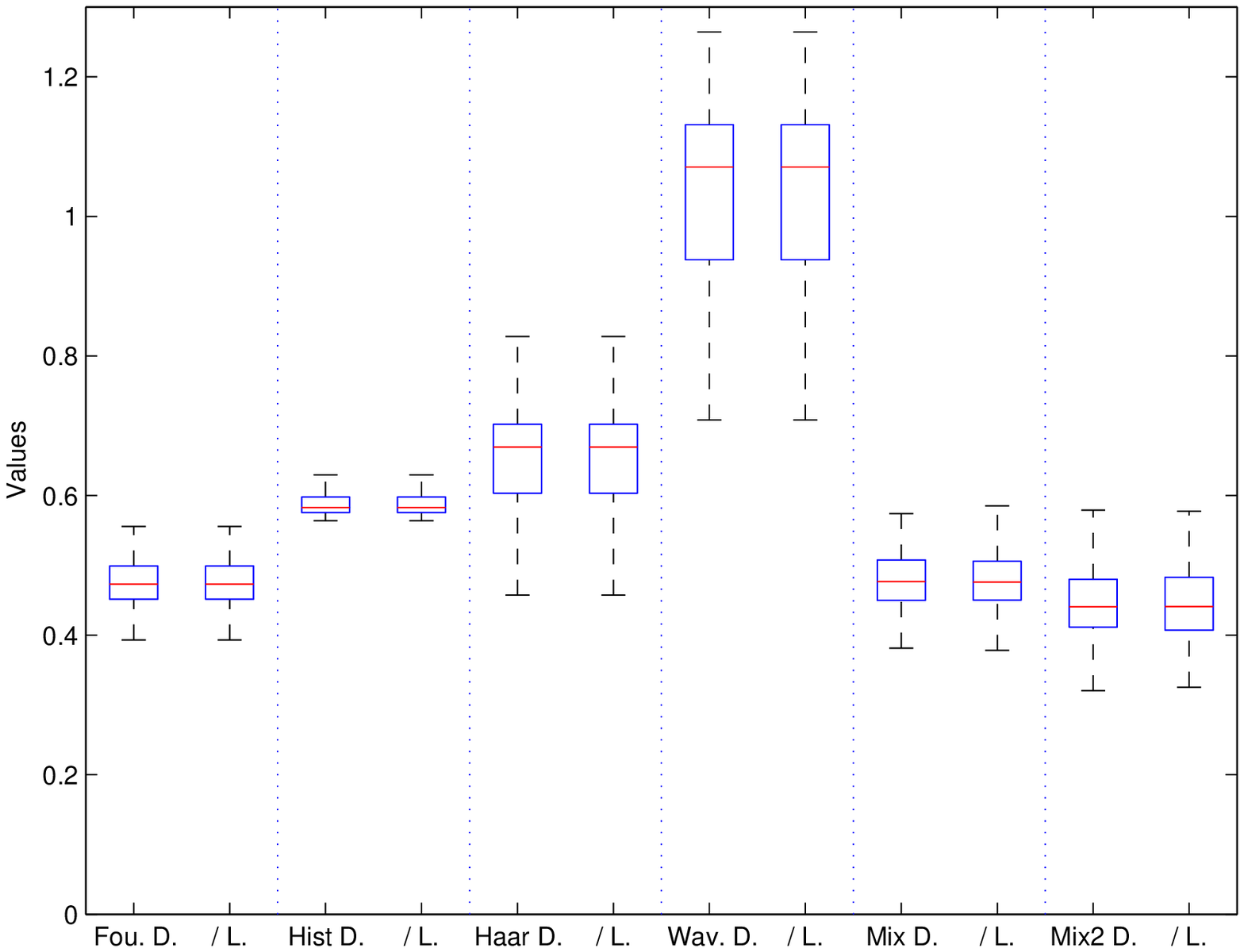} &
\parinc{\widthboxplot1}{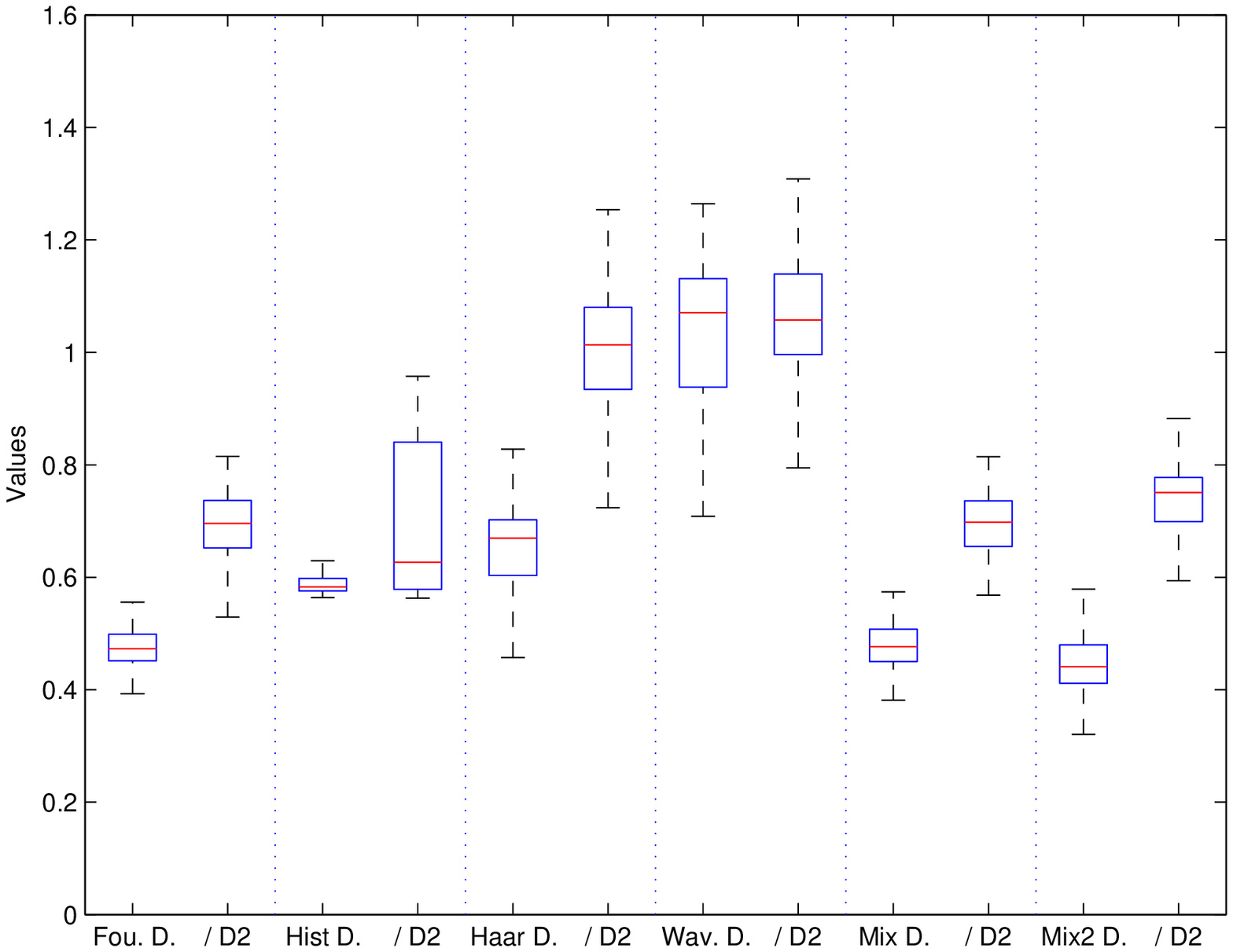} &
\parinc{\widthboxplot1}{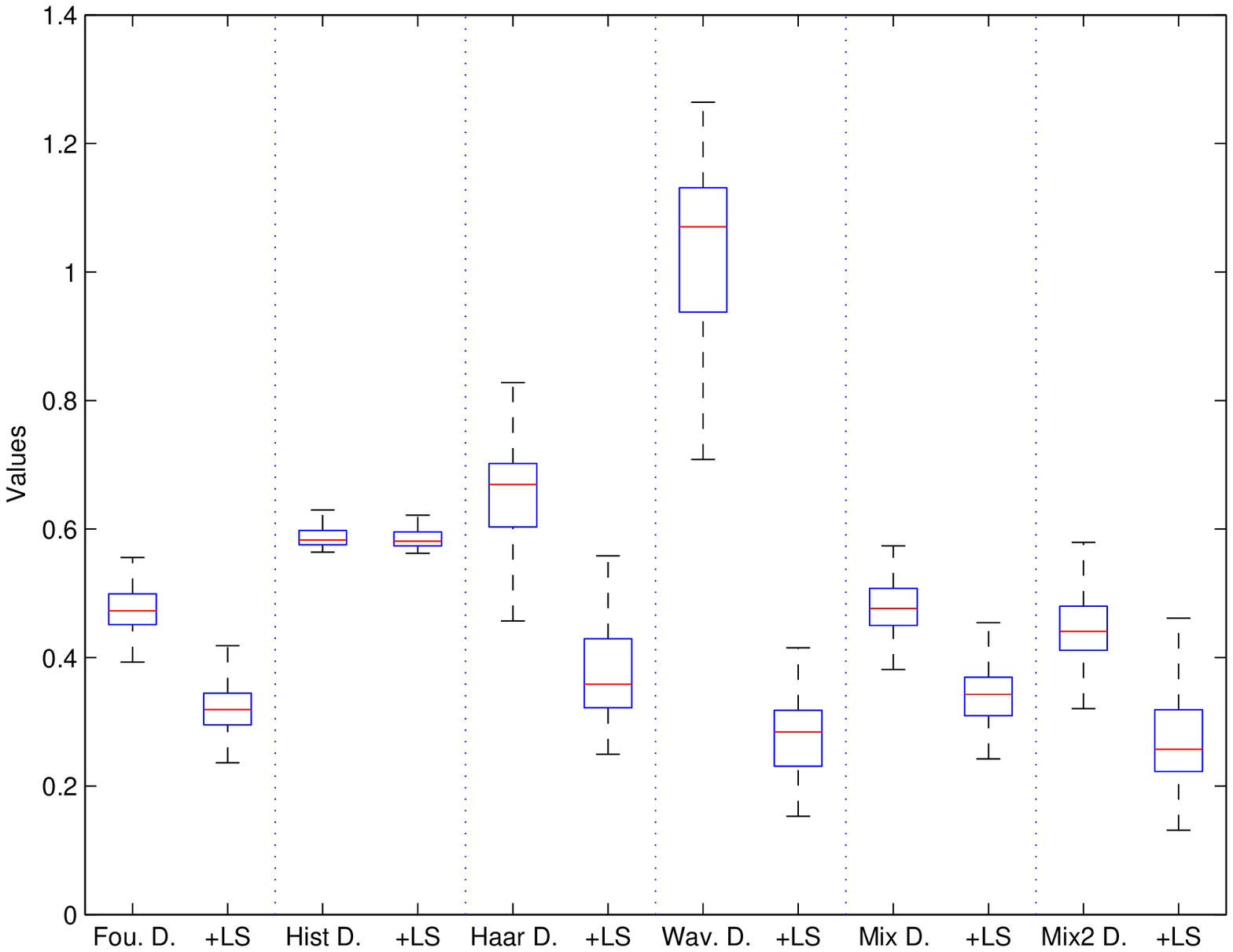}\\
$f_4$ & \parinc{\widthboxplot1}{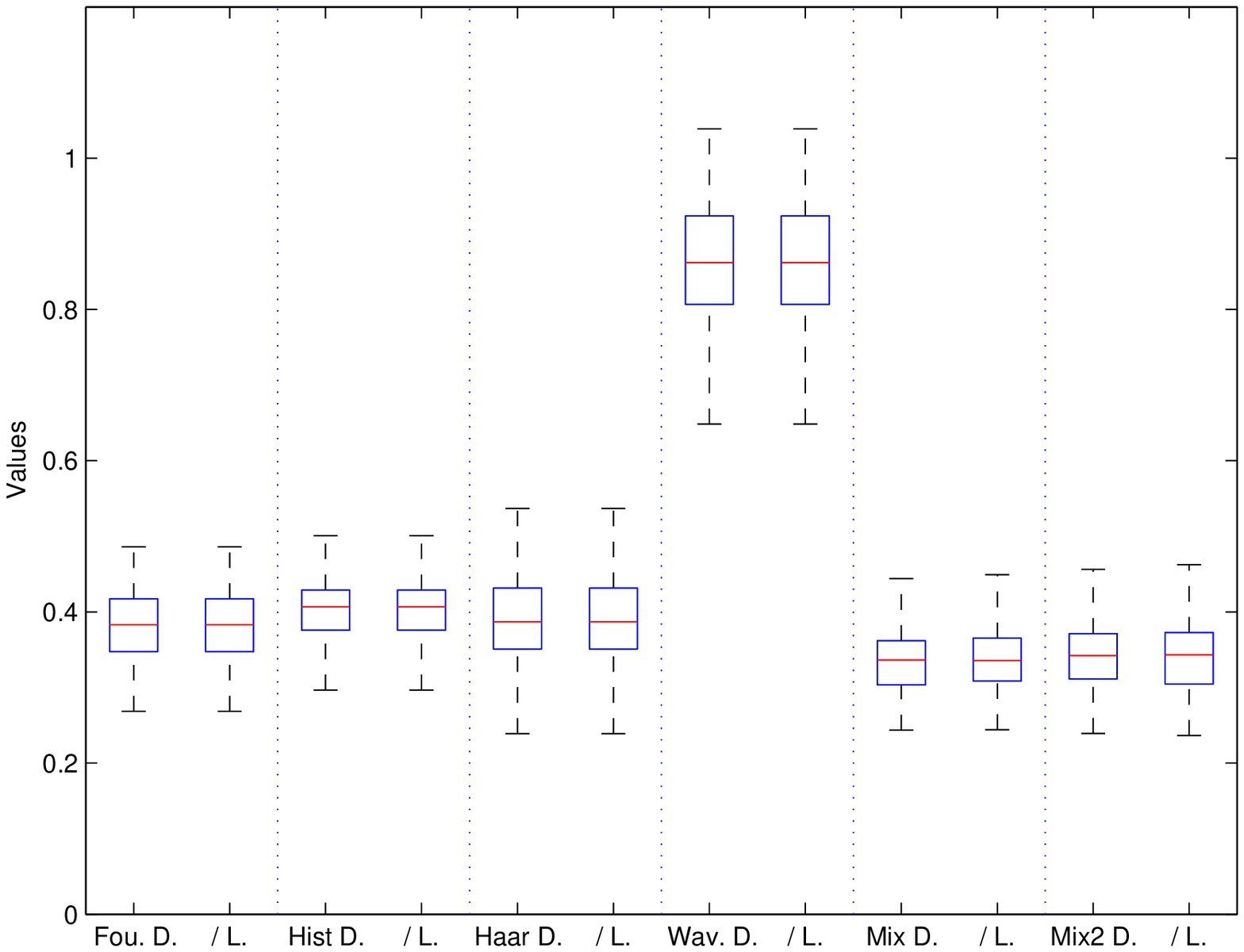} &
\parinc{\widthboxplot1}{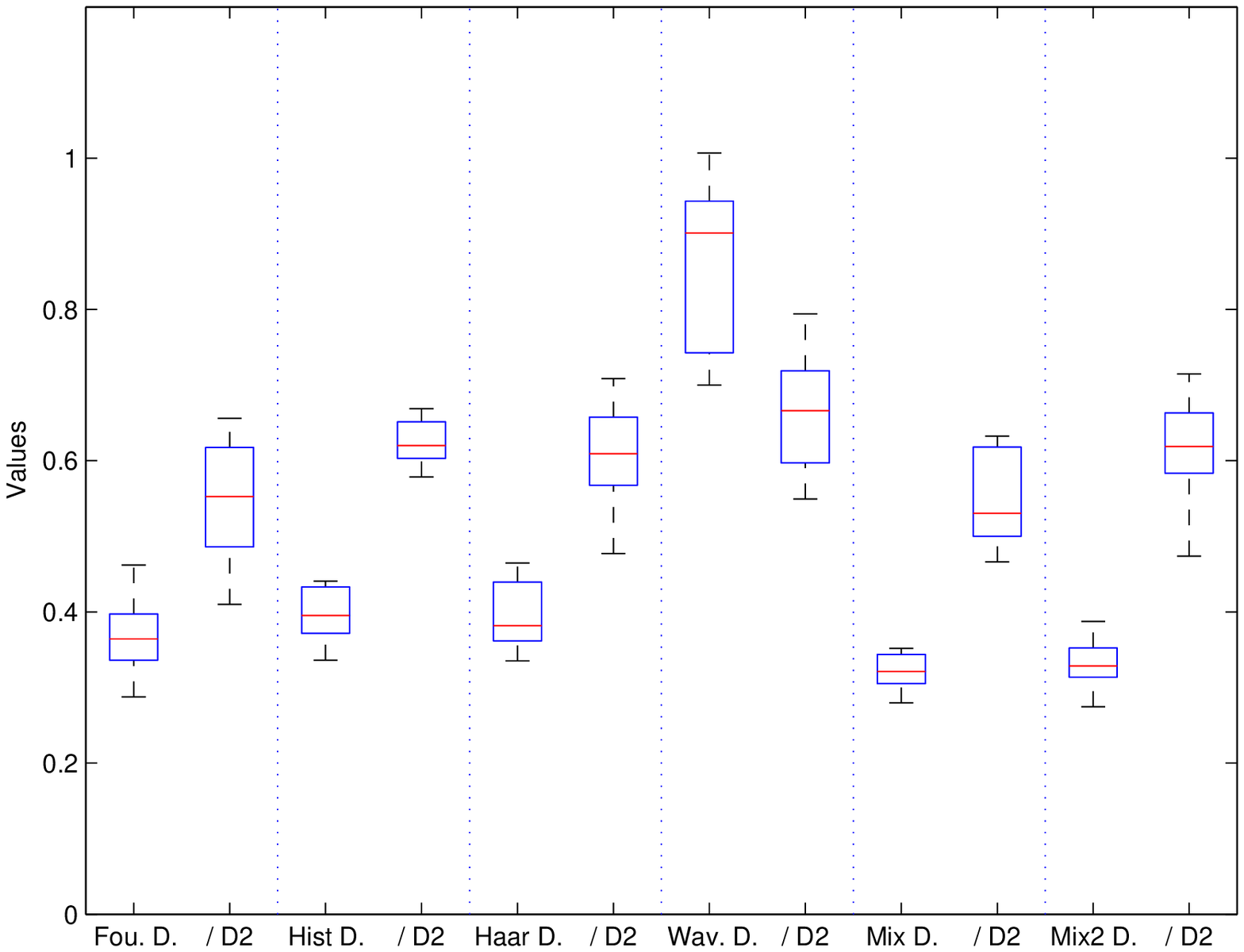} &
\parinc{\widthboxplot1}{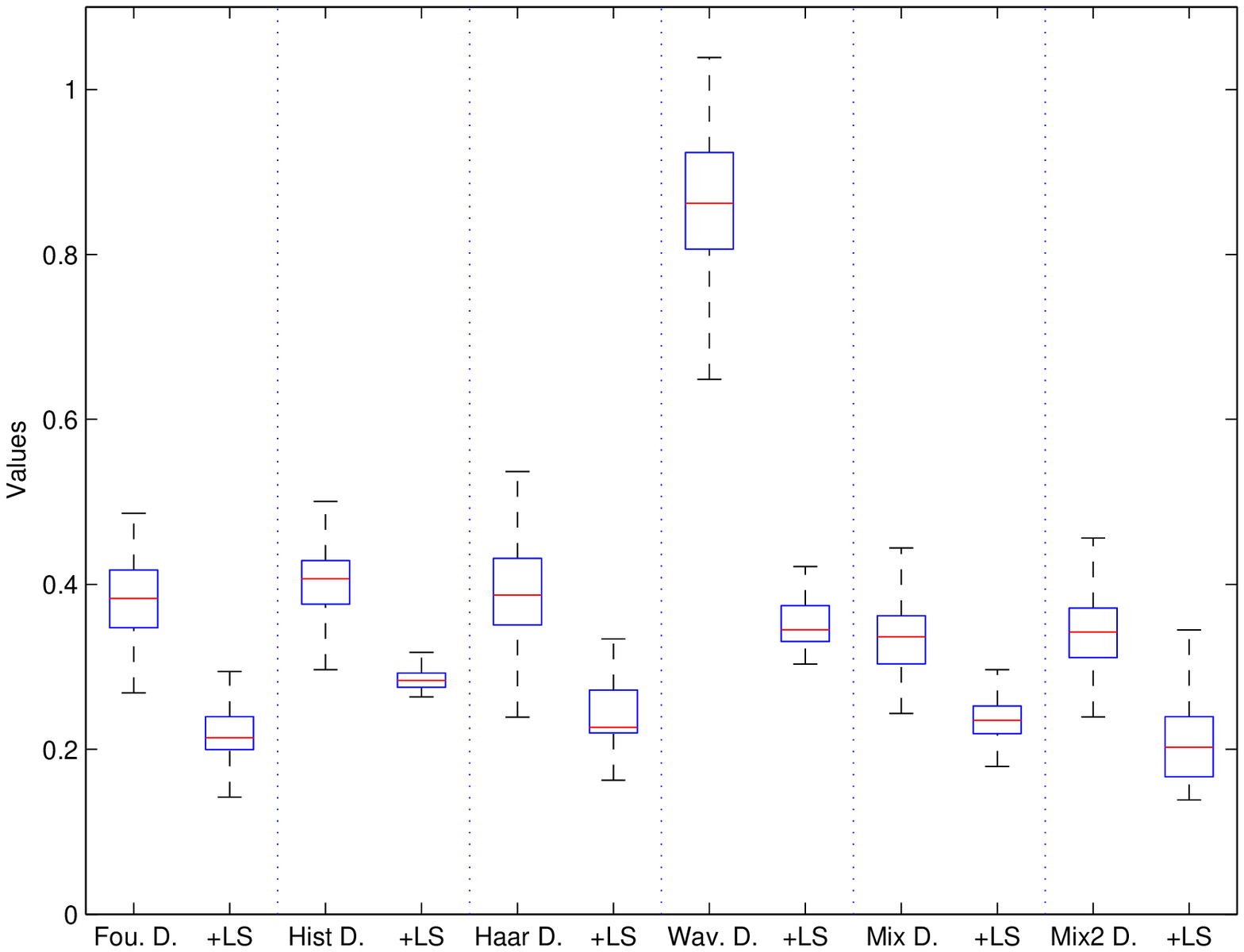}
  \end{tabular}

\caption{Boxplots for $n=500$.  Left column: Dantzig and Lasso estimates. Center column: Dantzig estimates associated with adaptive and non-adaptive constraints.  Right column: Our estimate and the two-step estimate.}
\label{fig:boxplots500}
\end{figure}

\begin{figure}
\def\widthboxplot1{4.75cm}
  \centering
  \begin{tabular}{cccc}
& Dantzig/Lasso & Dantzig/Non adapt. Dantzig & Dantzig/Dantzig+LS \\
$f_1$ & \parinc{\widthboxplot1}{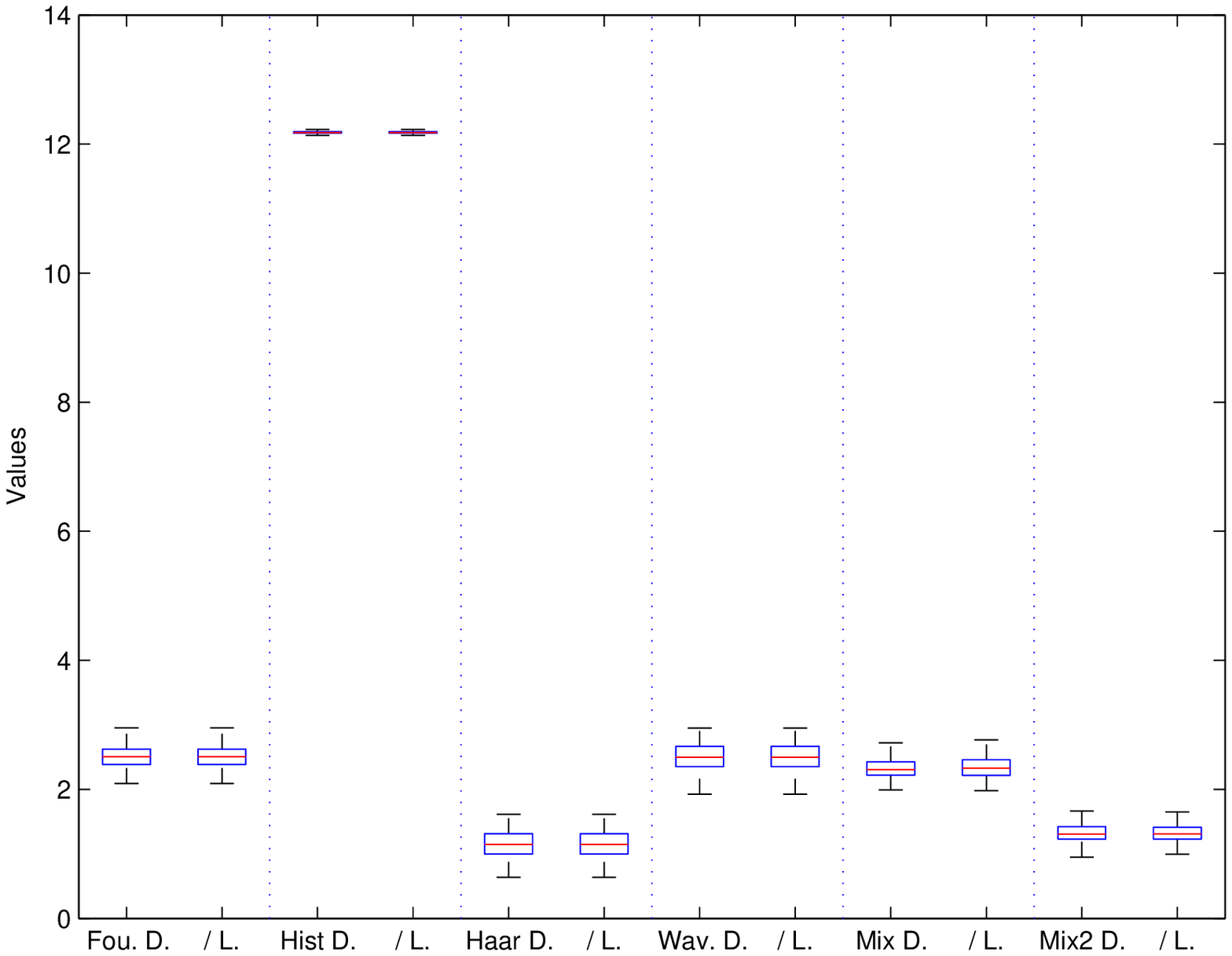} &
\parinc{\widthboxplot1}{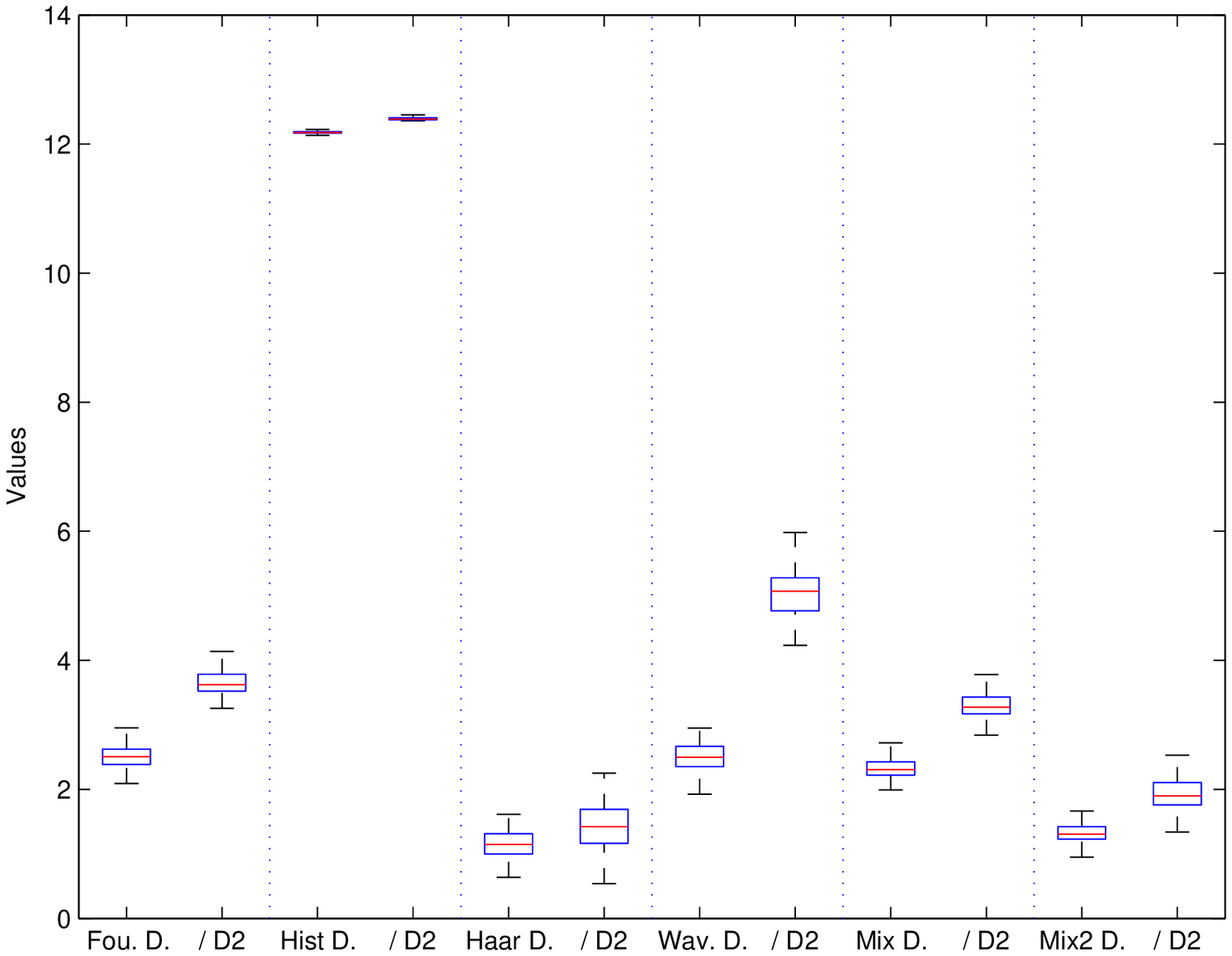} &
\parinc{\widthboxplot1}{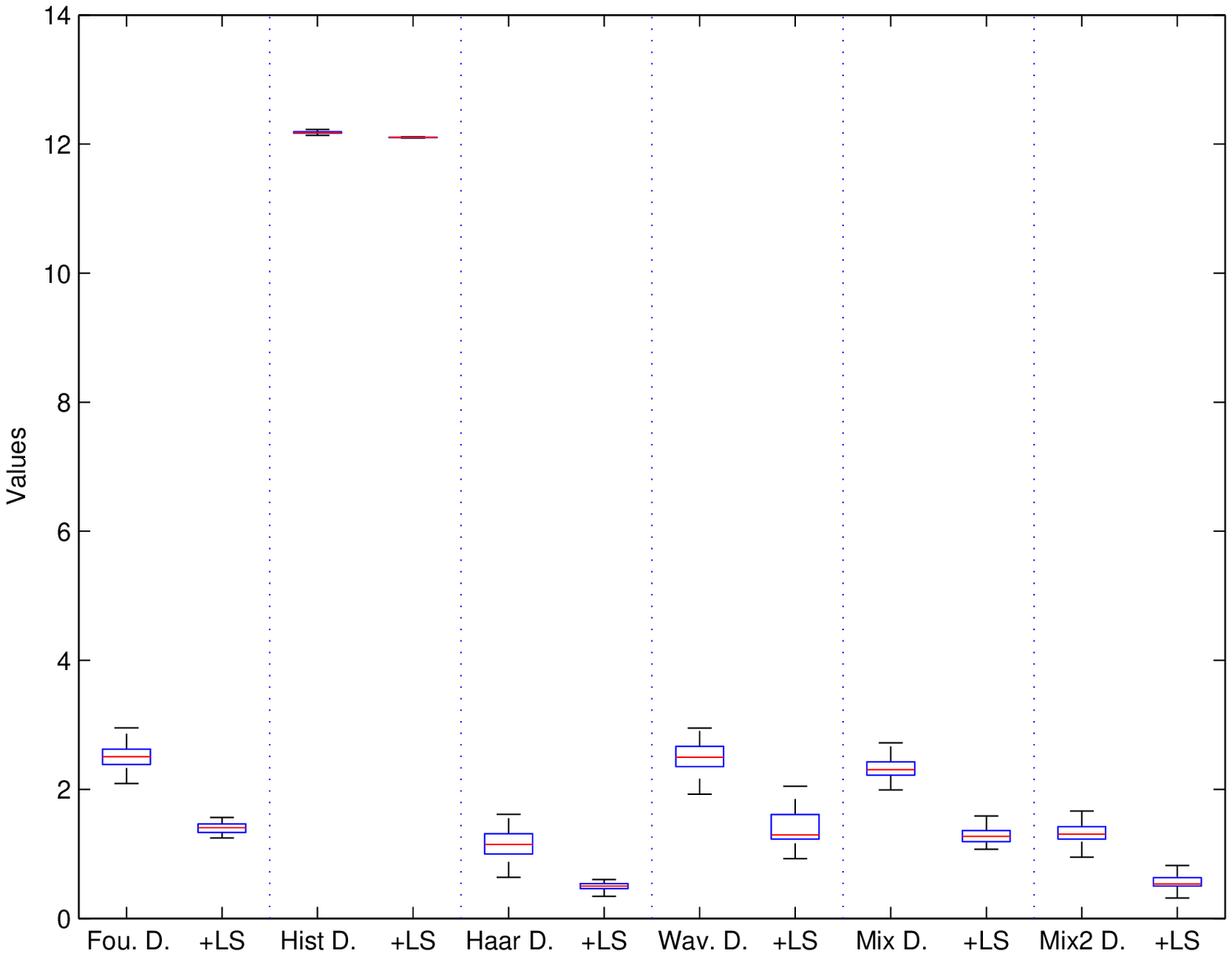}\\
$f_2$ & \parinc{\widthboxplot1}{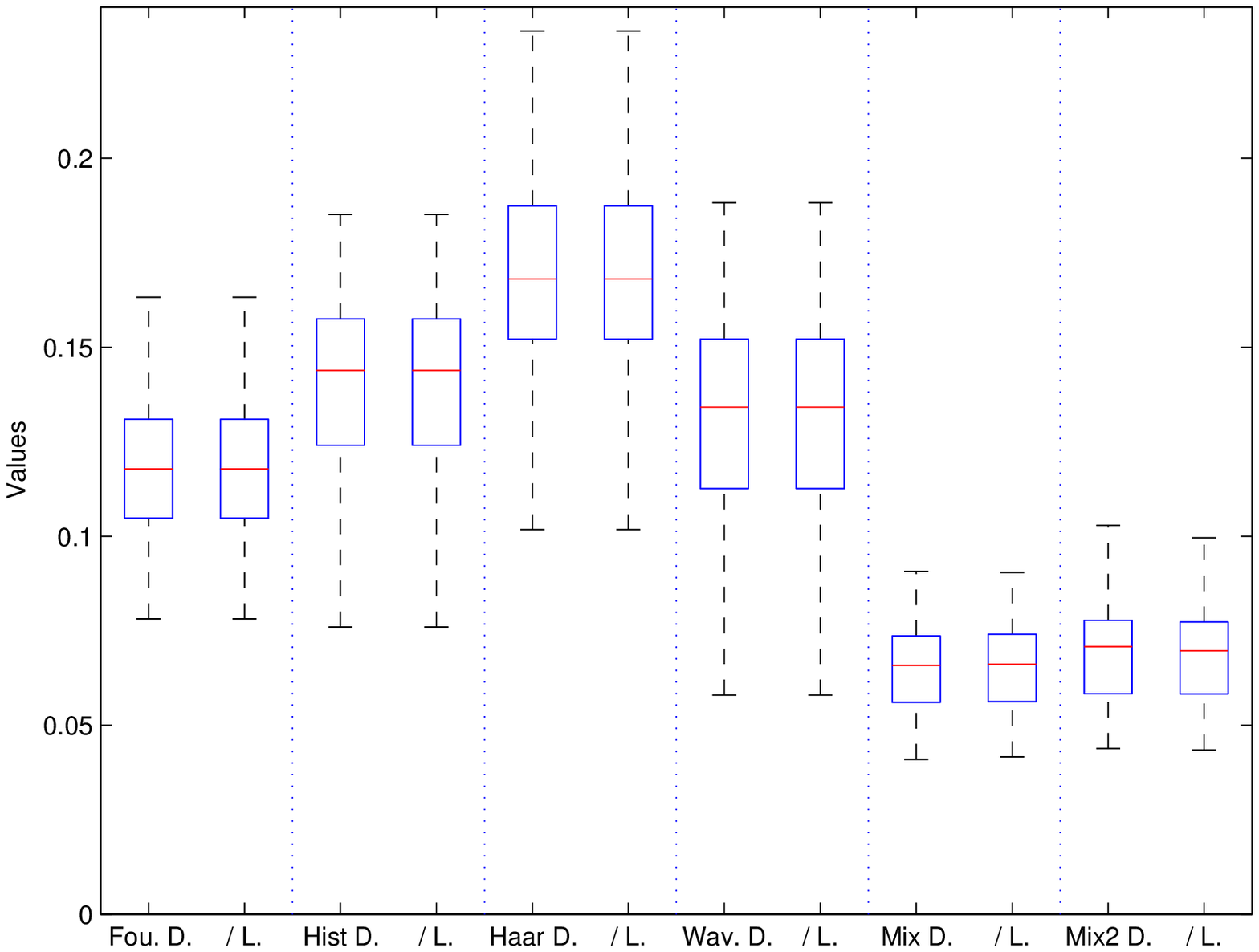} &
\parinc{\widthboxplot1}{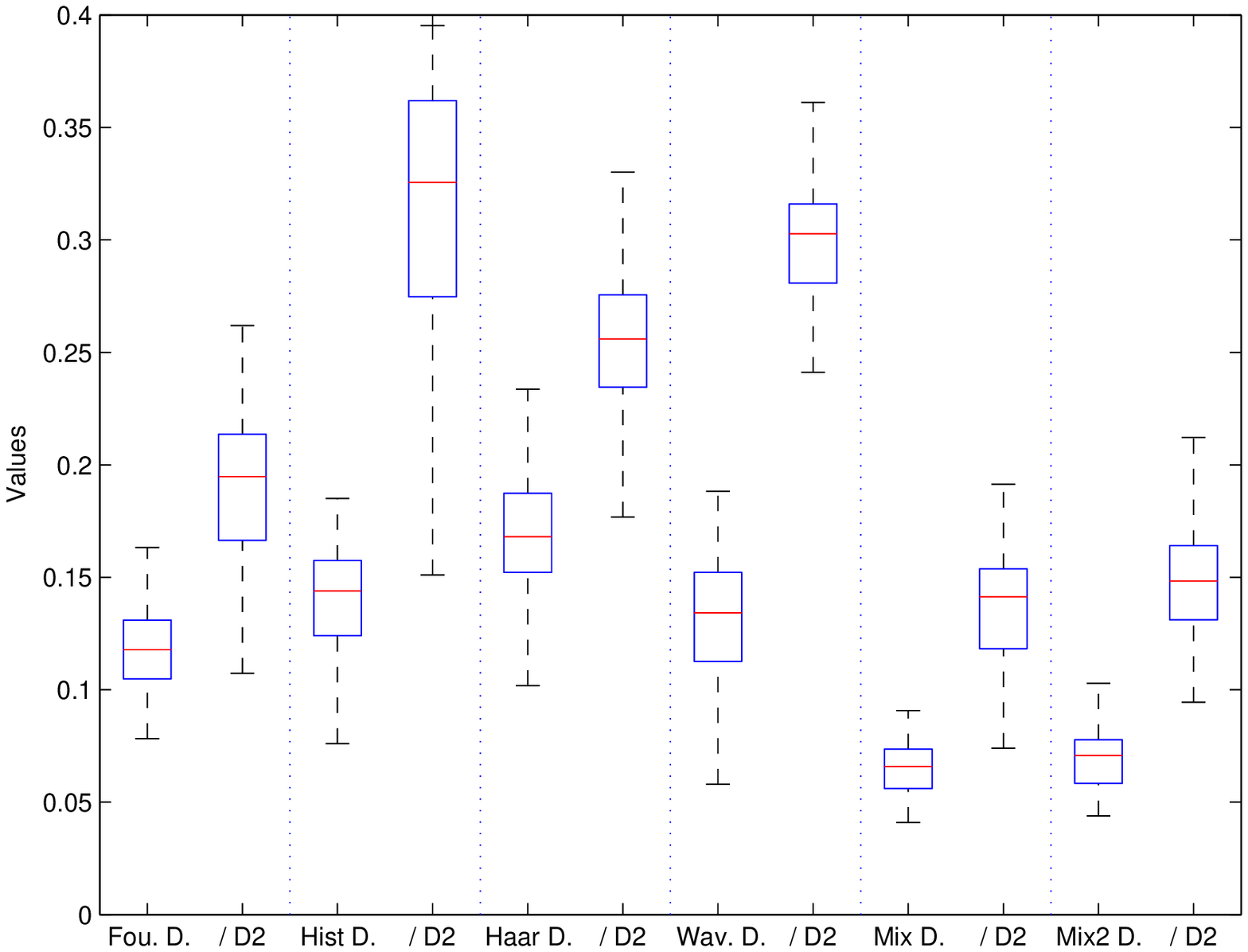} &
\parinc{\widthboxplot1}{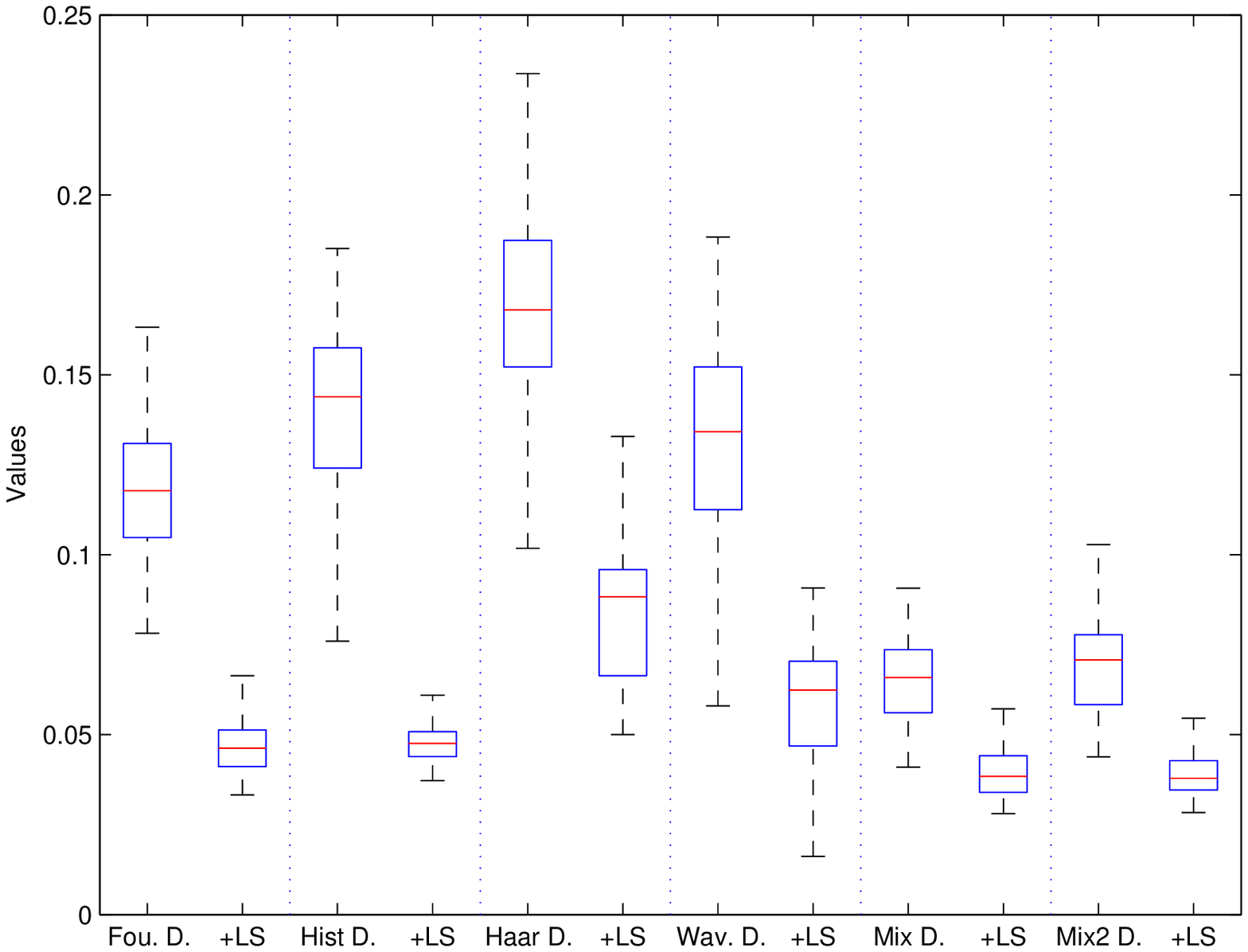}\\
$f_3$ & \parinc{\widthboxplot1}{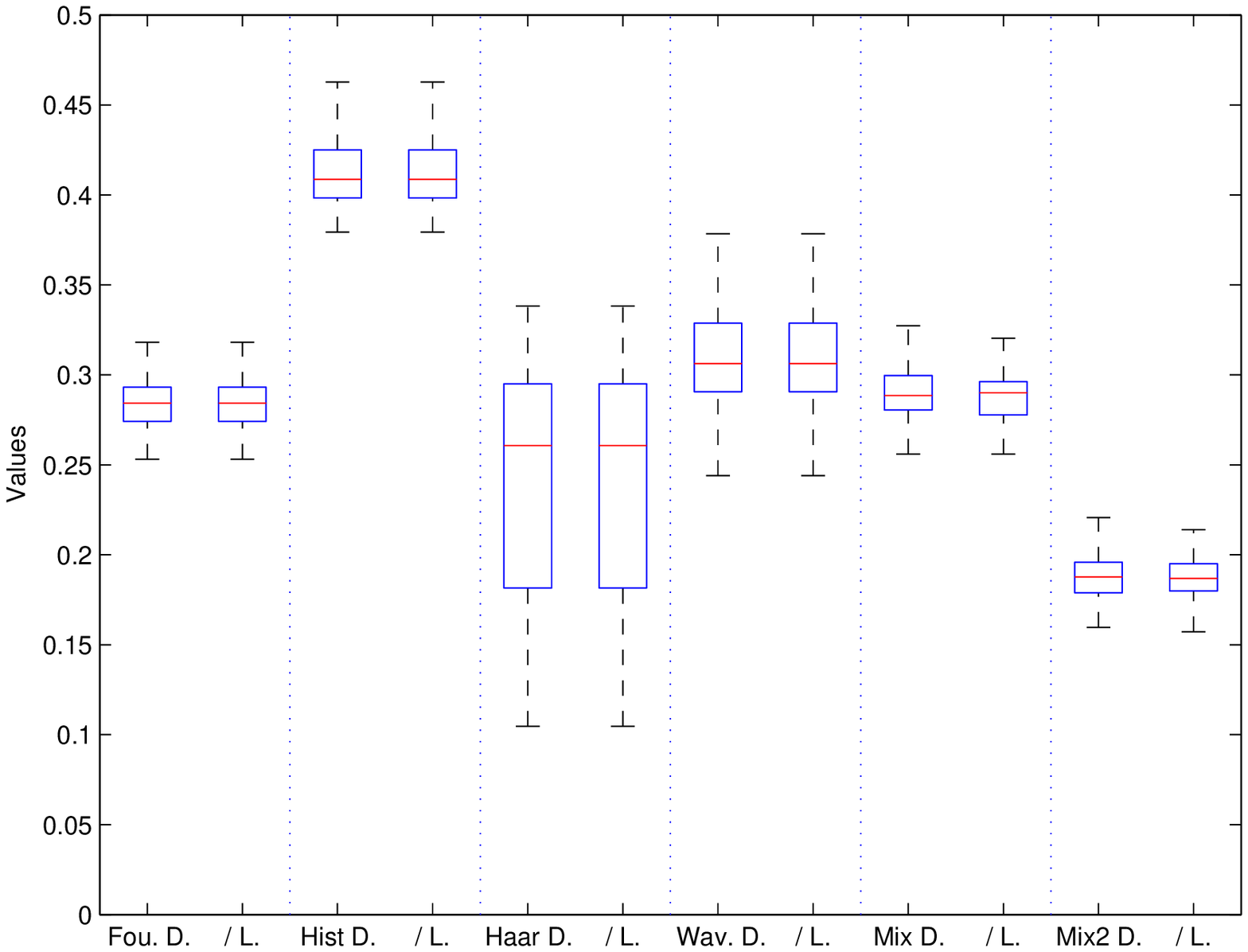} &
\parinc{\widthboxplot1}{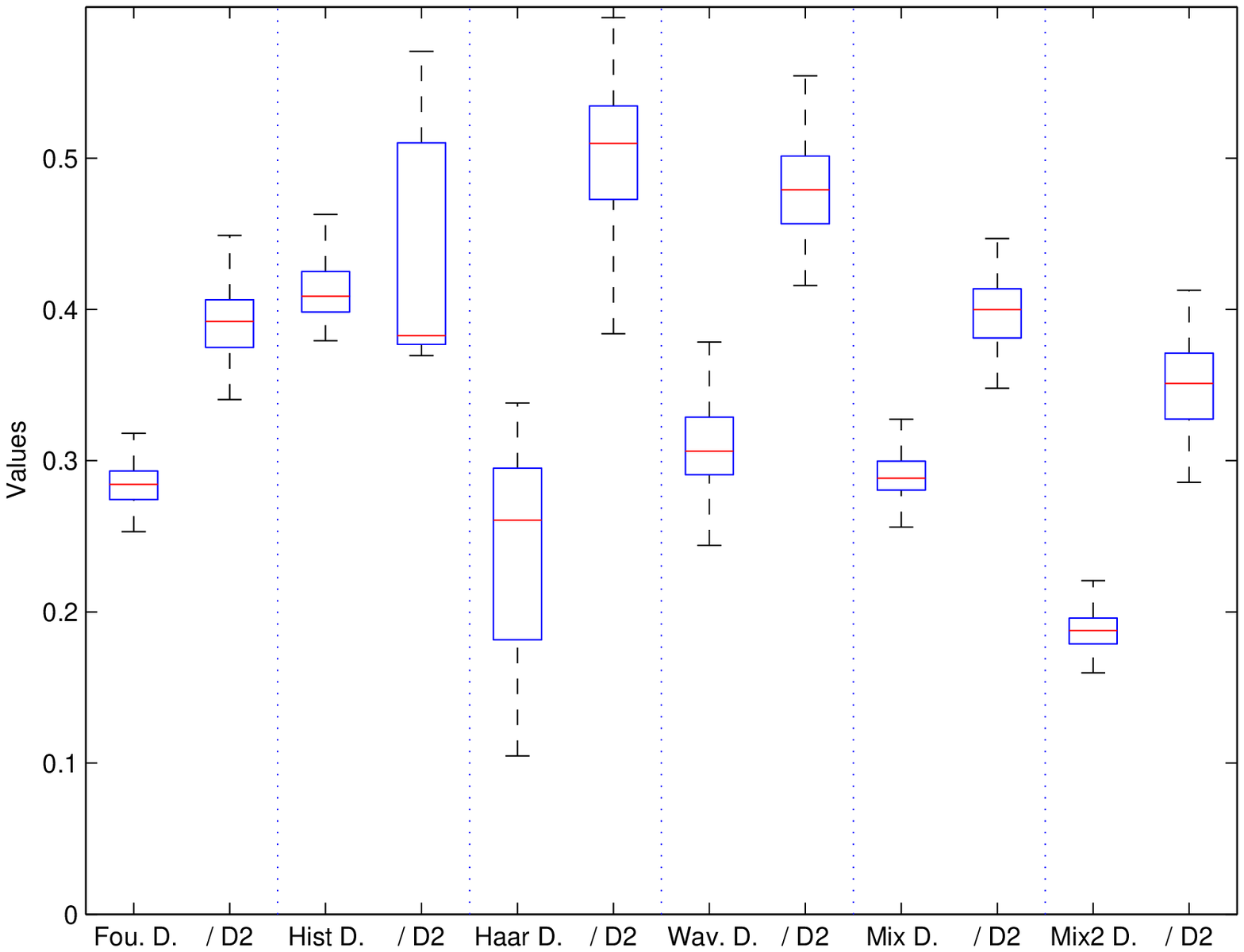} &
\parinc{\widthboxplot1}{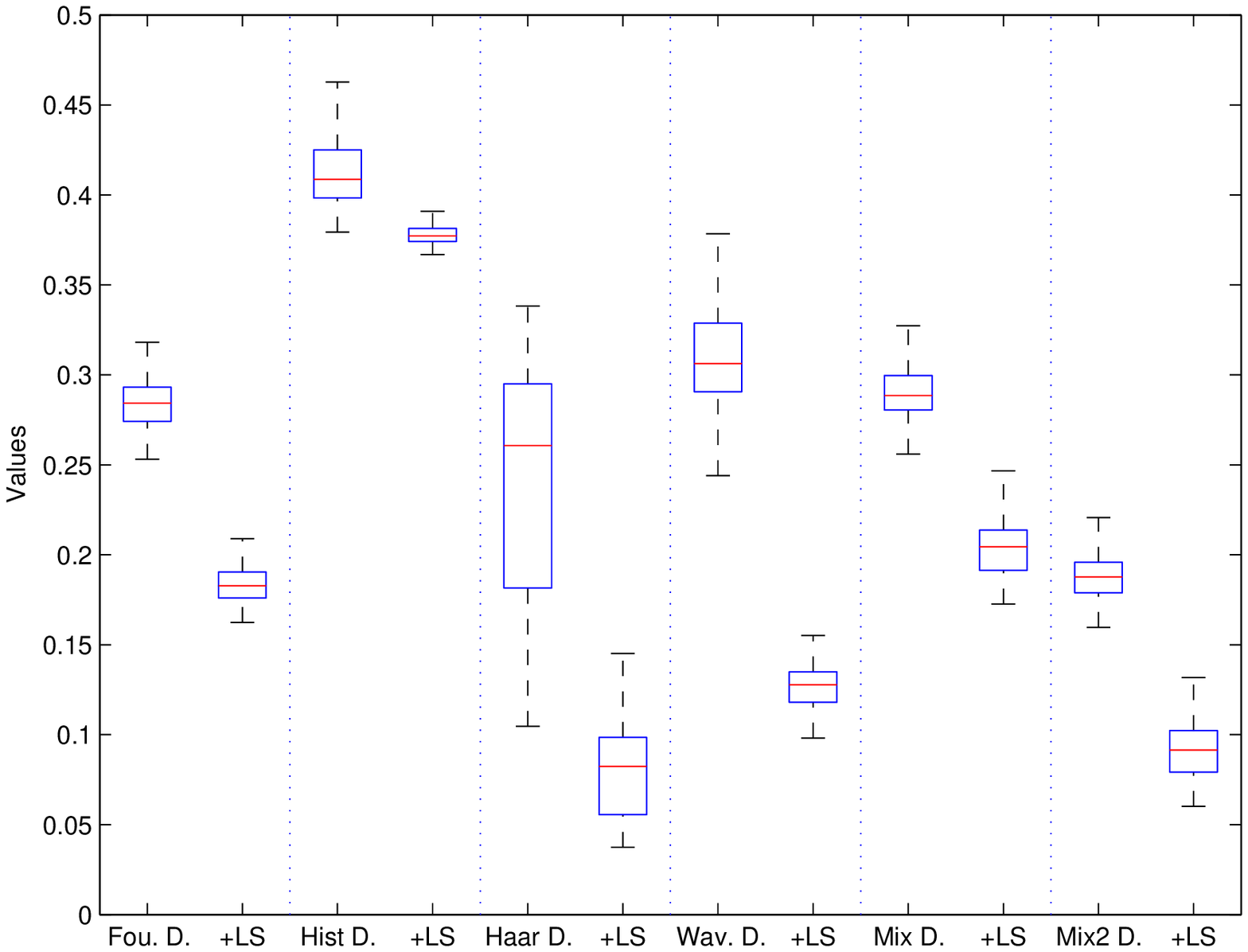}\\
$f_4$ & \parinc{\widthboxplot1}{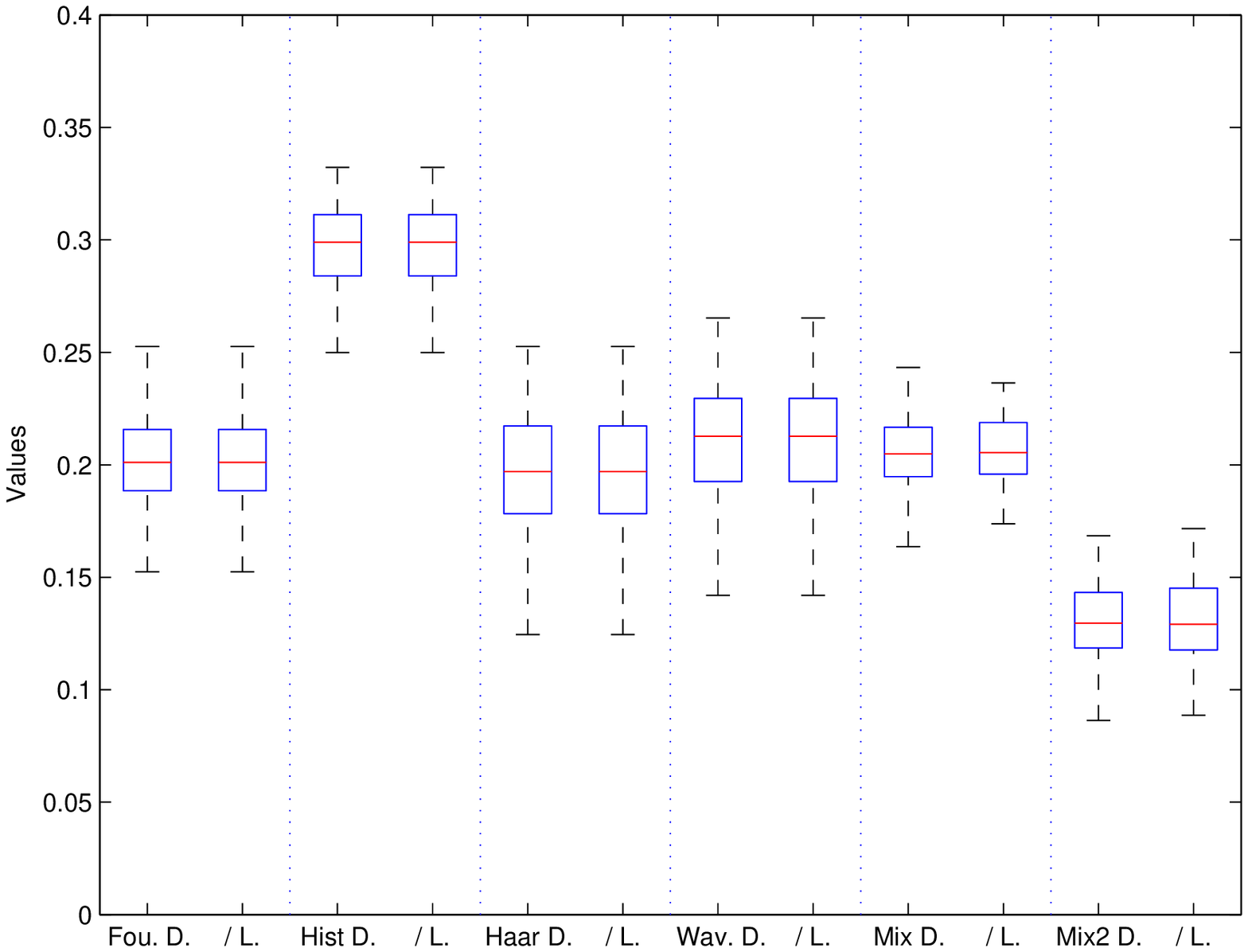} &
\parinc{\widthboxplot1}{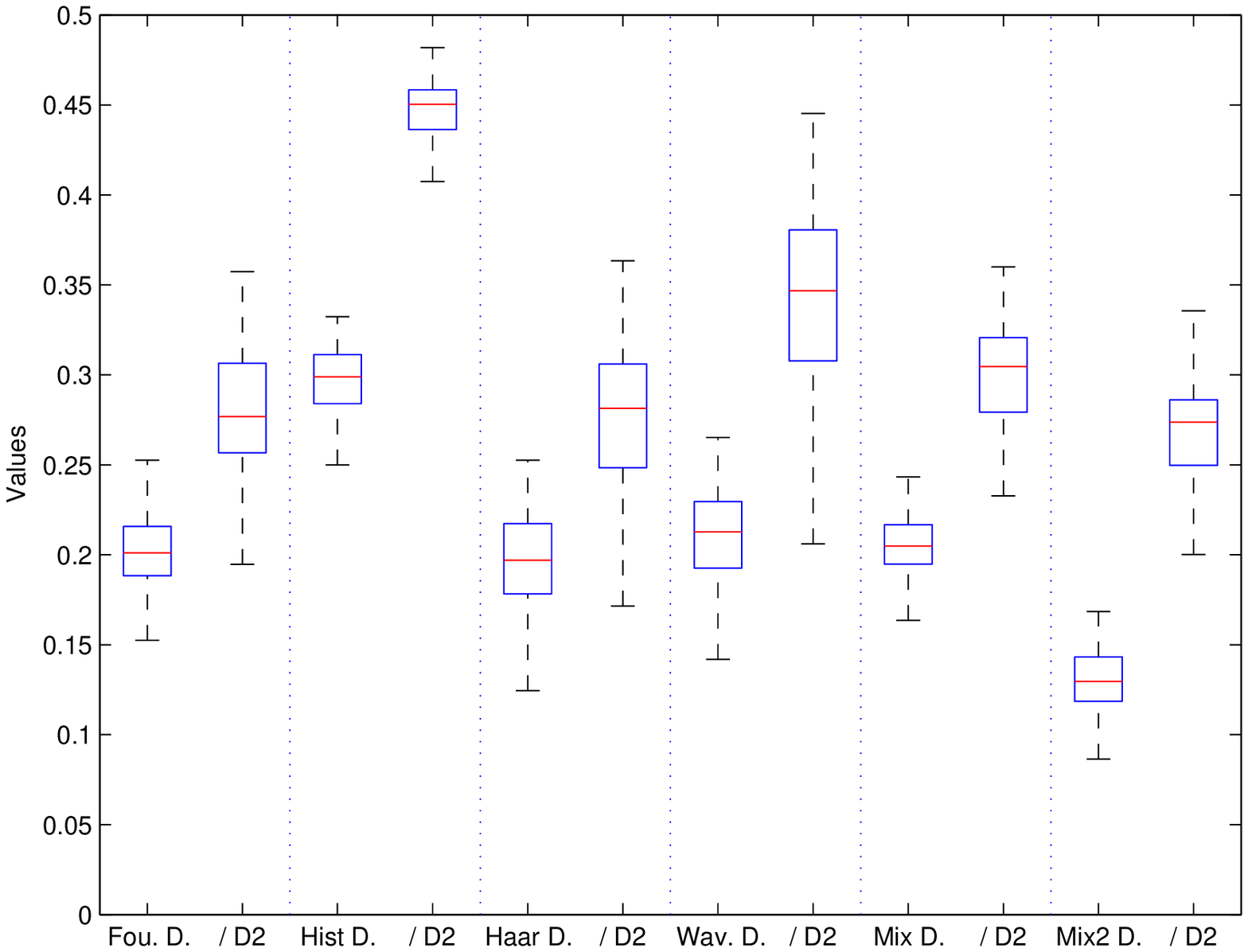} &
\parinc{\widthboxplot1}{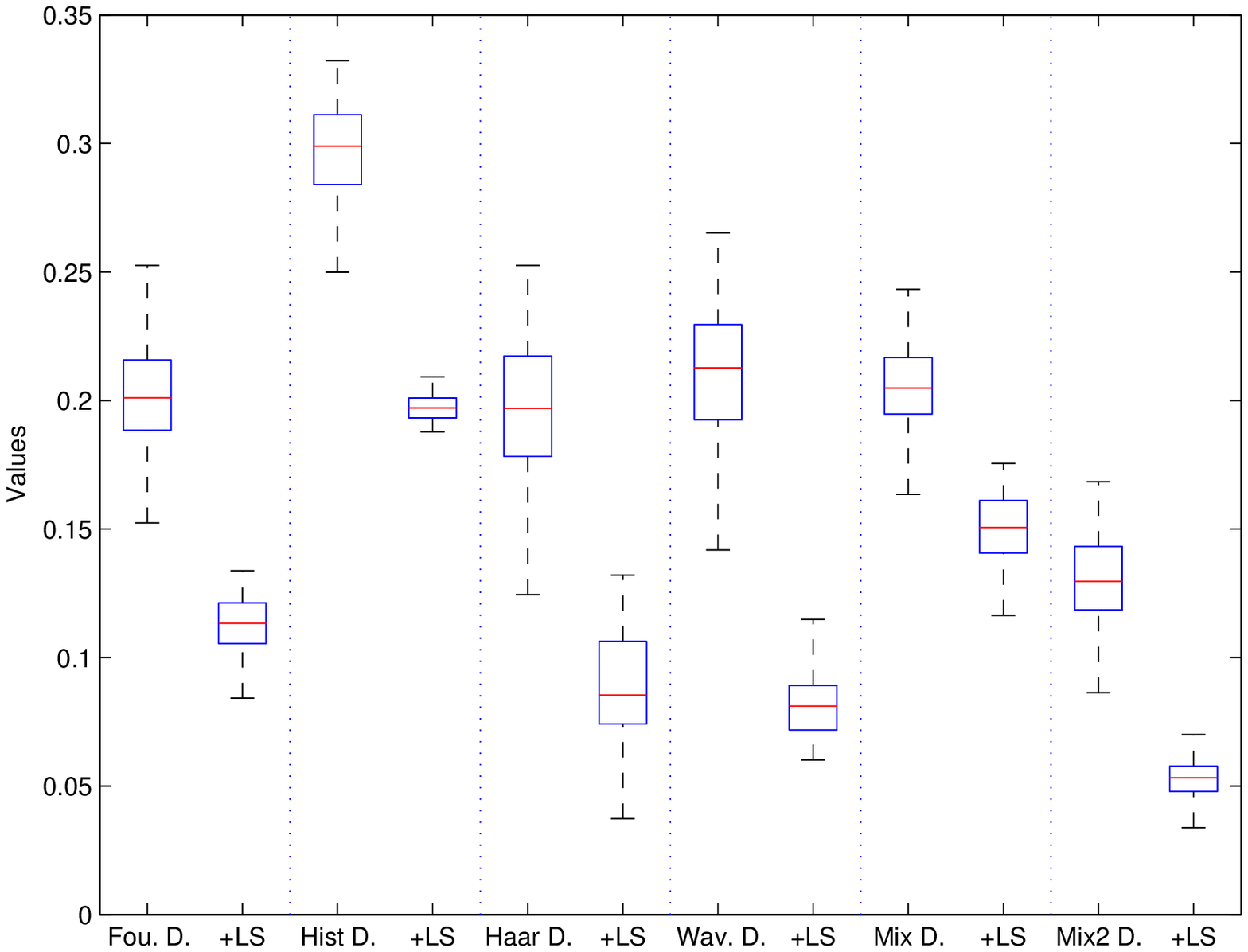}
  \end{tabular}

\caption{Boxplots for $n=2000$.  Left column: Dantzig and Lasso estimates. Center column: Dantzig estimates associated with adaptive and non-adaptive constraints.  Right column: Our estimate and the two-step estimate.}
\label{fig:boxplots2000}
\end{figure}

Boxplots of Figures~\ref{fig:boxplots500} and
\ref{fig:boxplots2000} summarize our numerical experiments for
$n=500$ and $n=2000$ and $100$ repetitions of the procedures. The
left column deals with the comparison between Dantzig and Lasso,
the center column shows the effectiveness of our data driven
constraint and the right column illustrates the improvement of the
two-step method. As expected, Dantzig and Lasso estimators are
strictly equivalent when the dictionary is orthonormal and very
close otherwise. For both algorithms and most of the densities,
the best solution appears to be the ``Mix2'' dictionary, except
for the density $f_1$ where the Haar wavelets are better for
$n=500$. This shows that the dictionary approach yields an
improvement over the classical basis approach. One observes also
that the ``Mix'' dictionary is better than the best of its
constituent, namely the Fourier basis and the histogram family,
which corroborates our theoretical results. The adaptive
constraints are much tighter than their non adaptive counterparts
and yield to much better numerical results. Our last series of
experiments shows the significant improvement obtained with the
least square step. As hinted by Cand\`es and Tao \cite{candes},
this can be explained by the bias common to $\ell_1$ methods which
is partially removed by this final least square adjustment.
Studying directly the performance of this estimator is a
challenging task.

\section{Proofs}\label{proofs}
\subsection{Proof of Theorem \ref{concentrationTh}}\label{proofconcentrationTh}
To prove the first part of Theorem \ref{concentrationTh}, we fix
$m\in\{1,\dots,M\}$ and we set for any $i\in\{1,\dots,n\}$,
\[
W_i=\frac{1}{n}\left(\p_m(X_i)-\be_{0,m}\right)
\]
 that satisfies almost surely
\[
|W_i|\leq\frac{2\norm{\p_m}_\infty}{n}.
\]
Then, we apply Bernstein's Inequality (see \cite{mas} on pages
24 and 26) with the variables $W_i$ and $-W_i$: for any $u>0$,
\begin{equation}\label{conccoeff}
\P\left(|\hb_m-\be_{0,m}|\geq \sqrt{\frac{2\si_{0,m}^2u}{n}} +\frac{2u\norm{\p_m}_\infty}{3n} \right)\leq 2e^{-u}.
\end{equation}
Now, let us decompose $\hat\si^2_m$ in two terms:
\begin{align*}
\hat\si^2_m& =  \frac{1}{2n(n-1)}\sum_{i\not
  =j}(\p_m(X_i)-\p_m(X_j))^2\\
\begin{split}
&=\frac{1}{2n}\sum_{i=1}^n(\p_m(X_i)-\be_{0,m})^2+\frac{1}{2n}\sum_{j=1}^n(\p_m(X_j)-\be_{0,m})^2
\\
&\mspace{40mu}-\frac{2}{n(n-1)}\sum_{i=2}^n\sum_{j=1}^{i-1}(\p_m(X_i)-\be_{0,m})(\p_m(X_j)-\be_{0,m})
\end{split}
\\
&= s_n-\frac{2}{n(n-1)} u_n
\end{align*}
with
\begin{equation}\label{sn}
s_n=\frac{1}{n}\sum_{i=1}^n(\p_m(X_i)-\be_{0,m})^2 \mbox{ and }
u_n=
\sum_{i=2}^n\sum_{j=1}^{i-1}(\p_m(X_i)-\be_{0,m})(\p_m(X_j)-\be_{0,m}).
\end{equation}
Let us first focus on $s_n$ that is the main term of
$\hat\si^2_m$ by applying again  Bernstein's Inequality with
\[
Y_i=\frac{\si^2_{0,m}-(\p_m(X_i)-\be_{0,m})^2}{n}
\] which satisfies
\[
Y_i\leq\frac{\si^2_{0,m}}{n}.
\] One has that for any $u>0$
\[
\P\left(\si^2_{0,m}\geq s_n +\sqrt{2 v_m u}+\frac{\si^2_{0,m}
    u}{3n}\right)\leq e^{-u}
\]
with
\[
v_m=\frac{1}{n}\E\left(\left[\si^2_{0,m}-(\p_m(X_i)-\be_{0,m})^2\right]^2\right).
\]
But we have
\begin{align*}
v_m&= \frac{1}{n}\left(\si^4_{0,m}+\E\left[(\p_m(X_i)-\beta_{0,m})^4\right]-2\si^2_{0,m}\E\left[(\p_m(X_i)-\beta_{0,m})^2\right]\right)\\
&=\frac{1}{n}\left(\E\left[(\p_m(X_i)-\beta_{0,m})^4\right]-\si^4_{0,m}\right)\\
&\leq\frac{\si^2_{0,m}}{n}\left(\norm{\p_m}_{\infty}+|\beta_{0,m}|\right)^2\\
&\leq \frac{4\si^2_{0,m}}{n}\norm{\p_m}_\infty^2 .
\end{align*}
Finally, with for any $u>0$
\[
S(u)=2\sqrt{2}\si_{0,m}\norm{\p_m}_\infty\sqrt{\frac{u}{n}}+\frac{\si^2_{0,m}u}{3n},
\]
we have
\begin{equation}
\label{concsn}
\P(\si^2_{0,m}\geq s_n +S(u))\leq e^{-u}.
\end{equation}
The term  $u_n$ is a degenerate U-statistics that satisfies for any $u>0$
\begin{equation}
\label{concustats}
\P(|u_n|\geq U(u))\leq 6 e^{-u},
\end{equation}
with for any $u>0$
\[
U(u)=\frac{4}{3}Au^2+\left(4\sqrt{2}+\frac{2}{3}\right)Bu^{\frac{3}{2}}+\left(2D+\frac{2}{3}F\right) u+2\sqrt{2}C\sqrt{u},\]
where $A$, $B$, $C$, $D$ and $F$
are constants not depending on $u$ that satisfy
\begin{align*}
A&\leq 4\norm{\p_m}_{\infty}^2,\\
B&\leq 2\sqrt{n-1}\norm{\p_m}_{\infty}^2,\\
C&\leq\sqrt{\frac{n(n-1)}{2}}\si^2_{0,m},\\
D&\leq\sqrt{\frac{n(n-1)}{2}}\si^2_{0,m},\\
\intertext{and}
F&\leq 2\sqrt{2}\norm{\p_m}_{\infty}^2\sqrt{(n-1)\log(2n)}
\end{align*}
(see \cite{reyriv2}). Then, we have for any $u>0$,
\begin{align*}
\frac{2}{n(n-1)}U(u)&\leq\frac{32}{3}\frac{\norm{\p_m}_{\infty}^2}{n(n-1)}u^2+\left(16\sqrt{2}+\frac{8}{3}\right)\frac{\norm{\p_m}_{\infty}^2}{n\sqrt{n-1}}u^{\frac{3}{2}}\\
&\mspace{40mu}+\left(2\sqrt{2}\frac{\si^2_{0,m}}{\sqrt{n(n-1)}}+\frac{8\sqrt{2}}{3}\frac{\sqrt{\log(2n)}\norm{\p_m}_{\infty}^2}{n\sqrt{n-1}}\right)u+\frac{4\si^2_{0,m}}{\sqrt{n(n-1)}}\sqrt{u}.
\end{align*}
Now, we take $u$ that satisfies
\begin{equation}\label{condu1}
u=o(n)
\end{equation}
and
\begin{equation}\label{condu2}
\sqrt{\log (2n)}\leq\sqrt{2u}.
\end{equation}
 Therefore, for any $\e_1>0$, we have for $n$ large enough,
\begin{equation}
\frac{2}{n(n-1)}U(u)\leq\e_1 \si^2_{0,m}+\left(16\sqrt{2}+8\right)\frac{\norm{\p_m}_{\infty}^2}{n\sqrt{n-1}}u^{\frac{3}{2}}+\frac{32}{3}\frac{\norm{\p_m}_{\infty}^2}{n(n-1)}u^2.\nonumber
\end{equation}
So, for $n$ large enough,
\begin{equation}\label{majoU}
\frac{2}{n(n-1)}U(u)
\leq\e_1 \si^2_{0,m}+C_1\norm{\p_m}_{\infty}^2\left(\frac{u}{n}\right)^{\frac{3}{2}},
\end{equation}
where $C_1=16\sqrt{2}+19$. Using Inequalities (\ref{concsn}) and (\ref{concustats}), we obtain
\begin{align*}
\P\left(\si^2_{0,m}\geq \hat\si^2_m+S(u)+\frac{2}{n(n-1)}U(u)\right)
&=\P\left(\si^2_{0,m}\geq s_n-\frac{2}{n(n-1)}u_n+S(u)+\frac{2}{n(n-1)}U(u)\right)\\
&\leq \P\left(\si^2_{0,m}\geq s_n+S(u)\right)+\P\left(u_n\geq U(u)\right)\\
&\leq  7e^{-u}.
\end{align*}
Now, using (\ref{majoU}), for any $0<\e_2<1$, we have for $n$ large enough,
\begin{align*}
\hat\si^2_m+S(u)+\frac{2}{n(n-1)}U(u)&=\hat\si^2_m+2\sqrt{2}\si_{0,m}\norm{\p_m}_\infty\sqrt{\frac{u}{n}}+\frac{\si_{0,m}^2 u}{3n}+\frac{2}{n(n-1)}U(u)\\
&\leq \hat\si^2_m+2\sqrt{2}\si_{0,m}\norm{\p_m}_\infty\sqrt{\frac{u}{n}}+\frac{\si_{0,m}^2 u}{3n}+\e_1 \si_{0,m}^2+C_1\norm{\p_m}_{\infty}^2\left(\frac{u}{n}\right)^{\frac{3}{2}}\\
&\leq \hat\si^2_m+2\sqrt{2}\si_{0,m}\norm{\p_m}_\infty\sqrt{\frac{u}{n}}+\e_2 \si^2_{0,m}+C_1\norm{\p_m}_{\infty}^2\left(\frac{u}{n}\right)^{\frac{3}{2}}.
\end{align*}
Therefore,
\begin{equation}\label{e2}
\P\left((1-\e_2)\si^2_{0,m}\geq \hat\si^2_m+2\sqrt{2}\si_{0,m}\norm{\p_m}_\infty\sqrt{ \frac{u}{n}}+C_1\norm{\p_m}_{\infty}^2\left(\frac{u}{n}\right)^{\frac{3}{2}}\right)\leq 7e^{-u}.
\end{equation}
Now, let us set
\[
a=1-\e_2,\quad b=\sqrt{2}\norm{\p_m}_\infty\sqrt{ \frac{u}{n}},\quad
c=\hat\si^2_m+C_1\norm{\p_m}_{\infty}^2\left(\frac{u}{n}\right)^{\frac{3}{2}}
\]
and consider the polynomial
\[
P(x)=ax^2-2bx-c,
\]
with roots $\frac{b\pm\sqrt{b^2+ac}}{a}$. So, we have
\begin{align*}
P(\si_{0,m})\geq 0&\Longleftrightarrow\si_{0,m}\geq\frac{b+\sqrt{b^2+ac}}{a}
\\
&\Longleftrightarrow\si^2_{0,m}\geq\frac{c}{a}+\frac{2b^2}{a^2}+\frac{2b\sqrt{b^2+ac}}{a^2}.
\end{align*}
It yields
\begin{align*}
\P\left(\si^2_{0,m}\geq\frac{c}{a}+\frac{2b^2}{a^2}+\frac{2b\sqrt{b^2+ac}}{a^2}\right)&\leq 7e^{-u},
\intertext{so,}
\P\left(\si^2_{0,m}\geq\frac{c}{a}+\frac{4b^2}{a^2}+\frac{2b\sqrt{c}}{a\sqrt{a}}\right)&\leq 7e^{-u},
\end{align*}
which means that for any $0<\e_3<1$, we have for $n$ large enough,
\begin{align*}
\P\left(\si^2_{0,m}\geq(1+\e_3)\left(\hat\si^2_m+C_1\norm{\p_m}_{\infty}^2\left(\frac{u}{n}\right)^{\frac{3}{2}}+8\norm{\p_m}_\infty^2 \frac{u}{n}+2\sqrt{2}\norm{\p_m}_\infty\sqrt{ \frac{u}{n}}\sqrt{\hat\si^2_m+C_1\norm{\p_m}_{\infty}^2\left(\frac{u}{n}\right)^{\frac{3}{2}}}\right)\right)&\leq 7e^{-u}.
\end{align*}
Finally, we can claim that for any $0<\e_4<1$, we have for $n$ large enough,
\begin{align*}
\P\left(\si^2_{0,m}\geq(1+\e_4)\left(\hat\si^2_m+8\norm{\p_m}_\infty^2 \frac{u}{n}+2\norm{\p_m}_\infty\sqrt{2\hat\si^2_m\frac{u}{n}}\right)\right)&\leq 7e^{-u}.
\end{align*}
Now, we take $u=\ga\log M$. Under Assumptions of Theorem \ref{concentrationTh}, Conditions (\ref{condu1}) and (\ref{condu2}) are satisfied.
The previous concentration inequality means that
\begin{align*}
\P\left(\si^2_{0,m}\geq(1+\e_4)\tilde\si^2_m\right)&\leq 7M^{-\ga}.
\end{align*}
Now, using (\ref{conccoeff}), we have for $n$ large enough,
\begin{align*}
\begin{split}
\P\left(|\be_{0,m}-\hb_m|\geq\eta_{\ga,m}\right)&=\P\left(|\be_{0,m}-\hb_m|\geq\sqrt{\frac{2\tilde\si^2_m\ga\log
M}{n}}+\frac{2\norm{\p_m}_\infty\ga\log M}{3n},\
\si^2_{0,m}<(1+\e_4)\tilde\si^2_m\right)
\\
&\mspace{40mu}+
\P\left(|\be_{0,m}-\hb_m|\geq\eta_{\ga,m},\
  \si^2_{0,m}\geq(1+\e_4)\tilde\si^2_m\right)
\end{split}
\\
\begin{split}
&\leq\P\left(|\be_{0,m}-\hb_m|\geq\sqrt{\frac{2\si^2_{0,m}\ga(1+\e_4)^{-1}\log
M}{n}}+\frac{2\norm{\p_m}_\infty\ga(1+\e_4)^{-1}\log M}{3n}\right)\\
&\mspace{40mu}+
\P\left(\si^2_{0,m}\geq(1+\e_4)\tilde\si^2\right)
\end{split}
\\
&\leq2M^{-\ga(1+\e_4)^{-1}}+7M^{-\ga}.
\end{align*}
Then, the first part of Theorem \ref{concentrationTh} is proved: for  any $\e>0$,
\[
\P\left(|\be_{0,m}-\hb_m|\geq\eta_{\ga,m}\right)\leq
C(\e,\delta,\gamma)M^{-\frac{\ga}{1+\e}},
\]
where $C(\e,\delta,\gamma)$ is a constant that depends on $\e,$ $\delta$ and $\gamma$.

For the second part of the result, we apply again  Bernstein's Inequality with  \[Z_i=\frac{(\p_m(X_i)-\be_{0,m})^2-\si^2_{0,m}}{n}\] which satisfies
\[
Z_i\leq\frac{(\p_m(X_i)-\beta_{0,m})^2}{n}\leq
\frac{4\norm{\p_m}_\infty^2}{n}.
\]
One has that for any $u>0$
\[
\P\left(s_n\geq\si^2_{0,m}  +\sqrt{2 v_m
    u}+\frac{4\norm{\p_m}_\infty^2 u}{3n}\right)\leq e^{-u}
\]
with
\[
v_m=\frac{1}{n}\E\left(\left[\si^2_{0,m}-(\p_m(X_i)-\be_{0,m})^2\right]^2\right)\leq\frac{4\si^2_{0,m}}{n}\norm{\p_m}_\infty^2.
\]
So, for any $u>0$,
\[
\P\left(s_n\geq\si^2_{0,m}
  +2\sqrt{2}\si_{0,m}\norm{\p_m}_\infty\sqrt{\frac{
      u}{n}}+\frac{4\norm{\p_m}_\infty^2 u}{3n}\right)\leq e^{-u}.
\]
Now, for any $\e_5>0$, for any $u>0$,
\[
\P\left(s_n\geq(1+\e_5)\si^2_{0,m} +\frac{\norm{\p_m}_\infty^2
    u}{n}\left(\frac{4}{3}+\frac{2}{\e_5}\right)\right)\leq e^{-u}.
\]
Using (\ref{concustats}), with
\[
\tilde S(u)=\frac{\norm{\p_m}_\infty^2
  u}{n}\left(\frac{4}{3}+\frac{2}{\e_5}\right),
\]
\begin{align*}
\P\left(\hat\si^2_m\geq(1+\e_5)\si^2_{0,m} +\tilde S(u)+
\frac{2}{n(n-1)}U(u)\right)&=\P\left(s_n-\frac{2}{n(n-1)}u_n\geq(1+\e_5)\si^2_{0,m} +\tilde S(u)+
\frac{2}{n(n-1)}U(u)\right)\\
&\leq \P\left(s_n\geq(1+\e_5)\si^2_{0,m} +\tilde S(u)\right)+\P\left(-u_n\geq U(u)\right)\\
&\leq e^{-u}+6e^{-u}=7e^{-u}.
\end{align*}
Using (\ref{majoU}),
\begin{align*}
\P\left(\hat\si^2_m\geq(1+\e_1+\e_5)\si^2_{0,m} +\tilde
  S(u)+C_1\norm{\p_m}_{\infty}^2\left(\frac{u}{n}\right)^{\frac{3}{2}}\right)&\leq 7e^{-u}.
\end{align*}
Since
\[
\eta_{\ga,m}=\sqrt{\frac{2\tilde\si^2_m\ga\log
M}{n}}+\frac{2\norm{\p_m}_\infty\ga\log M}{3n},
\]
with
\[
\tilde\si^2_m=\hat\si^2_m+2\norm{\p_m}_\infty\sqrt{\frac{2\hat\si^2_m\ga\log M}{n}}+\frac{8\norm{\p_m}_\infty^2\ga\log
M}{n},
\]
we have for any $\e_6>0$,
\begin{align*}\label{majoseuil}
\eta_{\ga,m}^2&\leq (1+\e_6)\left(\frac{2\tilde\si^2_m\ga\log M}{n}\right)+(1+\e_6^{-1})\left(\frac{4\norm{\p_m}_\infty^2(\ga\log
M)^2}{9n^2}\right)\nonumber\\
\begin{split}
&\leq(1+\e_6)\left(\frac{2\ga\log M}{n}\right)\left(\hat\si^2_m+2\norm{\p_m}_\infty\sqrt{\frac{2\hat\si^2_m\ga\log M}{n}}+\frac{8\norm{\p_m}_\infty^2\ga\log
M}{n}\right)
\\&\mspace{40mu}+\frac{4}{9}(1+\e_6^{-1})\left(\frac{\norm{\p_m}_\infty\ga\log
M}{n}\right)^2
\end{split}
\nonumber\\
\begin{split}
&\leq(1+\e_6)^2\hat\si^2_m\left(\frac{2\ga\log M}{n}\right)+4\e_6^{-1}(1+\e_6)\left(\frac{\norm{\p_m}_\infty\ga\log
M}{n}\right)^2
\\
&\mspace{40mu}+16(1+\e_6)\left(\frac{\norm{\p_m}_\infty\ga\log
M}{n}\right)^2+\frac{4(1+\e_6^{-1})}{9}\left(\frac{\norm{\p_m}_\infty\ga\log
M}{n}\right)^2.
\end{split}
\end{align*}
Finally, with $u=\ga\log M$, with probability larger than $1-7M^{-\ga}$,
\[
\hat\si^2_m<(1+\e_1+\e_5)\si^2_{0,m}                                  +\tilde
S(\ga\log M)+C_1\norm{\p_m}_{\infty}^2\left(\frac{\ga\log
    M}{n}\right)^{\frac{3}{2}},
\]
and
\begin{align*}
\begin{split}
\eta_{\ga,m}^2&<(1+\e_6)^2(1+\e_5+\e_1)\si^2_{0,m}\left(\frac{2\ga\log
M}{n}\right)+(1+\e_6)^2\left(\frac{\ga\log
M}{n}\right)^2\norm{\p_m}_\infty^2\left(\frac{8}{3}+\frac{4}{\e_5}\right)\\
&\mspace{40mu}+2C_1(1+\e_6)^2\norm{\p_m}_{\infty}^2\left(\frac{\ga\log
M}{n}\right)^{\frac{5}{2}}+\norm{\p_m}_{\infty}^2\left(\frac{\ga\log
M}{n}\right)^{2}\left(4\e_6^{-1}(1+\e_6)+16(1+\e_6)+\frac{4(1+\e_6^{-1})}{9}\right).
\end{split}
\end{align*}
%
Finally, with $\e_6=1$, $\e_1=\e_5=\frac{1}{2}$, for $n$ large enough,
\[
\P\left(\eta_{\ga,m} \geq
4\si_{0,m}\sqrt{\frac{\ga\log
M}{n}}+\frac{10\norm{\p_m}_\infty\ga\log M}{n}\right)\leq
7M^{-\ga}.
\]
Note that $\sqrt{32/3+32+8+32+8/9}=9.1409$.

For  the   last  part,   starting  from  (\ref{e2})   with  $u=\ga\log   M$  and
$\e_2=\frac{1}{7}$, we have  for $n$ large enough and  with probability larger
than $1-7M^{-\ga}$,
\begin{align*}
\frac{6}{7}\si^2_{0,m}&\leq \hat\si^2_m+2\sqrt{2}\si_{0,m}\norm{\p_m}_\infty\sqrt{
\frac{\ga\log               M}{n}}+C_1\norm{\p_m}_{\infty}^2\left(\frac{\ga\log
M}{n}\right)^{\frac{3}{2}}\\
&\leq \hat\si^2_m+\frac{2}{7}\si^2_{0,m}+7\norm{\p_m}_\infty^2\frac{\ga\log               M}{n}+C_1\norm{\p_m}_{\infty}^2\left(\frac{\ga\log
M}{n}\right)^{\frac{3}{2}}.
\end{align*}
So, for $n$ large enough,
\[
\frac{4}{7}\si^2_{0,m}\leq \hat\si^2_m+8\norm{\p_m}_\infty^2\frac{\ga\log M}{n}\leq\tilde\si^2_m
\]
and
\[
\eta_{\ga,m}>   \si_{0,m}\sqrt{\frac{8\ga\log
M}{7n}}+\frac{2\norm{\p_m}_\infty\ga\log M}{3n}.
\]
\subsection{Proof of Theorem \ref{oracleadmiloc}}
Let $\lambda=(\lambda_m)_{m=1,\ldots,M}$ and
set $\Delta=\lambda-\hat{\lambda}^D$.
We have
\begin{equation}\label{rela1}
\normD{f_\lambda-\fo}^2= \normD{\hat{f}^D-\fo}^2 +
\normD{f_\lambda-\hat{f}^D}^2+2\int(\hat{f}^D(x)-\fo(x))(f_\lambda(x)-\hat{f}^D(x))dx.
\end{equation}
We have $\normD{f_\lambda-\hat{f}^D}^2=\normD{f_\Delta}^2$.
Moreover, with probability at least
$1-C_1(\e,\delta,\gamma)M^{1-\frac{\gamma}{1+\e}}$, we have
\begin{align}
\left|\int(\hat{f}^D(x)-\fo(x))(f_\lambda(x)-\hat{f}^D(x))dx\right|=&
\left|\sum_{m=1}^M(\lambda_m-\hat{\lambda}_m^D)\left[(G\hat{\lambda}^D)_m-\beta_{0,m}\right]\right|\label{maj2}\\
\le & \normu{\Delta} 2\norminf{\eta_{\gamma}}\nonumber,
\end{align}
where the last line is a consequence  of the
definition of the Dantzig estimator and of Theorem \ref{concentrationTh}. Then, we have
\begin{equation*}
\normD{\hat{f}^D-\fo}^2
\leq
 \normD{f_\lambda-\fo}^2+4\norminf{\eta_{\gamma}}\normu{\Delta}-\normD{f_\Delta}^2.
\end{equation*}
We use then the following Lemma:
 \begin{Lemma}\label{dantcon}
Let $J\subset\{1,\ldots,M\}$. For any $\lambda\in\R^M$
\begin{equation*}
\normu{\Delta_{\Comp{J}}} \leq
\normu{\Delta_{J}}+2\normu{\lambda_{\Comp{J}}}+ \left(
\normu{\hat\la^D} -  \normu{\la} \right)_+ ,\end{equation*} where
$\Delta=\hat \lambda^D-\lambda$.
\end{Lemma}
\begin{proof}[Proof of Lemma~\ref{dantcon}] This lemma is based on the fact that
\[
\normu{\hat\la^D}\leq \normu{\la} + \left( \normu{\hat\la^D} -  \normu{\la} \right)_+,
\]
which implies that
\[
\normu{\Delta_{J}+\la_{J}}+\normu{\Delta_{\Comp{J}}+\la_{\Comp{J}}}\leq\normu{\la_{J}}+\normu{\la_{\Comp{J}}}
+ \left( \normu{\hat\la^D} -  \normu{\la} \right)_+,
\]
and thus
\[
\normu{\la_{J}}-\normu{\Delta_{J}}+\normu{\Delta_{\Comp{J}}}-\normu{\la_{\Comp{J}}}\leq\normu{\la_{J}}+\normu{\la_{\Comp{J}}}
+ \left( \normu{\hat\la^D} -  \normu{\la} \right)_+.
\]
\end{proof}
Note that if $\lambda$ satisfies the Dantzig condition then by
definition of $\hat\la^D$: $ \left( \normu{\hat\la^D} -
\normu{\la} \right)_+=0$.
Using the previous lemma, we have:
\begin{equation*}\left( \normu{\Delta_{\JoC}} - \normu{\Delta_{\Jo}} \right)_+
\leq
2\normu{\lambda_{\JoC}}
+ \left( \normu{\hat\la^D} -  \normu{\la} \right)_+
.\end{equation*}
Using now
$\Lam=\normu{\lambda_{\JoC}}+\frac{\left(
    \normu{\hat\la^D} -  \normu{\la} \right)_+}{2}$, so that
$\Lam=\normu{\lambda_{\JoC}}$ as
soon as $\lambda$ satisfies the Dantzig condition, we obtain
\begin{align*}
 \normD{f_\Delta} &\geq \kappa_{\Jo}\normd{\Delta_{\Jo}} - \mu_{\Jo}
 \left( \normu{\Delta_{\Comp{\Jo}}}-\normu{\Delta_{\Jo}}\right)_+\\
& \geq \kappa_{\Jo}\normd{\Delta_{\Jo}} - 2 \mu_{\Jo
}\Lam
\end{align*}
and thus
\begin{align*}
\normd{\Delta_{\Jo}} \leq \frac{1}{\kappa_{\Jo}}  \normD{f_\Delta} + 2
\frac{\mu_{\Jo}}{\kappa_{\Jo}} \Lam.
\end{align*}
We deduce thus
\begin{align*}\normu{\Delta}&\le 2\normu{\Delta_{\Jo}}+2\Lam\\
 &\le 2\sqrt{|\Jo|}\normd{\Delta_{\Jo}}+2\widetilde{\normu{\lambda_{\JoC}}}\\
 &\le\frac{2\sqrt{|\Jo|}}{\kappa_{\Jo}}\normD{f_\Delta}+2\Lam\left(1+\frac{2\mu_{\Jo}\sqrt{|\Jo|}}{\kappa_{\Jo}}\right)\end{align*}
and then since
\[4
\norminf{\eta_{\gamma}}
\frac{2\sqrt{|\Jo|}}{\kappa_{\Jo}}\normD{f_\Delta}\leq
\frac{16|\Jo| \norminf{\eta_{\gamma}} ^2 }{\kappa_{\Jo}^2}+
\normD{f_\Delta}^2
\]
we have
\begin{align*}
4\norminf{\eta_{\gamma}} \normu{\Delta}-\normD{f_\Delta}^2 &\le
\frac{16|\Jo|\norminf{\eta_{\gamma}}^2}{\kappa_{\Jo}^2}+8
\norminf{\eta_{\gamma}}
\Lam\left(1+\frac{2\mu_{\Jo}\sqrt{|\Jo|}
}{\kappa_{\Jo}}\right)\\
& \leq 16|\Jo|
\left(\frac{1}{\beta}+\frac{1}{\kappa_{\Jo}^2}\right)
\norminf{\eta_{\gamma}}^2
+\beta\frac{\Lam^2}{|\Jo|}\left(1+\frac{2\mu_{\Jo}\sqrt{|\Jo|}}
{\kappa_{\Jo}}\right)^2,
\end{align*}
which is the result of the theorem.
\subsection{Consequences of Assumptions~1 and 2}\label{localcons}
To prove Proposition \ref{lema0}, we establish Lemmas \ref{lema3bis} and \ref{lema4}. In the sequel, we consider two integers $s$ and $l$ such that $1\le s\le M/2$, $l\ge s$ and $s+l\le M$.
We first recall Assumptions~1 and 2. Assumption~1 is stated in a more general form, which allows to unify the statement of the subsequent results.
\begin{itemize}\item\textbf{Assumption 1}
\[
\phi_{\min}(s+l)> \theta_{l,s+l}.
\]
\item\textbf{Assumption 2}
\[
l\phi_{\min}(s+l)> s\phi_{\max}(l).
\]
\end{itemize}
In the sequel, we assume that Assumptions~1 and 2 are both true.
\begin{Lemma}\label{lema3bis}
Let $\Jo\subset\{1,\ldots,M\}$  with cardinality $|\Jo|=s$ and $\Delta\in\R^M$.
 We denote by $\JUn$ the subset of $\{1,\ldots,M\}$ corresponding to
the $l$ largest  coordinates of  $\Delta$ (in absolute value) outside
$\Jo$ and we set $\JoUn=\Jo\cup  \JUn$.
We denote by $P_{\JoUn}$ the projector on the linear space spanned by
$(\varphi_{m})_{m\in \JoUn}$. We have:
\begin{align*}
\normD{P_{\JoUn}f_\Delta} \geq
\sqrt{\phi_{\min}(s+l)} \normd{\Delta_{\JoUn}}
-
\min\left(
\mu_1, \mu_2
\right)\normu{\Delta_{J_0^C}},
\end{align*}
with
\[\mu_1=\frac{\theta_{l,s+l}}{\sqrt{l\phi_{\min}(s+l)}}\quad\mbox{and}\quad\mu_2=\sqrt{\frac{\phi_{\max}(l)}{l}}.
\]
\end{Lemma}
\begin{proof}
For $k>1$, we denote by $J_k$ the indices corresponding to the coordinates
of $\Delta$ outside $J_0$ whose absolute values are between the
$((k-1)\times l+1)$--th and the $(k\times l)$--th largest ones (in absolute value). Note that this
definition is consistent with the definition of $J_1$. Using this notation, we have
\begin{align*}
  \normD{P_{\JoUn} f_\Delta}
& \geq \normD{P_{\JoUn} f_{\Delta_{\JoUn}}} - \normD{\sum_{k\geq 2} P_{\JoUn}
f_{\Delta_{J_k}}} \\
& \geq \normD{f_{\Delta_{\JoUn}}} - \sum_{k\geq 2} \normD{ P_{\JoUn}
f_{\Delta_{J_k}}}.
\end{align*}
Since $\JoUn$ has  $s+l$ elements, we have
\[
\normD{f_{\Delta_{\JoUn}}} \geq \sqrt{\phi_{\min}(s+l)}
\normd{\Delta_{\JoUn}}.
\]
Note that $P_{\JoUn} f_{\Delta_{J_k}} = f_{C_{\JoUn}}$ for some
vector $C\in\R^M$. Since,  $$\langle
P_{\JoUn} f_{\Delta_{J_k}} - f_{\Delta_{J_k}},  P_{\JoUn}
f_{\Delta_{J_k}}\rangle=  0,$$ one obtains that
\begin{align*}
\normD{P_{\JoUn} f_{\Delta_{J_k}}}^2
&= \langle f_{\Delta_{J_k}}, f_{C_{\JoUn}} \rangle
\intertext{and thus}
\normD{P_{\JoUn} f_{\Delta_{J_k}}}^2
& \leq \theta_{l,s+l} \normd{\Delta_{J_k}}
\normd{C_{\JoUn}}
 \leq \theta_{l,s+l} \normd{\Delta_{J_k}}
\frac{\normD{f_{C_{\JoUn}}}}{\sqrt{\phi_{\min}(s+l)}}\\
&\leq \frac{\theta_{l,s+l}}{\sqrt{\phi_{\min}(s+l)}}
\normd{\Delta_{J_k}} \normD{P_{\JoUn} f_{\Delta_{J_k}}}.
\intertext{This implies that} \normD{P_{\JoUn} f_{\Delta_{J_k}}}
&\leq \frac{\theta_{l,s+l}}{\sqrt{\phi_{\min}(s+l)}}
\normd{\Delta_{J_k}} = \mu_1 \sqrt{l} \normd{\Delta_{J_k}}.
\end{align*}
Moreover, using that
$J_k$ has less than $l$ elements, we obtain that
\[
\normD{P_{\JoUn} f_{\Delta_{J_k}}} \leq \normD{f_{\Delta_{J_k}}}
\leq \sqrt{\phi_{\max}(l)} \normd{\Delta_{J_k}} = \mu_2 \sqrt{l}
\normd{\Delta_{J_k}}.
\]
Now using that $\normd{\Delta_{J_{k+1}}} \leq
\normu{\Delta_{J_{k}}} / \sqrt{l}$, we obtain
\begin{align*}
  \sum_{k\geq 2} \normD{P_{\JoUn}f_{\Delta_{J_k}}}
\leq
\min\left(
\mu_1,\mu_2
\right) \normu{\Delta_{J_0^C}}
\end{align*}
and finally
\begin{align*}
\normD{P_{\JoUn}f_\Delta} \geq \sqrt{\phi_{\min}(s+l)}
\normd{\Delta_{\JoUn}} - \min\left( \mu_1,\mu_2 \right)
\normu{\Delta_{J_0^C}}.
\end{align*}
\end{proof}
\begin{Lemma}\label{lema4}
We use the same notations as in Lemma \ref{lema3bis}. For $c\geq 0$, assume that
\begin{equation}\label{majoration}
\normu{\Delta_{\JoC}}\le \normu{\Delta_{\Jo}}+ c.
\end{equation} Then we have
\begin{align*}
\normD{P_{\JoUn}f_\Delta} \geq
 \max
 \left(
\kappa_1,\kappa_2
\right)
 \normd{\Delta_{\JoUn}} -
\min\left(
\mu_1,\mu_2
\right)
c,
\end{align*}
with \[\kappa_1=\sqrt{\phi_{\min}(s+l)}\left(1-\frac{\theta_{l,s+l}}{\phi_{\min}(s+l)}\sqrt{\frac{s}{l}}\right)\quad\mbox{and}\quad\kappa_2=\sqrt{\phi_{\min}(s+l)}\left(1-\sqrt{\frac{s\phi_{\max}(l)}{l\phi_{\min}(s+l)}}\right).\]
\end{Lemma}
\begin{proof}
Using Lemma~\ref{lema3bis} and (\ref{majoration}), we obtain that
\begin{align*}
\normD{P_{\JoUn}f_\Delta} &\geq
\sqrt{\phi_{\min}(s+l)} \normd{\Delta_{\JoUn}}
-
\min\left(
\mu_1,\mu_2
\right) ( \normu{\Delta_{J_0}} +
c ).
\intertext{Using $\normu{\Delta_{J_0}} \leq \sqrt{s}
  \normd{\Delta_{J_0}}$, we deduce that}
\begin{split}
\normD{P_{\JoUn}f_\Delta} &\geq
\left( \sqrt{\phi_{\min}(s+l)} -
\sqrt{s}
\min\left(
\mu_1,\mu_2
\right)\right)
\normd{\Delta_{\JoUn}}
-c
\min\left(
\mu_1,\mu_2
\right)
\end{split}\\
\begin{split}
 &\geq
 \max
 \left(
\kappa_1,\kappa_2
\right)
 \normd{\Delta_{\JoUn}}  - c\min\left(
\mu_1,\mu_2
\right).
\end{split}
\end{align*}
\end{proof}
\subsection{Proof of Theorem \ref{lower}}
The dictionary considered here is the Haar dictionary
$(\phi_{jk})_{j,k}$ and is double-indexed. As a consequence, in
the following, the quantity $\beta_{0,jk}$, $\hat{\beta}_{jk}$,
$\si_{0,jk}^2$ $\eta_{\ga,jk}$, $\tilde{\si}^2_{jk}$ and
$\hat{\si}^2_{jk}$ are defined as in (\ref{betazero}),
(\ref{empirestim}), (\ref{sigmazero}), (\ref{etam}),
(\ref{sigmatilde}) and (\ref{sigmachapeau}) where $\varphi_m$ is
replaced by $\phi_{jk}$. Note that, since $f_0=\indic_{[0,1]}$, we
have, for $j\neq -1$, $\beta_{0,jk}=0$ and for any $j$, $\si_{0,jk}^2=1$ if $k\in\{0,\ldots,2^j-1\}$ and $0$ otherwise.

The proof of (\ref{majorisk}) is provided by using the oracle inequality satisfied by hard thresholding given by Theorem 1 of \cite{reyriv2} and the rough control of the soft thresholding estimate by the hard one:
\[\left||\hat{\beta}_{jk}|-\eta_{\ga,jk}\right|1_{\{|\hat{\beta}_{jk}|\geq\eta_{\ga,jk}\}}\leq 2
|\hat{\beta}_{jk}|1_{\{|\hat{\beta}_{jk}|\geq\eta_{\ga,jk}\}}.\]
An alternative is directly obtained by adapting the oracle results derived for soft thresholding rules in the regression model considered by Donoho and Johnstone \cite{dj:94}.

To prove (\ref{minorisk}), we establish the following lemma.
\begin{Lemma}\label{lem:calibration}
Let $\ga<1.$ We consider $j\in\mathbb{N}$ such that
\begin{equation}\label{condj}\frac{n}{(\log n)^\al}\leq 2^j<\frac{2n}{(\log
n)^\al}, \end{equation} for some $\al>1$. Then for all $\e>0$ such
that $\ga+2\e<1$,
\[\sum_{k=0}^{2^j-1}
 \E\left(\hb_{jk}^2\indic_{|\hb_{jk}|\geq\eta_{\ga,jk}}\right)\geq\frac{2\ga(1+\e)e^{-2}}{\pi}(\log
n)^{1-2\al}n^{-(\ga+2\e)}(1+o_n(1)).\]
\end{Lemma}
Then, we use the following inequality. For $j$ that
satisfies (\ref{condj}), we have for $r>0$,
\begin{eqnarray*}
 \E(\norm{\hat{f}^D-\fo}^2_2)&\geq &\sum_{k=0}^{2^j-1}
 \E\left(\left(|\hb_{jk}|-\eta_{\ga,jk}\right)^2\indic_{|\hb_{jk}|\geq\eta_{\ga,jk}}\right)\\
&\geq &\sum_{k=0}^{2^j-1}
 \E\left(\left(|\hb_{jk}|-\eta_{\ga,jk}\right)^2\indic_{|\hb_{jk}|\geq(1+r)\eta_{\ga,jk}}\right)\\
&\geq &\left(\frac{r}{r+1}\right)^2\sum_{k=0}^{2^j-1}
 \E\left(\hb_{jk}^2\indic_{|\hb_{jk}|\geq(1+r)\eta_{\ga,jk}}\right)\\
&\geq &\left(\frac{r}{r+1}\right)^2\sum_{k=0}^{2^j-1}
 \E\left(\hb_{jk}^2\indic_{|\hb_{jk}|\geq\eta_{jk,(1+r)^2\ga}}\right).
\end{eqnarray*}
So, if $r$ and $\e$ are such that $(1+r)^2\ga+2\e<1$, then applying Lemma
\ref{lem:calibration}, Inequality (\ref{minorisk}) is proved for any
$\delta$ such that $(1+r)^2\ga+2\e<\delta<1$.
\\

\begin{proof} [Proof of Lemma~\ref{lem:calibration}]
Let $j$ that satisfies (\ref{condj}) and $0\le k\le2^j-1$. We have
\[
\tilde{\si}^2_{jk}=\hat{\si}^2_{jk}+2\norm{\phi_{j,k}}_\infty\sqrt{2\gamma\hat{\si}^2_{jk}\frac{\ln
n}{n}}+8\gamma \norm{\phi_{j,k}}_\infty^2 \frac{\ln n}{n}.
\]
So, for any $0<\e<\frac{1-\gamma}{2}<\frac{1}{2}$,
\[
\tilde{\si}^2_{jk}\leq(1+\e)\hat{\si}^2_{jk}+2\gamma
\norm{\phi_{j,k}}_\infty^2 \frac{\ln n}{n}\left(\e^{-1}+4\right)
.\] Now,
\begin{eqnarray*}
\eta_{\gamma,jk}&=&\sqrt{2\gamma \tilde{\si}^2_{jk}\frac{\ln n}{n}} +\frac{2\norm{\phi_{j,k}}_\infty\gamma \ln n}{3n}\\
&\leq&\sqrt{2\gamma \frac{\ln n}{n}\left((1+\e)\hat{\si}^2_{jk}+2\gamma \norm{\phi_{j,k}}_\infty^2 \frac{\ln n}{n}\left(\e^{-1}+4\right)\right)}+\frac{2\norm{\phi_{j,k}}_\infty\gamma \ln n}{3n}\\
&\leq&\sqrt{2\gamma (1+\e)\hat{\si}^2_{jk} \frac{\ln
n}{n}}+\frac{2\norm{\phi_{j,k}}_\infty\gamma \ln
n}{n}\left(\frac{1}{3}+\sqrt{4+\e^{-1}}\right).
\end{eqnarray*}
Furthermore, we have
\[ \hat{\si}^2_{jk}=s_{njk}-\frac{2}{n(n-1)}u_{njk},\]
where $s_{njk}$ and $u_{njk}$ are defined as in (\ref{sn}) with
$\p_m$ replaced by $\phi_{jk}$. This implies that
\[\eta_{\gamma, jk}\leq\sqrt{2\gamma (1+\e)\frac{\ln n}{n}s_{njk}}+\sqrt{2\gamma (1+\e)\frac{\ln n}{n}\times\frac{2}{n(n-1)}|u_{njk}|}+\frac{2\norm{\phi_{j,k}}_\infty\gamma \ln n}{n}\left(\frac{1}{3}+\sqrt{4+\e^{-1}}\right).\]
Using (\ref{concustats}), with probability larger than
$1-6n^{-2}$, we have
\[|u_{njk}|\leq U(2\ln n),\]
and, since $\si_{0,jk}^2=1$
\begin{eqnarray*}
\frac{2}{n(n-1)}U(2\ln n)&\leq&\frac{c_1}{n}\sqrt{\log n}+\frac{c_2}{n}\log n+c_3\norm{\phi_{j,k}}_\infty^2\left(\frac{\log n}{n}\right)^{\frac{3}{2}}+c_4\norm{\phi_{j,k}}_\infty^2\left(\frac{\log n}{n}\right)^{2}\\
&\leq& C_1\frac{\log
n}{n}+C_2\norm{\phi_{j,k}}_\infty^2\left(\frac{\log
n}{n}\right)^{\frac{3}{2}},
\end{eqnarray*}
where $c_1$, $c_2$, $c_3$, $c_4$, $C_1$ and $C_2$ are universal
constants. Finally, with probability larger than $1-6n^{-2}$, we
obtain that
\[\sqrt{2\gamma (1+\e)\frac{\ln n}{n}\times\frac{2}{n(n-1)}|u_{njk}|}\leq\sqrt{2\gamma (1+\e)C_1}\frac{\ln n}{n}+\sqrt{2\gamma (1+\e)C_2}\norm{\phi_{j,k}}_\infty\left(\frac{\ln n}{n}\right)^{\frac{5}{4}}.\]
So, since $\gamma<1$, there exists $w(\e)$, only depending on $\e$ such that with probability larger than $1-6n^{-2}$,
\begin{eqnarray*}
\eta_{\ga,jk}&\leq&\sqrt{2\gamma (1+\e)\frac{\ln
n}{n}s_{njk}}+w(\e)\norm{\phi_{jk}}_\infty \frac{\ln n}{n}.
\end{eqnarray*}
We set
\[\widetilde{\eta_{\ga,jk}}=\sqrt{2\ga(1+\e) s_{njk}\frac{\log n}{n}}+w(\e)\frac{2^{\frac{j}{2}}\log n}{n}\]
so $\eta_{\ga,jk}\leq \widetilde{\eta_{\ga,jk}}.$ Then, we have
\begin{eqnarray*}
s_{njk}&=&\frac{1}{n}\sum_{i=1}^n\left(\phi_{jk}(X_i)-\beta_{0,jk}\right)^2\\
&=&\frac{2^j}{n}\sum_{i=1}^n\left(\indic_{X_i\in
[k2^{-j},(k+0.5)2^{-j}[}-\indic_{X_i\in
[(k+0.5)2^{-j},(k+1)2^{-j}[}\right)^2\\
&=&\frac{2^j}{n}\left(N^+_{jk}+N^-_{jk}\right),
\end{eqnarray*}
with \[N^+_{jk}=\sum_{i=1}^n\indic_{X_i\in
[k2^{-j},(k+0.5)2^{-j}[},\quad N^-_{jk}=\sum_{i=1}^n\indic_{X_i\in
[(k+0.5)2^{-j},(k+1)2^{-j}[}.\]
We consider $j$ such that
\[\frac{n}{(\log n)^\al}\leq 2^j<\frac{2n}{(\log n)^\al},\quad \al>1.\]
In particular, we have \[\frac{(\log n)^\al}{2}< n2^{-j}\leq (\log
n)^\al.\] Now, we can write
\[\hb_{jk}=\frac{1}{n}\sum_{i=1}^n\phi_{jk}(X_i)= \frac{2^{\frac{j}{2}}}{n}(N^+_{jk}-N^-_{jk}),\]
that implies that
\begin{eqnarray*}
&&\hspace{-1cm} \sum_{k=0}^{2^j-1}
 \E\left(\hb_{jk}^2\indic_{|\hb_{jk}|\geq\eta_{\ga,jk}}\right)\\
&\geq        &        \sum_{k=0}^{2^j-1}
\E\left(\hb_{jk}^2\indic_{|\hb_{jk}|\geq\widetilde{\eta_{\ga,jk}}}\indic_{|u_{njk}|\leq U(2\ln n)}\right)\\
&\geq & \sum_{k=0}^{2^j-1}
\frac{2^j}{n^2}\E\left((N^+_{jk}-N^-_{jk})^2 \indic_{|\hb_{jk}|\geq
\sqrt{2\ga (1+\e)s_{njk}\frac{\log n}{n}}+w(\e)\frac{2^{j/2}\log n}{n}}\indic_{|u_{njk}|\leq U(2\ln n)}\right).\\
&\geq & \sum_{k=0}^{2^j-1}
\frac{2^j}{n^2}\E\left((N^+_{jk}-N^-_{jk})^2
\indic_{\frac{2^{\frac{j}{2}}}{n}|N^+_{jk}-N^-_{jk}|\geq               \sqrt{2\ga(1+\e)
\frac{2^j}{n}\left(N^+_{jk}+N^-_{jk}\right)\frac{\log
n}{n}}+w(\e)\frac{2^{j/2}\log n}{n}}\indic_{|u_{njk}|\leq U(2\ln n)}\right)\\
&\geq & \sum_{k=0}^{2^j-1}
\frac{2^j}{n^2}\E\left((N^+_{jk}-N^-_{jk})^2
\indic_{|N^+_{jk}-N^-_{jk}|\geq               \sqrt{2\ga(1+\e)
\left(N^+_{jk}+N^-_{jk}\right)\log
n}+w(\e)\log n}\indic_{|u_{njk}|\leq U(2\ln n)}\right)\\
&\geq & \frac{2^{2j}}{n^2}\E\left((N^+_{j1}-N^-_{j1})^2
\indic_{|N^+_{j1}-N^-_{j1}|\geq               \sqrt{2\ga(1+\e)
\left(N^+_{j1}+N^-_{j1}\right)\log n}+w(\e)\log
n}\indic_{|u_{njk}|\leq U(2\ln n)}\right).
\end{eqnarray*}
Now, we consider a bounded sequence $(w_n)_n$ such that for any
$n$, $w_n\geq w(\e)$ and such that $\frac{\sqrt{v_{nj}}}{2}$ is an
integer with
\[v_{nj}=\left(\sqrt{4\ga(1+\e)\tilde\mu_{nj}\log(n)  }+w_n\log(n)\right)^2\]
and $\tilde\mu_{nj}$ is the largest integer smaller or equal to $n2^{-j-1}$. We have
\[v_{nj}\sim 4\ga(1+\e)\tilde\mu_{nj}\log n\]
since
\[\frac{(\log n)^\al}{4}-1<  n2^{-j-1}-1 <\tilde\mu_{nj}\leq n2^{-j-1}\leq \frac{(\log n)^\al}{2}.\]
Now, set
\[l_{nj}=\tilde\mu_{nj}+\frac{1}{2}\sqrt{v_{nj}},\quad m_{nj}=\tilde\mu_{nj}-\frac{1}{2}\sqrt{v_{nj}},\]
that are positive for $n$ large enough. If $N^+_{j1}=l_{nj}$ and
$N^-_{j1}=m_{nj}$ then we have $N^+_{j1}-N^-_{j1}=\sqrt{v_{nj}}$.
Finally, we obtain that
\begin{align}
&\sum_{k=0}^{2^j-1}
 \E\left(\hb_{jk}^2\indic_{|\hb_{jk}|\geq\eta_{\ga,jk}}\right)\nonumber\\
&\geq \frac{2^{2j}}{n^2}v_{nj}\P\left( N^+_{j1}=l_{nj},\quad
N^-_{j1}=m_{nj},\quad |u_{njk}|\leq U(2\ln n)\right)\nonumber\\
&\geq v_{nj}(\log n)^{-2\al}\left[\P\left( N^+_{j1}=l_{nj},\quad
N^-_{j1}=m_{nj}\right)-\P\left(|u_{njk}|>U(2\ln n)\right)\right]\nonumber\\
&\geq v_{nj}(\log
n)^{-2\al}\left[\frac{n!}{l_{nj}!m_{nj}!(n-l_{nj}-m_{nj})!}p_j^{l_{nj}+m_{nj}}(1-2p_j)^{n-(l_{nj}+m_{nj})}-\frac{6}{n^2}\right]\nonumber\\
&\geq v_{nj}(\log
n)^{-2\al}\times\left[\frac{n!}{l_{nj}!m_{nj}!(n-2\tilde\mu_{nj})!}p_j^{2\tilde\mu_{nj}}(1-2p_j)^{n-2\tilde\mu_{nj}}-\frac{6}{n^2}\right],\label{terms}
\end{align}
where
\[p_j=\int\indic_{
[2^{-j},(1+0.5)2^{-j}[}(x) f_0(x)dx=\int \indic_{
[(1+0.5)2^{-j},2^{-j+1}[}(x) f_0(x)dx=2^{-j-1}.\] Now, let us
study each term of (\ref{terms}). We have
\begin{eqnarray*}
p_j^{2\tilde\mu_{nj}}&=&\exp\left(2\tilde\mu_{nj}\log(p_j)\right)\\
&=&\exp\left(2\tilde\mu_{nj}\log(2^{-j-1})\right),
\end{eqnarray*}
\begin{eqnarray*}
(1-2p_j)^{n-2\tilde\mu_{nj}}&=&\exp\left((n-2\tilde\mu_{nj})\log(1-2p_j)\right)\\
&=&\exp\left(-(n-2\tilde\mu_{nj})2^{-j}+o_n(1)\right)\\
&=&\exp\left(-n2^{-j}\right)(1+o_n(1)),
\end{eqnarray*}
and
\begin{eqnarray*}
(n-2\tilde\mu_{nj})^{n-2\tilde\mu_{nj}}&=&\exp\left(\left(n-2\tilde\mu_{nj}\right)\log\left(n-2\tilde\mu_{nj}\right)\right)\\
&=&\exp\left(\left(n-2\tilde\mu_{nj}\right)\left(\log n+\log\left(1-\frac{2\tilde\mu_{nj}}{n}\right)\right)\right)\\
&=&\exp\left(\left(n-2\tilde\mu_{nj}\right)\log n-\frac{2\tilde\mu_{nj}\left(n-2\tilde\mu_{nj}\right)}{n}\right)(1+o_n(1))\\
&=&\exp\left(n\log n-2\tilde\mu_{nj}-2\tilde\mu_{nj}\log
n\right)(1+o_n(1)).
\end{eqnarray*}
Then, using the Stirling relation, $n!=n^ne^{-n}\sqrt{2\pi
n}(1+o_n(1))$, we deduce that
\begin{eqnarray*}
\frac{n!}{(n-2\tilde\mu_{nj})!}p_j^{2\tilde\mu_{nj}}(1-2p_j)^{n-2\tilde\mu_{nj}}
&=&\frac{e^{n-2\tilde\mu_{nj}}}{e^n}\times
\frac{n^n}{(n-2\tilde\mu_{nj})^{n-2\tilde\mu_{nj}}}\times  p_j^{2\tilde\mu_{nj}}(1-2p_j)^{n-2\tilde\mu_{nj}}\times
(1+o_n(1))\\
&=&\exp\left(-2\tilde\mu_{nj}\right)\times
\frac{\exp\left(n\log n\right)}{(n-2\tilde\mu_{nj})^{n-2\tilde\mu_{nj}}}\times  p_j^{2\tilde\mu_{nj}}(1-2p_j)^{n-2\tilde\mu_{nj}}\times
(1+o_n(1))\\
&=&\exp\left(-2\tilde\mu_{nj}\right)\times\frac{\exp\left(n\log n+2\tilde\mu_{nj}\log(2^{-j-1})-n2^{-j}\right)}{\exp\left(n\log n-2\tilde\mu_{nj}-2\tilde\mu_{nj}\log n\right)}(1+o_n(1))\\
&=&\exp\left(2\tilde\mu_{nj}\log n+2\tilde\mu_{nj}\log(2^{-j-1})-n2^{-j}\right)(1+o_n(1)).
\end{eqnarray*}
It remains to evaluate $l_{nj}!\times m_{nj}!$:
\begin{eqnarray*}
l_{nj}!\times      m_{nj}!&=&\left(\frac{l_{nj}}{e}\right)^{l_{nj}}\left(\frac{m_{nj}}{e}\right)^{m_{nj}}\sqrt{2\pi
l_{nj}}\sqrt{2\pi m_{nj}}(1+o_n(1))\\
&=&\exp\left(l_{nj}\log l_{nj}+m_{nj}\log
m_{nj}-2\tilde\mu_{nj}\right)\times 2\pi\tilde\mu_{nj}(1+o_n(1)).
\end{eqnarray*}
If we set
\[x_{nj}=\frac{\sqrt{v_{nj}}}{2\tilde\mu_{nj}}=o_n(1),\]
then
\[l_{nj}=\tilde\mu_{nj}+\frac{\sqrt{v_{nj}}}{2}=\tilde\mu_{nj}(1+x_{nj}),\]
\[m_{nj}=\tilde\mu_{nj}-\frac{\sqrt{v_{nj}}}{2}=\tilde\mu_{nj}(1-x_{nj}),\]
and using that
\begin{eqnarray*}
(1+x_{nj})\log(1+x_{nj})&=&(1+x_{nj})\left(x_{nj}-\frac{x_{nj}^2}{2}+\frac{x_{nj}^3}{3}+O(x_{nj}^4)\right)\\
&=&x_{nj}-\frac{x_{nj}^2}{2}+\frac{x_{nj}^3}{3}+x_{nj}^2-\frac{x_{nj}^3}{2}+O(x_{nj}^4)\\
&=&x_{nj}+\frac{x_{nj}^2}{2}-\frac{x_{nj}^3}{6}+O(x_{nj}^4),
\end{eqnarray*}
we obtain that
\begin{eqnarray*}
l_{nj}\log l_{nj}&=&\tilde\mu_{nj}(1+x_{nj})\log\left(\tilde\mu_{nj}(1+x_{nj})\right)\\
&=&\tilde\mu_{nj}(1+x_{nj})\log(1+x_{nj})+\tilde\mu_{nj}(1+x_{nj})\log\left(\tilde\mu_{nj}\right)\\
&=&\tilde\mu_{nj}\left(x_{nj}+\frac{x_{nj}^2}{2}-\frac{x_{nj}^3}{6}+O(x_{nj}^4)\right)+\tilde\mu_{nj}(1+x_{nj})\log\left(\tilde\mu_{nj}\right).
\end{eqnarray*}
Similarly, we obtain that
\begin{eqnarray*}
m_{nj}\log
m_{nj}&=&\tilde\mu_{nj}\left(-x_{nj}+\frac{x_{nj}^2}{2}+\frac{x_{nj}^3}{6}+O(x_{nj}^4)\right)+\tilde\mu_{nj}(1-x_{nj})\log\left(\tilde\mu_{nj}\right),
\end{eqnarray*}
that implies that
\begin{eqnarray*}
l_{nj}\log l_{nj}+m_{nj}\log m_{nj}&=&\tilde\mu_{nj}\left(x_{nj}^2+O(x_{nj}^4)\right)+2\tilde\mu_{nj}\log\left(\tilde\mu_{nj}\right)\\
&\leq&\tilde\mu_{nj}x_{nj}^2+2\tilde\mu_{nj}\log(n2^{-j-1})+O(\tilde\mu_{nj}x_{nj}^4).
\end{eqnarray*}
Since
\[\tilde\mu_{nj}x_{nj}^2=\frac{v_{nj}}{4\tilde\mu_{nj}}\sim\ga(1+\e)\log n,
\]
we have, for $n$ large enough,
\[
\tilde\mu_{nj}x_{nj}^2+O(\tilde\mu_{nj}x_{nj}^4)\leq(\ga+2\e)\log n
\]
and
\[l_{nj}\log l_{nj}+m_{nj}\log m_{nj}\leq (\ga+2\e)\log n+2\tilde\mu_{nj}\log(n2^{-j-1}).\]
Finally, we have
\begin{eqnarray*}
l_{nj}!\times      m_{nj}!
&=&\exp\left(l_{nj}\log l_{nj}+m_{nj}\log m_{nj}-2\tilde\mu_{nj}\right)\times 2\pi\tilde\mu_{nj}(1+o_n(1))\\
&\leq&\exp\left((\ga+2\e)\log
n+2\tilde\mu_{nj}\log(n2^{-j-1})-2\tilde\mu_{nj}\right)\times
2\pi\tilde\mu_{nj}(1+o_n(1)).
\end{eqnarray*}
Since $0<\e<\frac{1-\gamma}{2}<\frac{1}{2}$, we conclude that
there exists $\delta<1$ such that
\begin{align*}
& \sum_{k=0}^{2^j-1}
 \E\left(\hb_{jk}^2\indic_{|\hb_{jk}|\geq\eta_{\ga,jk}}\right)\\&\geq v_{nj}(\log n)^{-2\al}\left[
\frac{\exp\left(2\tilde\mu_{nj}\log
n+2\tilde\mu_{nj}\log(2^{-j-1})-n2^{-j}\right)}{\exp\left((\ga+2\e)\log
n+2\tilde\mu_{nj}\log(n2^{-j-1})-2\tilde\mu_{nj}\right)\times
2\pi\tilde\mu_{nj}}-\frac{6}{n^2}\right](1+o_n(1))
\\
&\geq \frac{v_{nj}(\log n)^{-2\al}}{2\pi\tilde\mu_{nj}}
\left[\exp\left(-(\ga+2\e)\log n-2\right)-\frac{6}{n^2}\right](1+o_n(1))\\
&\geq\frac{2\ga(1+\e)e^{-2}}{\pi}(\log
n)^{1-2\al}n^{-(\ga+2\e)}(1+o_n(1))
\end{align*}
 and Lemma~\ref{lem:calibration} is proved.
\end{proof}

\end{document}